\documentclass{gen-j-l}

\usepackage{amsfonts,amsmath,amssymb}
\usepackage{amsthm}
\usepackage{enumerate}
\usepackage{graphics}
\usepackage{graphicx}
\usepackage[cmtip,all]{xy}

\theoremstyle{plain} %
  \newtheorem{theorem}{Theorem}[section] %
  
  \newtheorem{problem}[theorem]{Problem}
  \newtheorem{basicproblem}[theorem]{Basic Problem}

\theoremstyle{definition}  
  \newtheorem{defin}[theorem]{Definition}%
  \newtheorem{example}[theorem]{Example}%
  \newtheorem{remark}[theorem]{Remark}%
  \newtheorem{observation}[theorem]{Observation}

\newcommand{\invHom}[3]{\operatorname{Hom}_{#1}({#2},{#3})}

\numberwithin{equation}{section}

\newcommand{\Exterior}{\mathchoice{{\textstyle\bigwedge}}%
    {{\bigwedge}}%
    {{\textstyle\wedge}}%
    {{\scriptstyle\wedge}}}    
    
\makeatletter
\@namedef{subjclassname@2020}{%
  \textup{2020} Mathematics Subject Classification}
\makeatother

\begin{document}

\title{Recent advances in branching problems of representations}

\author{Toshiyuki Kobayashi}
\address{
Graduate School of Mathematical Sciences,
The University of Tokyo, 3-8-1 Komaba, Meguro, Tokyo, 153-8914 Japan}
\email{toshi@ms.u-tokyo.ac.jp}

\address{
Translated by Toshihisa Kubo}
\email{toskubo@econ.ryukoku.ac.jp}




\begin{abstract}
How does an irreducible representation of a group 
behave when restricted to a subgroup? 
This is part of \emph{branching problems}, 
which are one of the fundamental problems in representation theory,
and also interact naturally with other fields of mathematics.

This expository paper is an up-to-date account on some new directions in
representation theory highlighting the branching problems for real 
reductive groups and their related topics ranging from global analysis
of manifolds via group actions to the theory of discontinuous groups
beyond the classical Riemannian setting.

This article is an outgrowth of the invited lecture that 
the author delivered at the commemorative event for 
the 70th anniversary of the re-establishment 
of the Mathematical Society of Japan,
and originally appeared in Japanese in Sugaku \textbf{71} (2019).
\end{abstract}

\maketitle

How does an irreducible 
representation of a group behave when restricted to a subgroup?
How is its restriction decomposed  (in a broad sense) into irreducible representations of the subgroup?
These questions are part of  \emph{branching problems}, 
which naturally emerge from various fields of mathematics.
The decomposition of the tensor product of two (or more) 
representations is such an example.
The Plancherel-type theorem that gives the expansion 
of functions on a homogeneous space $X$ is often equivalent to 
a special case of a branching problem via ``hidden symmetry''. 
The theta correspondence, which plays a prominent role in number theory, 
is also a branching law in the broad sense.
Further, when one tries to understand the geometry of a submanifold through
the function space on it, the branching problems 
for a pair of transformation groups
emerge in a natural way. 

For finite-dimensional representations, branching laws are in principle computable.
Combinatorial methods for computing branching laws exist,
 and various algorithms have been further developed.
On the other hand, irreducible
representations of non-compact reductive Lie groups such as 
$GL(n,\mathbb{R})$ are mostly infinite-dimensional and
a ``general algorithm'' that computes branching laws is still far from being known. Indeed, it often happens that
representations cannot be controlled well by a subgroup even when it is maximal
(Section \ref{subsec:6.2}).
Moreover, it could also be the case that the irreducible representations, 
which play as ``receptacles'', of a subgroup is 
not well-understood. Indeed, historically, \lq\lq{new}\rq\rq\ irreducible
representations were sometimes discovered through the branching
laws of the restriction of the \lq\lq{known}\rq\rq\ representations of a 
larger group.
Via the study of branching problems of a representation of a larger group,
we may expect to explore deeper properties of irreducible representations
of subgroups.

In the mid-1980s when the author started to challenge branching problems of 
infinite-dimensional representations,
there seemed to be a widespread pessimism that it was ``hopeless'' 
to build a general theory of the restriction of infinite-dimensional representations of reductive Lie groups except for some specific cases.
When one tries to obtain branching laws for groups larger than 
$SL(2,\mathbb{R})$ or $SL(2,\mathbb{C})$,
some  ``bad phenomena'' 
such as infinite multiplicities
appear and 
it was not easy to find a promising direction
to develop a theory of the branching problems.
In a thorough analysis of such ``bad phenomena'', 
the author encountered not a few ``mysterious  and nice phenomena''. 
From there he was fortunate to have been able to 
reveal some general principles for branching problems.
Through new themes such as
``the theory of discretely decomposable restrictions'',
``the theory of visible actions'',
``the theory of finite/bounded/multiplicity-free multiplicities'', 
\lq\lq{global analysis of the minimal representations}\rq\rq,
and
``the construction of symmetry breaking operators'',
which were born in this way,
more mathematicians have become increasingly interested in 
these new directions and are advancing the study of branching problems.
Now the \textbf{theory of branching problems for infinite-dimensional 
representations of reductive Lie groups}, 
which used to be regarded as hopeless, proceeds to a completely new developing stage.
	
It may be helpful to provide at this stage 
a brief overview of these developments
in the last 20--30 years.	
In this expository article we give the ideas on
the general theories of
``restrictions of representations'' 
and give some new perspectives with
programs to advance further the study of branching problems.
We collect some references in this area, most of which
are closely related to the viewpoint in this article.
The author apologizes to the many people whose works are 
not mentioned here because of the author's ignorance.
He would like to notice to the reader in advance  that 
there are many other important topics that are not treated here.

\medskip

\section{Branching Laws of Representations---Introduction}\label{sec:Intro}

As part of commemorative events for the 70th anniversary 
of the re-establishment of the Mathematical Society of Japan\footnote{
The organization was established in 1877 as the Tokyo Sugaku Kaisha
(Tokyo Mathematics Society).
This was the first academic society in Japan, and was
reorganized in its present form in 1946.},
the author was requested to give a plenary talk on 
an area of mathematics that are remarkably developed recently.
The lecture should be addressed to the wide audience 
about \lq\lq{why}\rq\rq, \lq\lq{what}\rq\rq, and \lq\lq{how}\rq\rq\
have been achieved in the past together with
future prospects. 
As it is not realistic to cover all the interesting topics,
the lecture 
was focused on revisiting 
the theory of branching problems of representations of 
reductive Lie groups, an area in which the author himself involved
and which has seen significant developments in the last 20--30 years.
Besides, recent progress of two relevant topics, 
global analysis via representation theory
and the theory of 
discontinuous groups, was also discussed.
This article is written based on the lecture notes \cite{xkmsj70}.
The author has attempted to clarify the ideas without technicalities
and to write this article in a manner that would simultaneously 
be generally comprehensive to graduate students and scholors
of diverse backgrounds and yet be of some value to experts
desiring to advance this field.

\subsection{What is a branching law of a representation?} \label{subsec:branch}

Hereafter, as a general rule,
irreducible representations of groups and those of their subgroups 
are denoted by uppercase and lowercase Greek letters, respectively, such as 
$\Pi$ and $\pi$.
When a representation $\Pi$ of a group $G$ is defined on a vector space $V$,
we write $\Pi|_{G'}$ for the representation of its subgroup $G'$, which is obtained by
restricting the action to the subgroup $G'$ on the same representation space $V$. 
Namely, $\Pi|_{G'}$ is the composition of two group homomorphisms
\begin{equation}
\label{eqn:1}
   \Pi|_{G'} \colon G' \hookrightarrow G \overset \Pi \longrightarrow GL(V).  
\end{equation}
Even if $\Pi$ is an irreducible representation of $G$, the restriction $\Pi|_{G'}$ 
is in general not irreducible as a representation of the subgroup $G'$.\footnote{
It may happen that the restriction $\Pi|_{G'}$ 
of an irreducible representation $\Pi$ of $G$ remains irreducible 
as a representation of a subgroup $G'$,
although such cases are rare.
See \cite{Zuckerman60} for a list of such triplets $(G, G', \Pi)$
and some geometric explanations.}

We begin with a classical case, where $\Pi$ is a finite-dimensional 
representation. If the restriction $\Pi|_{G'}$ is completely reducible, then
it can be described as the finite direct sum of irreducible representations of $G'$:
\begin{equation}
\label{eqn:fbranch}
  \Pi|_{G'} \simeq \bigoplus_{\pi}
            (\underbrace{\pi \oplus \cdots \oplus \pi}_{m(\Pi, \pi)})
            = \bigoplus_{\pi} m({\Pi}, \pi)\pi.
\end{equation}
The irreducible decomposition \eqref{eqn:fbranch} is 
a prototype of the \textbf{branching law}  of the restriction $\Pi|_{G'}$.
The number of times that the irreducible representation $\pi$ of $G'$ appears 
in $\Pi$ is called the \textbf{multiplicity} of the branching law.
Under the assumption that the finite-dimensional representation $\Pi|_{G'}$ is 
completely reducible, the following identities hold:
\begin{equation}
\label{eqn:fmult}
  m({\Pi}, \pi)=\dim_{\mathbb{C}}
               {\operatorname{Hom}}_{G'}(\pi, \Pi|_{G'})
              =\dim_{\mathbb{C}}
               {\operatorname{Hom}}_{G'}(\Pi|_{G'}, \pi).  
\end{equation}

In this article we deal with the case that $\Pi$ is an infinite-dimensional representation.
In contrast to the case that the representation space $V$ is finite-dimensional,
it may happen that there is no irreducible subrepresentation of 
the subgroup $G'$ in $V$. 
The ``multiplicity'' of the branching law may also be infinite even when
$G'$ is a maximal subgroup of $G$.
When talking about a \textbf{branching problem} 
for an infinite-dimensional representation,
we consider not only finding the branching law
(the irreducible decomposition of the restriction)
but also the following broader problem,
aiming to understand the ``restriction of the representation'' itself.

\begin{problem}[Branching Problem]\label{prob:1.1}
How does the restriction $\Pi|_{G'}$ behave as a representation of a subgroup $G'$?
\end{problem}

\subsection{Restriction of a representation and examples of branching problems}

Let us give a few examples on a ``restriction of a representation''
that arise from different contexts.
We would like the reader to feel the richness and diversity of
the themes related to branching problems
by browsing these examples and possibly
by skipping unfamiliar topics if any. 

\begin{example}[Tensor product representation]\label{ex:fusion}
The tensor product representations of two representations 
$(\pi_j, V_j)$ $(j=1,2)$ of a group $H$ 
\[
  \pi_1 \otimes \pi_2 \colon H \to GL(V_1 \otimes V_2), 
\quad
 h\mapsto \pi_1 (h) \otimes \pi_2 (h)
\]
can be interpreted as an example of the ``restriction of a representation''.
Namely, the tensor product representation $\pi_1 \otimes \pi_2$ 
is identified with the restriction of 
the outer tensor product representation
\[
  \pi_1 \boxtimes \pi_2 \colon H \times H \to GL(V_1 \otimes V_2), 
\quad
 (h_1, h_2) \mapsto \pi_1 (h_1) \otimes \pi_2 (h_2)
\]
of the direct product $G := H \times H$
to a subgroup 
$G'=\operatorname{diag}(H):=\{(h,h):h \in H \}$, which is isomorphic to $H$.
Fusion rules in theoretical physics are 
the irreducible decompositions 
of tensor product representations $\pi_1 \otimes \pi_2$.
The Clebsch--Gordan rule 
and the Pieri rule
are special cases of fusion rules.
\end{example}

\begin{example}[Cartan--Weyl's highest weight theory]
Cartan--Weyl's theory gives
a classification of 
irreducible finite-dimensional representations of a connected
compact Lie group $G$ by their \textbf{highest weights}.
One may interpret it as part of branching problems,
where the highest weights arise as the \lq\lq{edges}\rq\rq 
of the irreducible decompositions (branching laws) 
of the representations when restricted to a maximal torus $T$ of $G$,
or equivalently, as the unique subrepresentation of the restriction to 
a Borel subalgebra of the complexified Lie algebra $\mathfrak{g}_{\mathbb{C}}$
(Example \ref{ex:CWhw}).
\end{example}

\begin{example}[Vogan's minimal $K$-type theory]
Vogan's classification theory of 
irreducible admissible representations $\Pi$ of a reductive group $G$
(Definition \ref{def:adm}) 
utilizes {\bf{minimal $K$-type}} and 
\textbf{$\mathfrak{u}$-cohomology} \cite{Vogan81}.
This algebraic approach
is different from the analytic approach that was taken in
Langlands' classification
(Sections \ref{subsec:Irr} and \ref{ex:CWhw}).
One may interpret Vogan's method as a branching problem,
where minimal $K$-types are the \lq\lq{edges}\rq\rq of the branching laws 
of the representations restricted to a maximal compact subgroup $K$ of $G$
and $\mathfrak{u}$-cohomologies are 
defined by the restriction to a nilpotent Lie subalgebra 
$\mathfrak{u}$ of $\mathfrak{g}_{\mathbb{C}}$ 
as variants of highest weight vectors,
see also Example \ref{ex:coh}.
\end{example}

\begin{example}[Character theory]
Let $G$ be a reductive Lie group and let $H_1$, $\ldots$, $H_k$
be the complete system of representatives of its Cartan subgroups.
The distribution character ${\operatorname{Trace}}(\Pi)$ 
of irreducible admissible representation $\Pi$ (Definition \ref{def:adm}) is 
a locally integrable function on $G$, and thus,
it is determined by the restrictions to Cartan subgroups 
$H_j$ $(j=1,\ldots,k)$ (Harish-Chandra).
The study of the characters ${\operatorname{Trace}}(\Pi)|_{H_j}$ is 
related to the understanding of the restriction of the representation $\Pi$ 
to the subgroups $H_j$.
When $G$ is compact, one has $k=1$ and 
the explicit formula
${\operatorname{Trace}}(\Pi)|_{H_1}={\operatorname{Trace}}(\Pi|_{H_1})$
is known as the Weyl character formula.
\end{example}

\begin{example}[Theta correspondence]\label{ex:Howe}
Let $\Pi$ be the Weil representation of metaplectic group  $G=Sp^{\sim}$
and let $G':= G_1' \cdot G_2'$ be a subgroup of $G$ consisting of 
a dual pair, that is, $G_1'$ and $G_2'$ are the centralizers of each other in $G$.
The restriction $\Pi|_{G'}$ yields the theta correspondence 
between irreducible representations of $G_1'$ and $G_2'$
(Howe  \cite{xhowe}),
which may be also thought of as a branching law (in a broad sense)
from $G$ to the subgroup $G' = G_1' \cdot G_2'$.
\end{example}

\begin{example}[Rankin--Cohen differential operator]\label{ex:RC}
The Rankin--Cohen differential operator \cite{xcmz97, xrankin}, which
explicitly constructs modular forms of higher weights from ones of lower weights,
is a \lq\lq{symmetry breaking operator}\rq\rq\ for the decomposition of 
the tensor product representation of holomorphic discrete series representations
of $SL(2, {\mathbb{R}})$ (Section \ref{subsec:SL2}).
\end{example}

\begin{example}[Cohomology of a representation]\label{ex:coh}
When $\operatorname{Hom}_{G'}(\Pi|_{G'}, \pi)$ equals to zero,
it is natural to consider
higher order cohomologies $\operatorname{Ext}_{G'}^{\ast}(\Pi|_{G'}, \pi)$.
Especially, 
when $G$ is a reductive Lie group and
$\pi={\bf{1}}$ (trivial representation),
the cohomology for
a maximal nilpotent Lie subalgebra ${\mathfrak{n}}$ instead of a subgroup $G'$
contains some information related to 
the asymptotic behavior of matrix coefficients on analytic representation theory.
\end{example}

\begin{example}[Gross--Prasad conjecture]\label{ex:GP}
For a pair of orthogonal groups $(G,G')=(O_n, O_{n-1})$ defined over real or 
$p$-adic fields,
the restriction of an irreducible admissible representation
$\Pi$ of $G$ to a subgroup $G'$ is multiplicity-free (\cite{AGRS2010, xsunzhu}).
The Gross--Prasad conjecture \cite{GP} describes 
the branching law (in a broad sense)
for tempered representations $\Pi$ (Definition \ref{def:coarseLX}),
and is extended to the Gan--Gross--Prasad conjecture
\cite{GGPW12}
including a pair of groups $(G,G')=(GL_n, GL_{n-1})$.
\end{example}

\begin{example}[Modular variety]
On a locally Riemannian symmetric space $X= \Gamma \backslash G/K$
obtained by taking a quotient of arithmetic subgroup $\Gamma$,
a subgroup $G'$ of $G$ defines a cycle called a {\bf{modular variety}}.
Understanding modular variety is closely related to 
the branching law of the restriction of 
automorphic representations of $G$ to the subgroup $G'$,
as a dual notion between geometry and functions on it.
\end{example}

As observed in these examples,
the problem of the branching law is deeply related to the structure of representations
and it also arises widely in various fields of mathematics.

\subsection{Branching laws of infinite-dimensional representations}

In this article, with emphasis on the connection to global analysis, 
we consider infinite-dimensional representations of Lie groups.
Loosely speaking, there are two approaches in representation theory;
an algebraic approach which handles a representation without 
a topology on a representation space $V$
and an analytic approach which handles a (continuous) representation by 
equipping $V$ with a topology.
A representation $\Pi$ on a topological vector space $V$ of a topological group $G$
is called a {\bf{continuous representation}} if the map
\[
  G \times V \to V, 
  \qquad
  (g,v) \mapsto \Pi(g) v
\]
is continuous. Hereafter, we consider continuous representations
on a topological vector space over $\mathbb{C}$, 
unless otherwise specified.
\begin{defin}[Unitary representation]
Let $V$ be a Hilbert space.
A continuous representation $\Pi \colon G \to G L(V)$ of a group $G$ 
defined on $V$ is called a {\bf{unitary representation}} if 
$\Pi(g)$ is a unitary operator for all $g \in G$ on $V$.
\end{defin}

 The advantage of unitary representations is that the concept of 
 \lq\lq{irreducible decompositions}\rq\rq\ makes sense.
 That is, if $\Pi$ is a unitary representation of a locally compact topological
 group $G$, then the restriction $\Pi|_{G'}$ can be 
decomposed into a direct integral of irreducible representations of $G$,
which may be considered as a branching law in the unitary case
(see Theorem \ref{thm:4.1} below).
In contrast to the case of $\dim_\mathbb{C} V < \infty$,
the ``multiplicity'' of the branching law is not necessarily finite 
and further a ``continuous spectrum'' may also appear.
On the other hand, in a more general case in which $\Pi$ is not necessarily 
a unitary representation, e.g.\ $V$ is a Fr{\' e}chet space,
irreducible decompositions have a less clear meaning.
In such a case we may study 
continuous $G'$-homomorphisms between the irreducible representation
$\Pi$ of $G$
and an irreducible representation $\pi$ of its subgroup $G'$, instead of 
the ``irreducible decomposition'' of the restriction $\Pi|_{G'}$. Then
the following two concepts
\begin{alignat*}{2}
&\text{the space of symmetry breaking operators}\quad
&& \operatorname{Hom}_{G'}(\Pi|_{G'}, \pi), 
\\ 
&\text{the space of holographic operators}\quad
&& \operatorname{Hom}_{G'}(\pi, \Pi|_{G'})
\end{alignat*}
become important research objects.
The dimensions of these spaces could largely vary depending on
a choice of the topologies, e.g.\ 
the space $\Pi^\infty$ of smooth vectors or the space 
$\Pi^{-\infty}$ of distribution vectors,
on the representation spaces
\cite 
{xKVogan2015}.

\subsection{Branching laws of representations of reductive Lie groups}

In completing this introduction, we shortly highlight 
the main theme of this article, ``branching problems of reductive Lie groups'',
by comparing the cases of finite-dimensional representations 
with those of infinite-dimensional representations.

\vskip 0.1in

\begin{itemize}

\item
Branching laws of finite-dimensional representations
\vskip 0.05in

\begin{itemize}

\item 
Irreducible representations, 
which appear as building blocks of branching laws,
were classified  in the early 20th century
(Cartan--Weyl's highest weight theory, Example \ref{ex:CWhw}).
\vskip 0.05in

\item
There exist algorithms that compute branching laws,
and combinatorial techniques have been further developed,
e.g.\ the Littlewood--Richardson rule and Littelman's path method.

\end{itemize}

\vskip 0.1in

\item
Branching laws of infinite-dimensional representations (unitary case)
\vskip 0.05in

\begin{itemize}

\item
The classification of irreducible unitary representations, 
which appear as building blocks of branching laws
of unitary representations,
has a rich history of study but is not completely understood
(Section \ref{subsec:Irr}).
\vskip 0.05in

\item
Algorithms that compute branching laws are not known except for 
some special cases (e.g.\ theta-correspondence
or Howe correspondence \cite{xhowe} and
the restriction of highest weight modules \cite{mf-korea}).
\vskip 0.05in

\item 
A \lq\lq{continuous spectrum}\rq\rq\ may appear in branching laws
(Section \ref{subsec:irrdeco}).
\vskip 0.05in

\item
The \lq\lq{grip}\rq\rq\ of 
a subgroup $G'$ may not be strong enough,
involving a phenomenon of infinite multiplicities in branching laws
(Section \ref{subsec:6.2}).

\end{itemize}

\vskip 0.1in

\item
Branching laws of infinite-dimensional representations (when no unitarity is imposed)
\vskip 0.05in

\begin{itemize}

\item
The classification of 
irreducible admissible representations, 
which appear as building blocks of (non-unitary) branching laws,
was established in the 1970s to early 1980s
(Section \ref{subsec:Irr}).
\vskip 0.05in

\item
Nevertheless, 
it is generally a difficult problem to determine when symmetry breaking operators
exist (e.g. the Gan--Gross--Prasad conjecture
for specific pairs $(G,G')$, see Example \ref{ex:GP}). 
Recently, new theories not only for the existence but also
for the construction and classification of
symmetry breaking operators are emerging (Section \ref{subsec:10.3}).

\end{itemize}

\end{itemize}

\vskip 0.1in

In the following we shall introduce new programs
in the branching problem and some recent developments,
while explaining the basic notion and 
facts mentioned here.

\section{Search for Fundamental Objects and 
\lq\lq{More}\rq\rq\ 
Fundamental Objects---the Classification of \lq\lq{Irreducibles}\rq\rq\
and Decomposition to \lq\lq{Irreducibles}\rq\rq}
\label{sec:as}

Let us consider what the possible roles of the branching problems 
are in the representation theory of Lie groups.
The notion of Lie groups (continuous groups) was 
introduced in 1870s by Sophus Lie (1842--1899)
in his attempt to develop an analogous theory for differential equations to
the Galois theory for algebraic equations.
Lie groups and their representation theory have been developed through
numerous intersections with analysis, geometry, algebra, 
and theoretical physics.
In this section, applying the philosophical notion ``analysis and synthesis'':
\begin{enumerate}
\item understanding the ``smallest objects'' (e.g.\ classification)  and
\item how things are built up from the ``smallest objects'',
\end{enumerate}
we give a brief account of the current status of
what problems have been solved and what problems remain open 
in the representation theory of Lie groups.
We then briefly describe the connection of these 
with other branches of mathematics.
Part of Section \ref{sec:as} overlaps with the earlier
exposition \cite{K94b} of the author.

\subsection{Lie groups and Lie algebras}\label{subsec:2.1new}\hfill\\[3pt]
\noindent{\bf{\underline{Lie algebras and their smallest objects}}}:
The \lq\lq{smallest objects}\rq\rq\ of Lie algebras are 
one-dimensional (abelian) Lie algebras and simple Lie algebras
which are of dimension $\geq 2$
and do not have nontrivial ideals.
Finite-dimensional simple Lie algebras over ${\mathbb{R}}$ are
classified as the 10 series of classical Lie algebras, namely,
${\mathfrak{sl}}(n,{\mathbb{R}})$, 
${\mathfrak{sp}}(n, {\mathbb{R}})$, 
${\mathfrak{su}}(p,q)$,
${\mathfrak{su}}^{\ast}(2n)$, 
${\mathfrak{so}}^{\ast}(2n)$, 
${\mathfrak{so}}(p,q)$, 
${\mathfrak{sp}}(p,q)$, 
${\mathfrak{sl}}(n,{\mathbb{C}})$,
${\mathfrak{so}}(n,{\mathbb{C}})$, and
${\mathfrak{sp}}(n,{\mathbb{C}})$, 
and 22 exceptional ones ({\'E}.~Cartan, 1914).

\vskip 0.05in
\par\noindent{\bf{\underline{Building up from the smallest objects}}}:
Any finite-dimensional Lie algebra is obtained by 
iterating extensions of simple Lie algebras or abelian Lie algebras.
When the extension is trivial, that is,
when
the Lie algebras are expressed as the direct sum of simple Lie algebras and
abelian Lie algebras, they are called {\bf{reductive Lie algebras}}.
When the Lie algebras are expressed as the direct sum of only simple Lie algebras,
they are called {\bf{semisimple Lie algebras}}.

\vskip 0.05in
\par\noindent{\bf{\underline{Lie groups}}}:
A Lie group is a group that carries a manifold structure
with continuous (equivalently, smooth) multiplication.
Typical examples of Lie groups include {\bf{algebraic groups}},
which are groups obtained as the zero sets of polynomials on $M(n,{\mathbb{R}})$.
Lie algebras are the infinitesimal algebraic structures of Lie groups
and all the local properties and some of global ones of 
Lie groups can be described in terms of Lie algebras
({\bf{Lie theory}}). 
Lie groups whose corresponding Lie algebras are 
simple, semisimple, and reductive are called simple Lie groups,
semisimple Lie groups, and {\bf{reductive Lie groups}}, respectively.

\vskip 0.05in
\par\noindent
{\bf{\underline{Reductive Lie groups}}}:
Reductive Lie groups are locally isomorphic to the direct product 
of abelian Lie groups and simple Lie groups.
{\bf{Classical groups}} such as 
$G L(n,{\mathbb{R}})$, $GL(n,{\mathbb{C}})$, 
$O(p,q)$, $S p(n,{\mathbb{R}})$, 
 $\dots$ are reductive algebraic Lie groups, that is, 
 algebraic groups as well as reductive Lie groups.
Any two maximal compact subgroups of a Lie group $G$
are conjugate to each other by an inner automorphism.
Let $G$ be a semisimple Lie group of finite center, 
and $K$ a maximal compact subgroup of $G$. Then
the homogenous space $G/K$ is simply connected, and 
has the structure of 
a Riemannian symmetric space \cite[Chap.\ VI]{Helgason2}.
The classification of simple Lie algebras is equivalent to
that of simply connected irreducible Riemannian symmetric spaces.

\subsection{Fundamental problems in representation theory}
Since representations are defined on vector spaces that has linearity,
the superposition principle may be considered in an equivariant fashion 
for group representations.
Then the viewpoint of the
\lq\lq{building up from the smallest units}\rq\rq\
given  in the beginning of Section \ref{sec:as}
raises the following two fundamental problems in representation theory:
\begin{enumerate}
\item classify the irreducible representations;
\item decompose given representations irreducibly.
\end{enumerate}
Branching problems are one of the main themes related to the  latter, 
``irreducible decomposition''.
However, not only that, it provides a useful method on the former,
``classification of irreducible representations''.
We shall illustrate this idea with some examples in Section \ref{subsec:2.7}.

\subsection{Irreducible decomposition}\label{subsec:irrdeco}
In general representations cannot always be decomposed into
the direct sum of irreducible representations.
Indeed, 
it sometimes happens that a representation $\Pi$
does not contain any irreducible submodule when
$\Pi$ is infinite-dimensional.
In a program for the branching problems described in Section \ref{sec:Stage},
we shall propose a new direction of research, that is,
the study of 
\lq\lq{symmetry breaking operators}\rq\rq\ 
between representations of two groups $G \supset G'$,
which are not necessarily unitary and 
may not be decomposed into irreducibles in the usual sense.
However, in this section, we first focus on a simpler situation in which
irreducible decompositions make sense.
For this we recall some classical results.

The best setting for irreducible decompositions is 
the case when the representations are unitary defined on a Hilbert space.
The ``smallest units'' in this case are irreducible unitary representations.
The \lq\lq{decomposition}\rq\rq\ of a Hilbert space is defined
by using the notion of the direct integral of 
a family of Hilbert spaces introduced by von Neumann.
Then a result of Mautner and Teleman
(Theorem \ref{thm:4.1} given below)
states that any unitary representation $\Pi$ can be decomposed into
irreducibles 
(\cite{Mautner50}, also \cite[Chapter 6]{xtatsuuma}).

\begin{defin}[Unitary dual]
The set $\widehat{G}$ of equivalence classes 
of irreducible unitary representations of a topological group $G$ 
is called the {\bf{unitary dual}} of $G$.
For a locally compact group $G$, 
the unitary dual $\widehat G$ carries a natural
topological structure called the Fell topology. 
\end{defin}

\begin{theorem}
[Irreducible decomposition of a unitary representation by the direct integral]
\label{thm:4.1}
For any unitary representation $\Pi$ of a locally compact group $G$,
there exist a Borel measure $\mu$ on the unitary dual $\widehat{G}$
and a measurable function 
$n_{\Pi} \colon \widehat{G} \to {\mathbb{N}} \cup \{\infty\}$ 
such that $\Pi$ is unitarily equivalent to the direct integral
of irreducible unitary representations:
\begin{equation}
\label{eqn:directdeco}
   \Pi \simeq \int_{\widehat{G}}^{\oplus} n_{\Pi}(\pi) \pi d \mu (\pi).
\end{equation}
\end{theorem}

Here $n_{\Pi}(\pi)\pi$ stands for a multiple of $(\pi, V_{\pi})$ as in 
\eqref{eqn:fbranch}, namely, a unitary representation of $G$ on the Hilbert space
$\mathcal{H}_\pi\otimes V_\pi$ where $\mathcal{H}_\pi$ is a Hilbert space of dimension $n_{\Pi}(\pi)$ with trivial $G$-action.
We may think of \eqref{eqn:directdeco} as the irreducible decomposition
of the unitary representation $\Pi$.
Theorem \ref{thm:4.1} includes 
the case that a  \lq\lq{continuous spectrum}\rq\rq\ appears
in the irreducible decomposition.

It should be noted that
Theorem \ref{thm:4.1} assures the \lq\lq{existence}\rq\rq\ of an irreducible
decomposition, but the uniqueness may fail in general. However, the 
irreducible decomposition \eqref{eqn:directdeco} of a unitary 
representation is unique up to unitary equivalence if the group $G$ is 
an algebraic Lie group, for which the proof of 
the uniqueness of the irreducible decomposition \eqref{eqn:directdeco} 
is reduced to the case where
$G$ is a real reductive Lie group, and 
this case was proved by Harish-Chandra by using the $K$-admissibility 
theorem (see Remark \ref{rem:3.4} (1) below),
see \cite[Thm.\ 14.6.10]{WaI} for instance.
We shall work mainly with reductive algebraic Lie groups, hence
the irreducible decomposition \eqref{eqn:directdeco}
is unique in our setting throughout the paper.
In this case the function $n_{\Pi}\colon \widehat{G}\to \mathbb{N}\cup \{\infty\}$
is well-defined up to a measure zero set, and is referred to as the 
multiplicity for the unitary representation $\Pi$.

\begin{defin}\label{def:coarseLX}
\phantom{a}
\vskip 0.05in
\begin{enumerate}
\item[{\rm{(1)}}]
(Support of an irreducible decomposition)\enspace
Let ${\operatorname{Supp}}_{\widehat G}(\Pi)$ be
the subset of the unitary dual $\widehat G$ defined as the support of 
the irreducible decomposition \eqref{eqn:directdeco} of a unitary
representation $\Pi$. 
If $\operatorname{Supp}_{\widehat G}(\Pi)$ is a countable set,
then $\Pi$ is decomposed into the direct sum of irreducible unitary representations
($\Pi$ is said to be {\bf{discretely decomposable}}).
\vskip 0.05in
\item[{\rm{(2)}}]
(Tempered representation)\enspace
When $\Pi$ is the regular representation $L^2(G)$ of $G$,
the set ${\operatorname{Supp}}_{\widehat G}(L^2(G))$ 
is denoted by $\widehat {G}_{\operatorname{temp}}$.
A unitary representation $\Pi$ of $G$ is called
a {\bf{tempered representation}} if 
$\Pi$ satisfies
${\operatorname{Supp}}_{\widehat G}(\Pi)
 \subset \widehat {G}_{\operatorname{temp}}$,
 equivalently, if $\Pi$ is weakly contained in the regular representation $L^2(G)$
 in the sense that any matrix coefficient is approximated by a linear combination
 of matrix coefficients of $L^2(G)$ on compact sets.
\end{enumerate}
\end{defin}

\begin{example}\label{ex:abeltemp}
If a Lie group $G$ is amenable, in particular, if $G$ is solvable or compact,
then $\widehat G_{\operatorname{temp}} = \widehat G$.
\end{example}

\begin{example}\label{ex:KZtemp}
When $G$ is a reductive Lie group,
irreducible tempered representations can be 
characterized by the asymptotic behavior of their matrix coefficients.
The classification of the set $\widehat {G}_{\operatorname{temp}}$ 
was accomplished by Knapp--Zuckerman \cite{KZ} 
in terms of the direct products of ${\mathbb{Z}}/2{\mathbb{Z}}$
called $R$ groups.
\end{example}

\subsection{Classification problem of the unitary dual---the orbit method 
and geometric quantization of symplectic manifolds}
\label{subsec:orbit}

Do we know all irreducible unitary representations of Lie groups?
Actually, for a general Lie group $G$,
the classification of the unitary dual $\widehat G$ 
has not been completely understood.
Postponing a brief summary of the current status 
to Sections \ref{subsec:Duflo} and \ref{subsec:Irr},
we begin with some typical examples in which
the unitary dual is completely classified.

\begin{example}[Abelian group]
All irreducible unitary representations of an abelian group $G$
are one-dimensional.
For instance, when $G={\mathbb{R}}$,
set $\chi_{\xi} \colon {\mathbb{R}} \to G L (1,{\mathbb{C}})$, 
 $x \mapsto e^{i x \xi}$. Then we have
$
  \widehat {\mathbb{R}} \simeq \{\chi_{\xi} : \xi \in {\mathbb{R}}\}.  
$
\end{example}

\begin{example}[Compact Lie group, Cartan--Weyl 1925]
\label{ex:CWhw}
For a compact group $G$,
all irreducible unitary representations $\Pi$ of $G$ are finite-dimensional.
Suppose $G$ is a connected compact Lie group with Lie algebra 
${\mathfrak{g}}$  and 
${\mathfrak{b}}$ is a Borel subalgebra of the complexified Lie algebra
${\mathfrak{g}}_{\mathbb{C}}={\mathfrak{g}} \otimes_{\mathbb{R}}{\mathbb{C}}$.
Then irreducible representations $\Pi$ can be classified by the unique 
one-dimensional subrepresentations
$\chi$ ({\bf{highest weight}})  of ${\mathfrak{b}}$.
\end{example}

\begin{example}
[$SL(n, {\mathbb{F}})$, ${\mathbb{F}}={\mathbb{R}}, {\mathbb{C}}, {\mathbb{H}}$]
\label{ex:SLunitary}
The unitary dual of 
$G=PSL(2,{\mathbb{R}})$ consists of 
the trivial representation, 
uncountable family of
infinite-dimensional representations
that are not equivalent to each other
(spherical principal series representations, complementary series representations),
and countable family of infinite-dimensional representations
(discrete series representations and the limit of discrete series representations)
(Bargmann 1947).
This result is generalized to $G=S L(n, {\mathbb{R}})$
($n=3$: Vakhutinski 1968, 
 $n=4$: Speh 1981, and
$n$ general: Vogan 1986).
The unitary dual of $G=SL(n,{\mathbb{C}})$ started with
the case of $n=2$ by Gelfand--Naimark (1947),
Tsuchikawa 1968 ($n=3$),
Duflo 1980 ($n\le 5$), and
Barbasch 1985 ($n$ general).
The unitary dual of $G=SL(n,{\mathbb{H}})$ is classified by
Hirai 1962 ($n=2$) and
Vogan 1986 ($n$ general).
\end{example}

\begin{example}[Nilpotent Lie group, Kirillov {\cite{kirillov62}}, 1962]
\label{ex:Kirillov}
For a simply connected nilpotent Lie group $G$, 
Kirillov discovered that
there exists a natural bijection
\begin{equation}
    \widehat G
    \,\,\simeq \,\,
    {\mathfrak {g}}^{\ast}/\operatorname{Ad}^{\ast}(G),
\label{eqn:orbit}
\end{equation}
where 
$\operatorname{Ad}^{\ast} \colon G \to GL_{\mathbb{R}}({\mathfrak {g}}^{\ast})$
is the coadjoint representation of the Lie group $G$, namely,
the contragredient representation of 
${\operatorname{Ad}}\colon G \to GL_{\mathbb{R}}({\mathfrak {g}})$.
This result is extended to 
exponential solvable Lie groups
\cite{F2010},
see also \cite{FL15}.
\end{example}

The unitary dual $\widehat{G}$ for non-compact non-commutative group
is \lq\lq{huge}\rq\rq in a sense, as it may contain uncountably many equivalence
classes of irreducible infinite-dimensional representations.
Example \ref{ex:Kirillov} above  asserts 
that when $G$ is a nilpotent Lie group, 
the unitary dual $\widehat G$, though it looks  \lq\lq{huge}\rq\rq,
can be parametrized by just one finite-dimensional representation, namely,
the coadjoint representation $(\operatorname{Ad}^{\ast}, {\mathfrak{g}}^{\ast})$.
The right-hand side of \eqref{eqn:orbit} is nothing but the set of 
coadjoint orbits, and the idea to understand the unitary dual
$\widehat{G}$ via coadjoint orbits is referred to as
the {\bf{orbit method}} or the {\bf{orbit philosophy}} 
initiated by Kirillov, and developed by Kostant, Duflo and others, see \cite{kirillov}.
For instance, when $G=G L(n,{\mathbb{R}})$,
the set $\mathfrak{g}^*/\operatorname{Ad}^*(G)$ of coadjoint orbits
may be identified with
the set of the Jordan normal forms of real square matrices of order $n$,
which is \lq\lq{supposed to approximate}\rq\rq\
the parameter set of irreducible unitary representations 
of $GL(n, {\mathbb{R}})$ according to the orbit philosophy.
Let us see that this is not just a coincidence 
and that the \lq\lq{correspondence}\rq\rq\ from coadjoint orbits
to irreducible unitary representations may
be interpreted as a 
\lq\lq{geometric quantization}\rq\rq\ at least to some extent.

The
\lq\lq{quantization}\rq\rq\
in physics
is the procedure from a classical understanding of physical phenomena
to a newer understanding known as quantization theory such as
\begin{center}
``classical mechanics $\rightsquigarrow$ quantum mechanics''.
\end{center}
As its mathematical analogue,
one may wish to define
\lq\lq{{\bf{geometric quantization}}}\rq\rq\
under suitable assumptions:
\begin{alignat}{2}
\text{symplectic manifold $M$} 
&\rightsquigarrow 
&&\text{Hilbert space ${\mathcal{H}}$};
\notag
\\
\text{symplectic transformations on $M$} 
&\rightsquigarrow\,\,
&&\text{unitary operators on ${\mathcal{H}}$}.
\notag
\end{alignat}
One may further expect the naturality and the functoriality
of this correspondence, in particular, the following properties 
are supposed to hold for a family of transformations:
\begin{alignat}{2}
\text{Hamiltonian actions on $M$} 
&\rightsquigarrow\,\,
&&\text{unitary representations on ${\mathcal{H}}$;}\label{eqn:gq}\\
\text{transitivity} 
&\rightsquigarrow\,\,
&&\text{irreducitiblity.}
\notag
\end{alignat}

Now, the right-hand side of  \eqref{eqn:orbit}  on the orbit method
is identified with the space of coadjoint orbits 
 ${\mathcal{O}}_{\lambda}=\operatorname{Ad}^{\ast}(G) \lambda$
 ($\lambda \in {\mathfrak {g}}^{\ast}$) of the Lie group $G$.
Each coadjoint orbit ${\mathcal{O}}_{\lambda}$ has 
a symplectic structure induced by the skew-symmetric bilinear form
\[
  {\mathfrak {g}}\times {\mathfrak {g}} \to {\mathbb{R}},
\quad
  (X,Y) \mapsto \lambda([X,Y])
\]
and the group action of $G$ is clearly Hamiltonian and transitive
(Kirillov--Kostant--Soureau).
Therefore, if the idea of the geometric quantization \eqref{eqn:gq} 
works, then one expects to obtain irreducible unitary representations of $G$ 
corresponding to the coadjoint orbits ${\mathcal{O}}_{\lambda}$,
namely, the correspondence from
the right-hand side to the left-hand side in
the orbit method  \eqref{eqn:orbit}.

For a real reductive Lie group $G$, it has been observed by experts 
that there is no natural bijection between the unitary dual $\widehat{G}$
and (a subset of ) $\mathfrak{g}^*/\mathrm{Ad}^*(G)$, and
the orbit method does not work perfectly.
For instance, it is notoriously difficult to give a natural interpretation
of complementary series representations from the orbit philosophy.\footnote{
There is some recent trial, e.g., a geometric construction of the full 
complementary series representations of $SO(n,1)$ is proposed in 
\cite[Thm.\ 3.4]{xhkm} along the orbit philosophy.}
Nevertheless, the set $\mathfrak{g}^*/\mathrm{Ad}^*(G)$ of coadjoint orbits 
gives a rough approximation of the unitary dual $\widehat{G}$. 
In particular, 
geometric quantization of semisimple coadjoint orbits is  quite satisfactory,
and this provides a considerable
part of the unitary dual $\widehat G$.
In fact, when the coadjoint orbit
$\mathcal{O}_\lambda$ is semisimple satisfying certain mild conditions
including integrality of $\lambda$,  one can define an irreducible unitary 
representation of $G$ as a \lq\lq{geometric quantization}\rq\rq\
of $\mathcal{O}_\lambda$, as one can see by combining
analytic/geometric results with
algebraic representation theory
\cite{xknappv,  Vogan81, Vogan84, WaI, Wong}
in the 1950s--1990s
in the representation theory of Lie groups.
Roughly speaking, the irreducible unitary representations obtained as 
the \lq\lq{geometric quantization}\rq\rq\ of semisimple orbits
${\mathcal{O}}_{\lambda} \in {\mathfrak{g}}^{\ast}/\operatorname{Ad}^{\ast}(G)$
are as follows:
\begin{alignat*}{3}
&\text{Spherical unitary principal series representation}
&& \cdots\,\,\,\,
&& \text{$\mathcal{O}_{\lambda}$ is a hyperbolic orbit;} 
\\
&\text{Zuckerman's derived functor module $A_{\mathfrak{q}}(\lambda)$}\,\,
&& \cdots\,\,
&& \text{$\mathcal{O}_{\lambda}$ is an elliptic orbit;} 
\\
&\text{Tempered representation}
&& \cdots\,\,
&& \text{the dimension of $\mathcal{O}_{\lambda}$ is maximal}
\\
&
&& 
&& (=\frac{1}{2}(\dim{\mathfrak{g}}-\operatorname{rank}{\mathfrak{g}})).
\end{alignat*}
Here, one conducts the geometric quantization by using 
the real polarization for hyperbolic orbits and 
the complex polarization for elliptic orbits.
In particular, when $G$ is a connected compact Lie group,
orbits $\mathcal{O}_{\lambda}$ are always elliptic orbits and further
they are compact K{\"a}hlar manifolds. In this case
the geometric quantization of ${\mathcal{O}}_{\lambda}$ 
that puts  ``integral conditions'' to $\lambda$ can be constructed by
the Borel--Weil--Bott theorem. This is equivalent to the opposite manipulation 
of taking highest weights (Example \ref{ex:CWhw}).
More details such as terminology, an explicit construction, 
and some delicate issues for singular parameters 
may be found in the expository paper \cite[Sect.\ 2]{K94b},
see also \cite[Commentary]{Sugiura}.

On the other hand, 
\lq\lq{geometric quantizations}\rq\rq\ 
of nilpotent orbits
have yet to be fully elucidated.
For minimal nilpotent orbits, there have been some new progress
in constructing geometric models
of the \lq\lq{corresponding}\rq\rq\ unitary representations 
({\bf{minimal representations}}) including
the $L^2$-model 
(the ``Schr{\"o}dinger model'')
on the Lagrangian submanifold in the nilpotent orbits and 
the global analysis with motif in minimal representations
in the last 20 years, 
in which the author himself has been also involved
\cite{BKO, frenkel60, xkmanoAMS}, see Section \ref{subsec:6.4}.
A recent progress of the geometric quantization
of nilpotent orbits shows some interesting interactions
with other areas of mathematics, 
which would deserve a separate survey.

\subsection{Classification of irreducible unitary representations: a role of reductive Lie groups}
\label{subsec:Duflo}

The Mackey theory analyzes (irreducible) unitary representations of a group
via those of normal subgroups and of their quotients. The following
theorem is proven by using the Mackey theory for group 
extension and by the structure of algebraic Lie groups.
 
\begin{theorem}[Duflo \cite{Duflo82}]\label{thm:Duflo}
The classification of irreducible unitary representations of 
real algebraic groups reduces to the classification problem of 
the unitary duals of real reductive Lie groups.
\end{theorem}

\par
The reductive Lie groups required in Theorem \ref{thm:Duflo} 
are abelian when $G$ is a nilpotent Lie group as it 
is obtained as an iteration of group extensions by abelian Lie groups.
Thus, any irreducible unitary representation of a nilpotent Lie group $G$ 
is built up from irreducible unitary representations of abelian groups,
namely, from one-dimensional unitary representations, and this is how 
Kirillov classified the unitary dual of simply-connected nilpotent Lie groups 
(Example \ref{ex:Kirillov}).
For the classification of the unitary dual for a general Lie group $G$,
it is sufficient to determine the unitary dual when $G$ is a simple Lie group
by Theorem \ref{thm:Duflo}.
We will review the current status of this long-standing 
problem in the next subsection.

\subsection{Classification theory of irreducible representations of reductive Lie groups}
\label{subsec:Irr}

In this subsection we 
overview what has been known about
the classification problem for 
irreducible representations of simple Lie groups (or slightly more  generally
reductive Lie groups).
Although we highlight
the unitary dual $\widehat G$ of a reductive Lie group
$G$, we also consider a smaller set $\widehat {G}_{\operatorname{temp}}$
(irreducible tempered representations, Definition \ref{def:coarseLX})
and a larger set $\operatorname{Irr}(G)$
(irreducible admissible representations, Definition \ref{def:adm}).
\[
   \widehat {G}_{\operatorname{temp}}
   \subset
   \widehat {G}
   \subset
   \operatorname{Irr}(G).
\]

\begin{example}
For $G={\mathbb{R}}$,
we set $\chi_{\xi}\colon G\to \mathbb{C}^\times$, 
$x\mapsto e^{\sqrt{-1}x\xi}$
for $\xi \in \mathbb{C}$. Then we have the following bijections:
$\widehat G_{\operatorname{temp}}
=\widehat G \simeq \{\chi_{\xi} : \xi \in {\mathbb{R}}\} \simeq {\mathbb{R}}$, 
${\operatorname{Irr}}(G)\simeq \{\chi_{\xi} : \xi \in {\mathbb{C}}\} \simeq {\mathbb{C}}$.
\end{example}

In order to define ${\operatorname{Irr}}(G)$,
the set of equivalence classes of
\lq\lq{irreducible admissible representations}\rq\rq,
we need to clarify what ``admissibility'' means.
There are several ways for the characterization of
admissibility (\cite[Chap,\ 3]{WaI}, see also Remark \ref{rem:2.14} (1) below).
We choose one of the equivalent definitions of admissibility,
which fits with the branching law, a main theme of this article.
For this, we recall that there are at most
countably many equivalence classes of irreducible representations of 
compact Lie groups. We also recall that
any continuous representation of a compact Lie group defined over 
a suitable topological vector space (for instance,
a Hilbert space,
or more generally,
a complete locally convex topological vector space)
contains the algebraic direct sum of irreducible representations
as a dense subspace ({\bf{discrete decomposition}}),
but the number of appearances of the same irreducible representations  
(\textbf{multiplicity}) may vary from 0 to infinity.

\begin{defin}
[Harish-Chandra's admissible representation]
\label{def:adm}
A continuous representation $\Pi$ of a reductive Lie group $G$ is said to be
\emph{admissible} when the restriction $\Pi|_{K}$ to a maximal compact subgroup $K$ of $G$
contains each irreducible representation of $K$ with at most finite multiplicity.
\end{defin}

In Section \ref{subsec:mfbranch}, we shall extend the notion 
\lq\lq{admissibility}\rq\rq\ to
the restriction with respect to a more general (non-compact) subgroup, 
and refer to Harish-Chandra's admissibility as $K$-admissibility.

In analysis one may consider various function spaces  and equip them with 
natural topologies such as the Banach space topology for $L^p$-functions 
or the Fr{\'e}chet space topology for  $C^{\infty}$-functions on a manifold.
Analogously, in analytic representation theory, 
when a continuous representation $\Pi$ of a Lie group $G$
is defined over a Banach space $V$
(or more generally, a complete locally convex topological vector space),
one may consider the space $V^\infty$ of {\bf{$C^{\infty}$-vectors}}:
\[
   V^{\infty}:=\{v \in V : G \to V, \;
\text{$g \mapsto \Pi(g)v$ is a $C^{\infty}$-map}\}.
\]
The space $V^\infty$ is a $G$-invariant dense subspace of $V$
and also the differential representation $d\Pi$ of the Lie algebra $\mathfrak{g}$
can be defined on $V^\infty$.
Then $V^\infty$ is complete with respect to the family of 
seminorms given by $\| d\Pi(X_1)\cdots d\Pi(X_k)v \|$ 
for $X_1,\ldots, X_k \in \mathfrak{g}$.
With respect to this topology the map $G \times V^\infty \to V^\infty$,
$(g, v) \mapsto \Pi(g)v$ is continuous, hence
one obtains a continuous representation of $G$ on the Fr{\'e}chet space $V^\infty$, to be denoted by $\Pi^\infty$.

A vector $v \in V$ is called $K$-finite if 
$\dim_\mathbb{C} \text{$\mathbb{C}$-span}\{\Pi(k)v\}<\infty$.
We write $V_K$ for  the space of $K$-finite vectors. 
If $(\Pi,V)$ is an admissible representation on a Banach space $V$,
then $V_K \subset V^\infty$, hence one can define the 
action of the Lie algebra $\mathfrak{g}$ as well as that of the 
maximal compact subgroup $K$ on $V_K$,
to which we refer to as the \emph{underlying $(\mathfrak{g},K)$-module}
of $(\Pi,V)$.

In this article, we adopt the Fr{\'e}chet representation
of \emph{moderate growth}  in defining 
the set ${\operatorname{Irr}}(G)$, which fits with the study of 
symmetry breaking operators in later sections.
Some justification
is given in Remark \ref{rem:2.14} via (essentially) equivalent definitions.

\begin{defin}[Irreducible admissible representation]\label{def:irradm}
Suppose $G$ is  a reductive Lie group.
If $(\Pi,V)$ is an irreducible admissible representation of $G$ on a 
Banach space $V$, then  the Fr{\'e}chet representation $(\Pi^\infty, V^\infty)$
is also irreducible, and called of \emph{moderate growth}
from the behavior of the matrix coefficients, see \cite[Chap.\ 11]{WaI}.
We denote by ${\operatorname{Irr}}(G)$ the set of 
equivalence classes of such irreducible admissible Fr{\'e}chet representations
$(\Pi^{\infty}, V^{\infty})$.
\end{defin}

\begin{remark}\label{rem:2.14}
{\rm{(1)}}
(Harish-Chandra)
\enspace
If $(\Pi,V)$ is an irreducible representation realized on a Banach space $V$,
then the following are equivalent.
\begin{equation}\label{eqn:adm}
\begin{aligned}
&\text{The center of the universal enveloping algebra $U({\mathfrak{g}})$ acts on
$V^{\infty}$ by scalar.}\\
&  \iff
\text{$(\Pi,V)$ is an admissible representation.}
\end{aligned}
\end{equation}
\par\noindent
{\rm{(2)}}\enspace
If $(\Pi_1, V_1)$ and $(\Pi_2, V_2)$ are two irreducible admissible
representations of $G$ on Banach spaces $V_1$ and $V_2$.
Then the underlying $(\mathfrak{g},K)$-modules $(\Pi_1)_K$ and $(\Pi_2)_K$
are isomorphic to each other as $(\mathfrak{g},K)$-modules
if and only if $(\Pi_1)^\infty$ and $(\Pi_2)^\infty$
are isomorphic as Fr{\'e}chet representations of $G$. We note that 
$\Pi_1$ and $\Pi_2$ are not necessarily isomorphic 
as Banach representations of $G$.
\par\noindent
{\rm{(3)}}\enspace
For every irreducible $(\mathfrak{g},K)$-module $X$, there exists
$\Pi \in {\operatorname{Irr}}(G)$ such that $\Pi_K \simeq X$, hence
one has a natural bijection between
${\operatorname{Irr}}(G)$
and the set of equivalence classes of irreducible $({\mathfrak{g}},K)$-modules.
\par\noindent
{\rm{(4)}}\enspace
Via the correspondence $(\Pi, V) \rightsquigarrow (\Pi^{\infty}, V^{\infty})$,
the unitary dual $\widehat G$ can be regarded as a subset of ${\operatorname{Irr}}(G)$.\end{remark}

The statements (2) and (3) in Remark \ref{rem:2.14} are part of
Casselman--Wallach's globalization theory \cite[Chap.\ 11]{WaI}.

\vskip 1pc
\par\noindent
{\bf{The classification of irreducible unitary representations of 
reductive Lie groups}}: 
The classification problem of the unitary dual $\widehat G$
of a reductive Lie group $G$ has a history of over 70 years.
It is completed for some special cases such as 
$SL(n,{\mathbb{F}})$ (${\mathbb{F}}={\mathbb{R}}, {\mathbb{C}}, {\mathbb{H}}$)
explained in Example \ref{ex:SLunitary},
complex Lie groups $SO(n,{\mathbb{C}})$ and $Sp(n,{\mathbb{C}})$
or some real reductive Lie groups of low rank.
Moreover, recently, 
the \lq\lq{Atlas project}\rq\rq\ led by 
Adams, van Leeuwen, Trapa, Vogan, and others,
which aims to give a \lq\lq{description}\rq\rq\ of $\widehat G$ 
by finite algorithms, has been also developed \cite{ALTV20}.
However, the unitary duals of simple Lie groups remain to be 
fully elucidated even for classical cases such as indefinite orthogonal groups 
$O(p,q)$ with general $p$, $q$.

As we mention below in Example \ref{ex:KZtemp},
the classification of $\widehat G_{\operatorname{temp}}$
(tempered representations),
which is a subset of the unitary dual $\widehat G$, 
is completed by Knapp--Zuckerman.

\par\noindent
{\bf{The classification of irreducible representations of reductive Lie groups
(when no unitarity is imposed)}}:
In contrast to the long-standing problem on the classification of the unitary dual $\widehat G$,
the classification of $\operatorname{Irr}(G)$, 
which is larger than $\widehat G$,
was completed from 1970s to the early 1980s.
(In other words, 
the current status of the classification of $\widehat G$ is that 
although the classification of $\operatorname{Irr}(G)$ is completed,
we do not fully understand the whole picture of the subset
$\widehat{G}$ of $\operatorname{Irr}(G)$.

The classification of the set of infinitesimal equivalence classes of irreducible
admissible representations (equivalently, that of $\operatorname{Irr}(G)$ as 
Fr{\'e}chet modules, see Remark \ref{rem:2.14} (3)) is
reduced to that of irreducible $({\mathfrak{g}},K)$-modules.
There are three approaches to the classification of irreducible 
$(\mathfrak{g},K)$-modules.

\begin{enumerate}
\item[$\bullet$] 
 (Langlands classification)
This is an analytic approach that 
reduces to the classification of $\widehat G_{\operatorname{temp}}$ by 
focusing on the asymptotic behavior of matrix coefficients.
The classification of $\widehat G_{\operatorname{temp}}$
(irreducible tempered representations) was accomplished by Knapp--Zuckerman.
\vskip 0.05in

\item[$\bullet$]
(Vogan's classification) This employs a purely algebraic method
that uses minimal $K$-type theory, Zuckerman's derived functor modules
(an algebraic generalization of the Borel--Weil--Bott theory),
and the Lie algebra cohomology 
(a generalization of highest weights).
\vskip 0.05in

\item[$\bullet$]  
(${\mathcal{D}}$-module approach)
This method has a geometric feature: it reduces the classification
of irreducible $(\mathfrak{g},K)$-modules (with regular infinitesimal characters)
to that of irreducible modules of the ring of twisted differential operators
over flag manifolds (Beilinson--Bernstein, Brylinski--Kashiwara)
and to the geometry of $K_{\mathbb{C}}$-orbits.
\end{enumerate}

\subsection{A role of the branching law on the classification of irreducible representations}
\label{subsec:2.7}

The concept of the \lq\lq{smallest units}\rq\rq\ varies
depending on the viewpoint.
Taking subgroups changes the viewpoint in representation theory:
irreducible representations of a group are no longer 
\lq\lq{smallest units}\rq\rq\ 
when the action is restricted to subgroups.
The {\bf{theory of branching law}} aims to elucidate
structures of representations
from the viewpoint of subgroups. Conversely,
its idea plays a useful role
in the classification theory of irreducible representations itself
\cite{deco-euro}, as mentioned in Section \ref{sec:Intro}.
Indeed, the idea of branching laws is used
to define \lq\lq{invariants}\rq\rq\ of irreducible representations,
e.g.,
in Cartan--Weyl's highest weight theory
for the classification of irreducible finite-dimensional representations
and also in
Vogan's theory for that of irreducible $({\mathfrak {g}}, K)$-modules, where
the restriction to a maximal torus and a maximal compact subgroups
respectively, is discretely decomposable and 
the \lq\lq{edge}\rq\rq\ of the branching law (with respect to a certain partial order)
gives \lq\lq{invariants}\rq\rq\ in the classification.
Further, sometimes, 
\lq\lq{new irreducible unitary representations}\rq\rq\
have been discovered in the process of finding branching laws
such as via the theta correspondence (Example \ref{ex:Howe}).

\section{From Local to Global---Global Analysis on Manifolds with Indefinite Metric}
\label{sec:1}

The theory of discontinuous groups beyond the Riemannian setting
is another young field that has remarkably developed in the last 30 years.
The ideas from discontinuous groups inspired the author
at several turning points in inventing the theory of 
the restriction of 
infinite-dimensional representations.
In this section we shall shed light on the parts in which both theories are related.

\subsection{Mysterious phenomena on the global analysis for indefinite metrics}
\label{subsec:2.1}

A {\textit{pseudo-Riemannian manifold}} is 
a manifold $M$ equipped with
a non-degenerate quadratic form $g_x$
at the tangent space $T_xM$ depending
smoothly on $x \in M$.
When $g$ is positive definite, $(M,g)$ is a Riemannian manifold.
When only one eigenvalue of  $g$ is negative, 
it is called a Lorentzian manifold,
which appears as a geometric structure of spacetime in general relativity.
The group of isometries of a pseudo-Riemannian manifold is automatically
a Lie group.
Semisimple symmetric spaces (in particular, irreducible symmetric spaces)
are examples of pseudo-Riemannian manifolds,
on which semisimple Lie groups act transitively as isometric transformations.

The motif 
\[
\text{local property $\rightsquigarrow$ global form}
\] 
has been one of the main streams in geometry since the 20th century,
especially in Riemannian geometry.
The question that 
when a local structure is fixed,
``how flexible the global form is, 
or conversely, what kind of limitations exists in the global form?'' 
is a prototype of the motif, which would bring us to deformation theory
and rigidity theorems, respectively.

The study of ``local to global'' interacts with various fields of mathematics,
depending on the types of local properties of interest.
If one highlights \lq\lq{local homogeneity}\rq\rq\ 
as a local property,
Lie theory and number theory will enter naturally in this study through 
an algebraic structure called {\bf{discontinuous groups}},
which controls the global nature of such manifolds.
In the classical case of Riemannian manifolds (where the metric $g$ is positive definite),
the study of discontinuous groups
has already entered its golden age by the 1950s.
This study interacts with various fields, including
Riemannian symmetric spaces, Lie theory, number theory,
differential geometry, and topology.

On the other hand, the study of
pseudo-Riemannian geometry seemed to be behind the trend of 
``local to global''.

An early work is due to Calabi--Markus \cite{CM} in 1960s for 
de Sitter manifolds.
A systematic study of discontinuous groups 
for pseudo-Riemannian homogeneous manifolds
started with \cite{kob89} in the late 1980s.
As an introduction to this theme, we begin this section with
three phenomena
in pseudo-Riemannian (locally homogeneous) geometry
that look strange from the \lq\lq{classical point of view}\rq\rq\ 
in Riemannian geometry, see also the research survey paper \cite{xk2001} 
for more details of the first two of them.

\vskip 1pc
\par\noindent
\ref{subsec:2.1}.(1)\enspace {\bf{Curvature and global form}:}
Curvatures are typical examples of local invariants of (pseudo-)Riemannian
manifolds.
To see what kind of constraints the curvature (local property)
gives to the global form of a manifold,
let us compare
classical results on Riemannian manifolds
(Theorems \ref{thm:2.1} and \ref{thm:hypcpt})
and 
some different features on pseudo-Riemannian manifolds
(Theorems \ref{thm:cm} and \ref{thm:AdScpt}).
A Riemannian manifold, or more generally, 
a pseudo-Riemannian manifold is called a \emph{space form}
if its sectional curvature $\kappa$ is constant.
In the Riemannian case, it is called a hyperbolic manifold
when $\kappa < 0$. In the Lorentzian case, it is 
a {\bf{de Sitter manifold}} or
an {\bf{anti-de Sitter manifold}}  when $\kappa > 0$ or $\kappa < 0$,
respectively.
We highlight the contrast between Riemannian and Lorentzian cases:
Theorems \ref{thm:2.1} and  \ref{thm:cm} are for positive 
curvatures, whereas Theorems \ref{thm:hypcpt} and \ref{thm:AdScpt}
are for negative curvatures.

\begin{theorem}[Myers]\label{thm:2.1}
Every complete Riemannian manifold with Ricci curvature $\ge \varepsilon (>0)$ is compact.
\end{theorem}

\begin{theorem}
[Calabi--Markus phenomenon {\cite{CM}}]
\label{thm:cm}
Every de Sitter manifold is non-compact.
\end{theorem}

\begin{theorem}\label{thm:hypcpt}
There exists a compact hyperbolic manifold for every dimension.
\end{theorem}

\begin{theorem}\label{thm:AdScpt}
Compact anti-de Sitter manifolds exist if and only if the 
dimension is odd.
\end{theorem}

\vskip 1pc
\par\noindent
\ref{subsec:2.1}.(2)\enspace 
{\bf{Rigidity of discontinuous groups and \lq\lq{deformability}\rq\rq:}}
Can we \lq\lq{deform}\rq\rq\ discontinuous groups for pseudo-Riemannian
symmetric spaces $G/H$? Since automorphisms of $G$ induce
\lq\lq{uninteresting}\rq\rq\ deformation of discontinuous groups,
we consider \lq\lq{deformation}\rq\rq\ up to automorphisms of $G$.
Let us compare a classical rigidity theorem
(Theorem \ref{thm:2.4})
when the metric tensor $g$ is positive definite 
with a discovery on the \lq\lq{flexibility of discontinuous groups}\rq\rq\ 
(Theorem \ref{thm:K98d})
when $g$ is indefinite.
\vskip 0.8pc
\begin{theorem}[Selberg--Weil's local rigidity \cite{Weil64}]\label{thm:2.4}
Any cocompact discontinuous group for an irreducible Riemannian symmetric space of
dimension $>2$
does not allow non-trivial continuous deformation.
\end{theorem}

\begin{theorem}[Kobayashi {\cite{K1998d}}]\label{thm:K98d}
There exists an irreducible symmetric space of 
arbitrarily higher dimension having a cocompact discontinuous group
that allows non-trivial continuous deformation.
\end{theorem}

\vskip 0.8pc

The local rigidity theorem by Selberg and Weil in
Theorem \ref{thm:2.4}
has brought a sequence of revolutionary works on
rigidity theorems by Mostow, Margulis, and Zimmer among others,
whereas in deformation theory of discontinuous groups for 
the Poincar{\'e} disk (the two-dimensional case)
the last century
has seen the development of an abundant theory of 
{\bf{Teichm{\"u}ller spaces}} on moduli spaces of Riemann surfaces.
Theorem \ref{thm:K98d} concerning indefinite metrics may be thought of as 
providing a new theme such as
\lq\lq{higher dimensional Teichm{\"u}ller theory}\rq\rq\ for 
locally semisimple symmetric spaces with indefinite metric.

\vskip 1pc
\par\noindent
\ref{subsec:2.1}.(3)\enspace {\bf{Rigidity of spectrum}:}
So far we have seen distinguished aspects in the geometry of discontinuous 
groups for pseudo-Riemannian manifolds.
We now address an analytic question: 
Suppose that a discontinuous group $\Gamma$ for 
a (pseudo-)Riemannian manifold $X$
admits a non-trivial continuous deformation. 
Do the eigenvalues of the Laplacian on
the quotient space 
$\Gamma \backslash X$ 
vary according to deformation of $\Gamma$, or are there 
\lq\lq{stable eigenvalues}\rq\rq?
We consider this problem 
in the setting where $X$ is a space form, i.e., 
a pseudo-Riemannian manifold
with  sectional curvature $-1$.
Theorem \ref{thm:2.6} is a classical result for Riemannian manifolds
of dimension two
and Theorem \ref{thm:2.3} is a new phenomenon for the 
Lorentzian manifold of dimension three
discovered in a joint work with
F.~Kassel \cite{Adv16}. See also \cite{xkIntrinsic}.

\begin{theorem}[Wolpert {\cite{Wolpert}}]\label{thm:2.6}
There does not exist a \lq\lq{stable eigenvalue}\rq\rq($>1/ 4$)
on closed Riemann surfaces.
\end{theorem}

\begin{theorem}[Kassel--Kobayashi \cite{Adv16}]\label{thm:2.3}
There exist infinitely many 
positive
\lq\lq{stable eigenvalues}\rq\rq
on three-dimensional compact anti-de Sitter manifolds.
\end{theorem}

Here we say $\lambda$ is a \textit{stable eigenvalue} if there exists
a non-zero $L^2$-function $f$ (depending on $\varphi$) on the quotient
manifold $\varphi(\Gamma)\backslash X$ satisfying
$\Delta f = \lambda f$ (weak sense) for every injective homomorphism
$\varphi \colon \Gamma \to G$ that is sufficiently close to the original
embedding $\iota \colon \Gamma \hookrightarrow G$.

\begin{remark}
The existence of 
\lq\lq{stable eigenvalues}\rq\rq\ in Theorem \ref{thm:2.3}
is proved in \cite{Adv16} by constructing $\Gamma$-periodic eigenfunctions 
that are obtained as the \lq\lq{$\Gamma$-average}\rq\rq\ of 
non-periodic eigenfunctions.
For the proof of the convergence and the non-vanishing of the 
\lq\lq{$\Gamma$-average}\rq\rq,
one uses a geometric estimate such as the 
\lq\lq{counting}\rq\rq\
of the $\Gamma$-orbit in a pseudo-ball of
a pseudo-Riemannian symmetric space
(Section \ref{subsec:cpt}) and 
analytic estimate of eigenfunctions, see Section \ref{subsec:L2X}
for the non-commutative harmonic analysis.
On the other hand, by using the branching law (Section \ref{sec:A}) of 
infinite-dimensional representations, it can be shown that 
there also exist infinitely many eigenvalues ($<0$) that \lq\lq{vary}\rq\rq\ 
according to deformation of a discontinuous group,
see \cite{KaKo20, K16b} for a precise formulation.
\end{remark}

\subsection{Inspiration from the theory of discontinuous groups}
Let us explain a loose
idea that connects 
the ``theory of discontinuous groups of pseudo-Riemannian manifolds''
for which some mysterious phenomena are
described in Section \ref{subsec:2.1} with
the ``theory of the restriction of infinite-dimensional representations''
which is the main theme of this article.

First, we review basic notion on group actions.
Suppose a topological group $\Gamma$ 
acts continuously on a locally compact space $X$.
We define a subset 
$\Gamma_S \subset \Gamma$ for a given subset $S \subset X$
as
\[
   \Gamma_S :=\{\gamma \in \Gamma : \gamma S \cap S \ne \emptyset \}.
\]
When $S$ is a singleton $\{x\}$, the subset $\Gamma_S$ is a subgroup.
Let us recall the following basic concepts.

\begin{defin}\label{def:proper}
\par\indent
{\rm{(1)}}\enspace
An action is  {\bf{properly discontinuous}}
$\iff$
$\# \Gamma_S < \infty$ for any compact set $S$.
\par\indent
{\rm{(2)}}\enspace
An action is {\bf{proper}}
$\iff$
$\Gamma_S$ is compact for any compact set $S$.
\par\indent
{\rm{(3)}}\enspace
An action is {\bf{free}}
$\iff$
$\Gamma_{\{x\}} = \{e\}$ for any $x \in X$.
\end{defin}

One may think of each of the above three concepts as a kind of 
properties that $\Gamma_S$ is reasonably ``small'' whenever $S$ is ``small''.

If $\Gamma$ acts on a manifold $X$ properly discontinuously and freely,
then the quotient space $\Gamma \backslash X$ 
is a Hausdorff space with respect to 
quotient topology and,  further, there exists a unique manifold structure
in $\Gamma \backslash X$ for which the quotient map 
$X \to \Gamma \backslash X$ is a smooth covering.
Conversely, the fundamental group $\pi_1(M)$  of a manifold $M$
acts properly discontinuously and freely 
on the universal covering space $\widetilde M$ as 
covering transformations and the original manifold $M$ can be recovered 
as $\pi_1(M) \backslash \widetilde M \simeq M$.

Suppose that $X$ is a pseudo-Riemannian manifold with
metric tensor $g$, and that $G$ is the group of isometries of $X$.
Then $G$ is always a Lie group. For
 a subgroup $\Gamma$ of $G$,
the following equivalence does \emph{not} hold generally
even when $\Gamma$ acts freely on $X$
(the implication $\Leftarrow$ always holds):
\[
\text{$\Gamma$ is discrete (in $G$)
\,\,
$
\Longleftrightarrow
$ 
\,\,
the action of $\Gamma$ on $X$ is properly discontinuous}.
\]
This is a significant difference from Riemannian manifolds
with $g$ positive definite, where the equivalence automatically holds,
and 
the failure of the implication $\Rightarrow$ in the pseudo-Riemannian case
is one of the causes of \lq\lq{mysterious phenomena}\rq\rq\ 
described in Section \ref{subsec:2.1}.
Thus it is crucial to gain a profound understanding
(not a formal paraphrase of the definition
of proper discontinuity)
when 
a discrete subgroup of $G$ acts on a non-Riemannian homogenous space.
An explicit properness criterion is known for reductive Lie groups $G$, 
which we recall now.

\vskip 1pc

\par\indent
Let 
$G=K \exp({\mathfrak{a}})K$ 
be the Cartan decomposition 
of a reductive Lie group $G$
and 
$\mu \colon G \to {\mathfrak{a}}/W$
the corresponding Cartan projection. 
Here $W$ stands for the Weyl group 
of the restricted root system of the Lie algebra $\mathfrak{g}$ with 
respect to the maximal abelian split subalgebra 
${\mathfrak{a}}$.
The next theorem extends the properness criterion
given by the author \cite{kob89} in the setting where both
$\Gamma$ and $H$ are reductive subgroups.

\begin{theorem}
[Criteria for proper discontinuity; Benoist \cite{Bn1996}, Kobayashi \cite{kob96}]
\label{thm:3.2}
Let $\Gamma$ be a discrete subgroup of a reductive Lie group $G$
and $H$ a closed subgroup of $G$. Then
the following two conditions on $(\Gamma, G, H)$ are equivalent:
\begin{enumerate}
\item[\textnormal{(i)}]
the action of $\Gamma$ on $G/H$ is properly discontinuous;
\vskip 0.05in
\item[\textnormal{(ii)}]
the set $\mu(\Gamma)\cap\mu(H)_{\varepsilon}$ is a finite set
for any $\varepsilon>0$, where
$\mu(H)_{\varepsilon}$ is a tubular neighborhood of $\mu(H)$
in $\mathfrak{a}/W$.
\end{enumerate}
\end{theorem}

\par
Theorem \ref{thm:3.2} plays a key role to prove
the necessary and sufficient condition for 
the Calabi-Markus phenomenon for reductive homogeneous spaces
$G/H$ (cf.\ Theorem  \ref{thm:cm}) in \cite{kob89} and
the existence problem of compact pseudo-Riemannian 
locally homogeneous spaces such as space forms
(Theorem  \ref{thm:AdScpt} for the Lorentzian case),
which were introduced in Section \ref{subsec:2.1} 
as \lq\lq{mysterious phenomena}\rq\rq\ for pseudo-Riemannian manifolds,
and also Theorem \ref{thm:K98d} (the deformation of discontinuous groups
in higher dimension).

The criterion for discrete decomposability for the branching
law of unitary representations (Section \ref{sec:A}) was inspired
by the properness criterion \cite{kob89} in its formulation.
Although the techniques to prove Theorems 
\ref{thm:3.2} (topology) and \ref{thm:2.2} 
(decomposition of a Hilbert space, micro-local analysis)
are quite different, 
there are some common characteristics
in both cases, that is, one searches a setting  in which 
non-compact subgroups behave {\bf{as if they were compact groups}}.\footnote{Not just
\lq\lq{apparent similarities}\rq\rq,
a close relationship between the
\lq\lq{proper action}\rq\rq\
and the
\lq\lq{discretely decomposable restriction}\rq\rq\
is recently elucidated in \cite{K16b} under the assumption
of certain \lq\lq{hidden symmetry}\rq\rq\ on spherical varieties.}

\subsection{Application from the theory of discontinuous groups:
from qualitative theorems to quantitative estimates}
\label{subsec:cpt}

Properness (or proper discontinuity) of the action is initially
a qualitative property. As a second step, we may deepen such 
qualitative properties to {\bf{quantitative}} estimates,
from which a further connection to another field
of mathematics emerges.
Let us give two such examples.

\vskip 0.1in

\noindent
{\bf{(1)}}\enspace{\bf{Quantitative estimates of proper discontinuity}}

If a discrete group $\Gamma$ acts on $X$ properly discontinuously
(Definition \ref{def:proper}), then the following inequality holds for any 
compact subset $S \subset X$:
\[
   \# \Gamma_S < \infty.
\]
Then we may consider a generalization of the classical counting of lattice points in 
non-Riemannian geometry as follows.

\begin{problem}\label{prob:3.11}
When a compact set $S$ is gradually increased,
how does $\# \Gamma_S$ increase?
\end{problem}

For instance,
if $X$ is a Riemannian manifold and $S$ is a ball $B(R)$ of radius $R$
centered at $o \in X$, then $\Gamma_{B(R)}$ coincides with the set
$\{ \gamma \in \Gamma : \gamma \cdot o \in B(2R)\}$, and 
Problem \ref{prob:3.11} is a classical counting problem of a lattice
which asks about the asymptotic behavior of 
$\#\Gamma_{B(R)}$
as the radius $R$ tends to infinity.
Such an estimate for a pseudo-Riemannian manifold is 
utilized in the proof of
the existence (Theorem \ref{thm:2.3}) of 
\lq\lq{stable eigenvalues}\rq\rq\ for the Laplacian on
a locally pseudo-Riemannian symmetric space \cite{Adv16}.

\vskip 0.1in

\noindent
{\bf{(2)}}\enspace
{\bf{Quantitative estimates of non-proper actions}}

We may also formalize ``quantification'' of ``non-properness'' 
in the opposite situation where
the action of the group is not proper.
Along this direction we propose a new approach
in non-commutative harmonic analysis by using
an idea of dynamical system instead of the traditional method of
differential equations for the study of the unitary representation
on $L^2(G/H)$, see Section \ref{subsec:tempered}.

\section{Program for Non-commutative Harmonic Analysis}

In contrast to local analysis, 
global analysis involves many difficult problems.
Further, 
in order to establish global results, certain 
\lq\lq{appropriate structure\rq\rq\ (e.g., curvature pinching)
should be imposed on the manifold $X$ unless otherwise easy
\lq\lq{counterexamples\rq\rq\ show up because of the non-compactness of $X$.
Such a ``structure on $X$'' may be provided by means of
a non-compact transformation group $G$ of $X$. This is an idea to 
investigate global analysis successfully in the following scheme:
\[
\text{global analysis on $X$ $\longleftrightarrow$ representation theory of $G$.}
\] 

In this section we focus on global analysis via transformation groups. 
After giving a brief overview on
some major achievements 
about non-commutative harmonic analysis
in the 20th century in Section \ref{subsec:L2X},
we describe a program for a new line of studies in
Sections \ref{subsec:grip}--\ref{subsec:6.4}.
Its connection to the theory of branching laws will be dealt with
in Sections \ref{sec:Stage}--\ref{sec:C}.

\subsection{Non-commutative harmonic analysis 
on semisimple symmetric spaces
\hspace{1in}---perspectives and achievements in the 20th century}
\label{subsec:L2X}
In global analysis on a manifold $X$,
{\bf{non-commutative harmonic analysis}} 
concerns the space of functions rather than individual functions, 
and uses the regular representation of the transformation group $G$
defined on the space of all functions.

Let us review basic notions.
When a group $G$ acts on a manifold $X$,
a linear action of $G$ on 
the space $\Gamma(X)$ of functions  
($\Gamma =C^{\infty}$, $C$, ${\mathcal{D}}'$, $\ldots$)
is naturally defined. That is, for each $g \in G$,
one transforms functions $f$ on $X$ to 
different functions
$
  \pi(g) f:= f(g^{-1} \cdot)
$
by the pull-back of the geometric action.
Since the resulting family of
the linear map $\pi(g) \colon \Gamma(X) \to \Gamma(X)$ satisfy
 $\pi(g_1 g_2)=\pi(g_1)\pi(g_2)$
with respect to compositions of maps,
$\pi$ defines a representation of the group $G$ on 
the function space $\Gamma(X)$, which is referred to as
{\bf{regular representation}}.
Further, if there exists a $G$-invariant Radon measure on $X$,
then $\pi(g)$ preserves the $L^2$-norm and therefore one can define
a unitary representation on the Hilbert space $L^2(X)$
of square integrable functions on $X$.

\begin{remark}
Even when 
a manifold $X$ is too small to allow a non-trivial action of a group
(no symmetry in geometry),
it may be possible to define a representation of the group on
the space of functions (symmetry in analysis). 
This viewpoint leads us to
 ``global analysis with minimal representations as a motif''
\cite{sato50},
where one utilizes \lq\lq{hidden symmetry}\rq\rq\ on the space of functions.
This new direction of global analysis 
will be mentioned in Section \ref{subsec:6.4}.
\end{remark}

Before explaining what \lq\lq{non-commutative harmonic analysis}\rq\rq\
means, we first review classical \lq\lq{commutative}\rq\rq\ harmonic analysis.
The Fourier transform on the Euclidean space
\begin{equation}
\label{eqn:F}
{\mathcal{F}}\colon  C_c({\mathbb{R}})\to C({\mathbb{R}}), 
\qquad
({\mathcal{F}} f)(\xi)
:=
\frac{1}{\sqrt{2 \pi}} \int_{-\infty}^{\infty} f(x)e^{-ix \xi} d x
\end{equation}
initially defined, for instance, in the space $C_c(\mathbb{R})$ of compactly 
supported continuous functions, extends to a unitary operator on
the Hilbert space $L^2({\mathbb{R}})$ of square-integrable functions
({\bf{the Plancherel theorem}}).
On the other hand, if we put $\pi(t)f:=f(\cdot -t)$, then
the Fourier transform ${\mathcal{F}}$ also satisfies the following
algebraic relations
\[
{\mathcal{F}}(\pi(t) f)(\xi)
=
e^{-i t \xi}({\mathcal{F}} f)(\xi)
\qquad
({}^{\forall} t \in {\mathbb{R}}).
\]
These properties of the Fourier transform ${\mathcal{F}}$ 
may be reinterpreted from the standpoint of the representation theory of groups
as follows.
For later purpose, we separate the role of the transformation group $G$
from that of the geometry $X$, and consider $G$ acts on $X$ and $L^2(X)$,
even though $G=X=\mathbb{R}$ here.
\begin{alignat*}{2}
&&\text{Algebraic relations}\cdots 
&\text{The Fourier transform $f \mapsto {\mathcal{F}}f(\xi)$ gives 
a $G$-homomorphism}\\
&& &\text{for each $\xi \in \mathbb{C}$ from 
the regular representation $\pi$ on $C_c(X)$ of} \\
&& &\text{the additive group $G$ to the representation space 
$\mathbb{C}$ of the}\\
&& &\text{irreducible representation 
$\chi_{-\xi} \colon G \to GL(1,{\mathbb{C}})$,  
$t \mapsto e^{-i t \xi}$ of $G$}.
\\
&&\text{$L^2$-theory\phantom{mmmmi}}
\cdots
& 
\text{The Plancherel theorem decomposes the regular representation}\\
&& &\text{$L^2(X)$ into irreducible unitary representations of $G$.}
\end{alignat*}

When a transformation group $G$ is the abelian group $\mathbb{R}$,
the representation space of 
the irreducible representation $\chi_\xi$ of $G$ is one-dimensional,
hence the above group-theoretic interpretation may seem
stress too much on formalism.
However, thanks to this Weyl's point of view,
the representation theory of transformation groups $G$ of manifolds $X$
have been incorporated with global analysis in a broad sense and 
led to a great development of its own.
The irreducible decomposition of the unitary
representation $L^2(X)$ of $G$ is called 
the {\bf{Plancherel-type theorem}}. 
The harmonic analysis on the Euclidean space is extended to that on
abelian locally compact groups $G$ (\lq\lq{commutative harmonic analysis}\rq\rq)
by Pontryagin in 1930s, 
where all irreducible unitary representations of $G$
are one-dimensional. In contrast, this analysis is called
\lq\lq{non-commutative harmonic analysis}\rq\rq\ if the transformation 
group $G$ is non-commutative.
When the transformation group $G$ is compact, 
then all irreducible representations 
are finite-dimensional
and $L^2(X)$ decomposes discretely as in the case of Fourier series
expansions in the case $G=X=S^1$. 
When the transformation group $G$ is non-abelian and non-compact,
the theory of ``infinite-dimensional irreducible representations'' of $G$
may be used for global analysis on the manifold $X$
as a powerful tool. 
Here are successful cases
for the Plancherel-type theorem for $G$-spaces $X$ with emphasis on
real reductive Lie groups $G$.

\begin{alignat*}{2}
&\text{$X=G$\enspace (group manifold)}
&&
\\
&
\text{\qquad
Peter--Weyl (1927)}
\qquad
&&\text{$G$ is a compact group.}
\\
&\text{\qquad
Pontryagin (1934)}
&&
\text{$G$ is an abelian locally compact group.}
\\
&\text{\qquad Gelfand school, Harish-Chandra (1950s)}
\quad
&&
\text{$G$ is a complex semisimple Lie group.}
\\
&\text{\qquad Harish-Chandra (1976)}
\qquad
&&\text{$G$ is a semisimple Lie group.}
\\
&\text{
$X=G/H$\enspace (symmetric space)}
&&
\\
&\text{\qquad T.\ Oshima (1980s), Delorme,} &&\\
&\text{\qquad van den Ban--Schlichtkrull (late 1990s)}
&&\text{$X$ is a semisimple symmetric space.}
\end{alignat*}

These great achievements in the 20th century 
were not limited to theorems in representation theory,
but also served as the driving forces for the development
of analysis such as functional analysis and algebraic analysis.
On the other hand, 
the beautiful and successful theory on global analysis on symmetric spaces
(including group manifolds) seemed to give an impression that these objects
live in a ``closed world''.
Here we note that group manifolds $G$ may be seen as a special
case of symmetric spaces $(G\times G)/\operatorname{diag}(G)$.
In fact, many of the techniques used in global analysis in there were based on
the structure theory of these spaces (e.g.\ the proof for group manifolds
or symmetric spaces can be often reduced to that for 
smaller groups or smaller symmetric spaces), and thus
it was not obvious to foresee a promising direction of global analysis
beyond symmetric spaces
in the 1980s around the time 
when the Plancherel-type theorem for reductive symmetric spaces
was announced by T.~Oshima \cite{Oshima2002}.

\subsection{``Grip strength'' of representations on global analysis --- new perspective, part 1}
\label{subsec:grip}
To find a nice framework beyond symmetric spaces, let us start 
from scratch, and 
consider the very basic question 
whether representation theory is \lq\lq{useful}\rq\rq\ 
for global analysis in the first place.
As a guiding principle to explore this question,
we propose the following perspective.

\begin{basicproblem}
[``Grip strength'' of representations \cite{Ksuron}]
\label{prob:6.1}
When a Lie group $G$ acts on $X$,
can the space of functions on $X$ be ``sufficiently controlled'' by the representation 
theory of $G$?
\end{basicproblem}

The vague words, ``sufficiently controlled'', or conversely,``uncontrollable'',
need to be formulated rigorously as mathematics.
Suppose $V$ is the space of functions of a $G$-manifold $X$.
There are (sometimes uncountably many) 
inequivalent irreducible subrepresentations 
in the infinite-dimensional $G$-module $V$.
Moreover the multiplicity of each irreducible representation 
can range from finite to infinite.

Confronting such a general situation, 
we emphasize the following principle:
\vskip 0.1in
\begin{itemize}
\item
even though there are infinitely many inequivalent irreducible 
subrepresentations,
the group action \textbf{can} distinguish the inequivalent parts;
\vskip 0.05in
\item
the group action \textbf{cannot} distinguish the parts where 
the same irreducible representations occur with multiplicities.
\end{itemize}
\smallskip

This observation suggests us to think of 
the multiplicity of irreducible representations as
the quantity measuring the \lq\lq{grip strength of a group}\rq\rq.
For an irreducible representation $\Pi$ of a group $G$,
the multiplicity of $\Pi$ in $C^{\infty}(X)$ is defined by
\begin{equation}
\label{eqn:multfn}
  \operatorname{dim}_{\mathbb{C}} \operatorname{Hom}_G
  (\Pi,C^{\infty}(X)) \in {\mathbb{N}} \cup \{\infty\}.
\end{equation}
We formulate Basic Problem \ref{prob:6.1} as follows.

\begin{problem}
[Grip strength of representations on global analysis]
\label{prob:6.2}
Let $X$ be a manifold on which a Lie group $G$ acts.
\begin{enumerate}
\item[{\rm{(1}})]
Find a necessary and sufficient condition on the pair $(G,X)$ for which
the multiplicity of each irreducible representation $\Pi$ of $G$ in the 
regular representation $C^{\infty}(X)$ is always \textbf{finite}.
\vskip 0.05in
\item[{\rm{(2)}}]
Determine a condition on the pair $(G,X)$ for which
the multiplicity is \textbf{uniformly bounded} with respect to
all irreducible representations $\Pi$.
\end{enumerate}
\end{problem}

Since the condition (1) concerns individual finiteness with no
constraint on the dependence of irreducible representations $\Pi$,
we may think that the group $G$ has \lq\lq{stronger grip power}\rq\rq\ in (2).
The case that the group action of $G$ is transitive
is essential in Problem \ref{prob:6.2},
which we shall assume in the following.
Then Problem \ref{prob:6.2} is completely solved by 
Kobayashi--Oshima \cite{xktoshima}
when $G$ is a reductive algebraic Lie group.
To state the necessary and sufficient condition,
we prepare some terminology.

\begin{defin}
[Spherical variety and real spherical variety]\label{def:6.3}
\phantom{a}
\vskip 0.05in
\begin{enumerate}
\item[{\rm{(1)}}]
Suppose that a complex reductive Lie group $G_{\mathbb{C}}$ acts 
biholomorphically 
on a connected complex manifold $X_{\mathbb{C}}$.
We say $X_{\mathbb{C}}$ is {\bf{spherical}} or
{\bf{$G_{\mathbb{C}}$-spherical}}
if a Borel subgroup of $G_{\mathbb{C}}$ has 
an open orbit in $X_{\mathbb{C}}$.
\vskip 0.05in
\item[{\rm{(2)}}]
Suppose that a reductive Lie group $G$ acts continuously 
on a connected real manifold $X$.
We say $X$ is {\bf{real spherical}} or {\bf{$G$-real spherical}}
if a minimal parabolic subgroup of $G$ has an open orbit in $X$.
\end{enumerate}
\end{defin}

The terminology of \lq\lq{real sphericity}\rq\rq\ was introduced
by the author \cite{Ksuron} in the early 1990s 
for the study of Basic Problem \ref{prob:6.1}.
By definition, if $X_{\mathbb{C}}$ is $G_{\mathbb{C}}$-spherical,
then $X_{\mathbb{C}}$ is also $G_{\mathbb{C}}$-real spherical as a real manifold.

\begin{example}
\label{ex:spherical}
Let $X$ be a homogeneous space of a reductive algebraic Lie group $G$ and 
$X_{\mathbb{C}}$ its complexification.
\vskip 0.05in
\begin{enumerate}
\item[{\rm{(1)}}]
The following implications hold
(Aomoto, Wolf, and Kobayashi--Oshima \cite[Prop.\ 4.3]{xktoshima}).
\begin{align*}
&\textnormal{$X$ is a symmetric space.}\\
&\qquad \Downarrow \textnormal{Aomoto, Wolf}\\
&\textnormal{$X_{\mathbb{C}}$ is $G_{\mathbb{C}}$-spherical.}\\
&\qquad \Downarrow \textnormal{Kobayashi--Oshima} \\
&\textnormal{$X$ is $G$-real spherical.}\\
&\qquad \Uparrow \textnormal{obvious}\\
&\textnormal{$G$ is compact.}
\end{align*}
When $X$ admits a $G$-invariant Riemannian structure, then $X_\mathbb{C}$ 
is $G_\mathbb{C}$-spherical if and only if $X$ is
a weakly symmetric space in the sense of Selberg, see 
Vinberg \cite{Vinberg01} and Wolf \cite{Wolf07}.
\vskip 0.05in

\item[{\rm{(2)}}]
The irreducible symmetric spaces were classified by 
Berger \cite{ber} at the level of Lie algebras.
\vskip 0.05in

\item[{\rm{(3)}}]
The classification theory of spherical varieties $X_{\mathbb{C}}$
is given by Kr{\"a}mer \cite{xkramer}, 
Brion \cite{xbrion}, and
Mikityuk \cite{xmik}.
\vskip 0.05in

\item[{\rm{(4)}}]
The homogeneous space 
$(G \times G \times G)/\operatorname{diag} (G)$
is not a symmetric space.
It was determined by the author \cite[Ex.\ 2.8.6]{Ksuron} 
when it becomes real spherical in the study of 
multiplicities when decomposing tensor product representations
(Example \ref{ex:8.6}). 
A generalization of this 
will be described in Example \ref{ex:mfSZ} (2) in connection with 
branching problems for symmetric pairs.
\end{enumerate}
\end{example}

The solutions to Problem \ref{prob:6.2}, which is a 
reformalization of Basic Problem \ref{prob:6.1}, 
are given by the following two theorems.
For simplicity of the exposition, 
we assume $H$ to be reductive, see 
Remark \ref{rem:48} for more general cases.

\begin{theorem}
[Criterion for the finiteness of multiplicity]
\label{thm:HP}
Let $G$ be a reductive algebraic Lie group and
$H$ a reductive algebraic subgroup of $G$.
We set $X=G/H$.
Then the following two conditions on the pair $(G,H)$
are equivalent.
\vskip 0.05in
\begin{enumerate}
\item[{\rm{(i)}}]
{\rm{(representation theory)}}\enspace
$\dim_{\mathbb{C}} \operatorname{Hom}_{G}(\Pi, C^{\infty}(X))< \infty$
 $({}^{\forall} \Pi \in {\operatorname{Irr}}(G))$.  
\vskip 0.05in
\item[{\rm{(ii)}}]
{\rm{(geometry)}}\enspace
$X$ is $G$-real spherical.
\end{enumerate}
\end{theorem}

In \cite{xktoshima}, we have proved
not only a qualitative result 
(Theorem \ref{thm:HP})
but also quantitative results, namely, an upper estimate of the multiplicity
by using a boundary problem of partial
differential equations and a lower estimate by generalizing the 
classical Poisson transform 
\cite[Ch.\ II]{Helgason}, see \cite[Sect.\ 6.1]{xkProg2014}.
These estimates from the above and below yield
a criterion of the uniform boundedness of multiplicity as in
the following theorem, where the equivalence 
(ii) $\iff$ (iii) is classically known, and the main theme 
here is a connection with the representation theoretic property (i).

\begin{theorem}
[Criterion for the uniform boundedness of multiplicity {\cite{xktoshima}}]
\label{thm:HB}
Let $G$ be a reductive algebraic Lie group,
$H$ a reductive algebraic subgroup of $G$,
and $X=G/H$.
Then the following three conditions on the pair $(G,H)$
are equivalent.
\vskip 0.05in
\begin{enumerate}
\item[{\rm{(i)}}]
{\rm{(representation theory)}}\enspace
There exists a constant $C$ such that 
\begin{equation*}
\dim_{\mathbb{C}} \operatorname{Hom}_{G}(\Pi, C^{\infty}(X)) \le C
\quad
({}^{\forall} \Pi \in {\operatorname{Irr}}(G)).
\end{equation*}
\item
[{\rm{(ii)}}]
{\rm{(complex geometry)}}\enspace
The complexification $X_{\mathbb{C}}$ of $X$ is 
$G_{\mathbb{C}}$-spherical. 
\vskip 0.05in

\item[{\rm{(iii)}}]
{\rm{(ring theory)}}\enspace
The ring of $G$-invariant differential operators on $X$ is commutative.
\end{enumerate}
\end{theorem}

\begin{remark}\label{rem:48}
{\rm{
Theorems \ref{thm:HP} and \ref{thm:HB} give solutions 
to Problem \ref{prob:6.2} (1) and (2), respectively.
More generally, these theorems hold not only for 
the space  $C^{\infty}(X)$ of smooth functions
but also for the space of distributions/hyperfunctions and 
the space of sections of equivariant vector bundles.
Furthermore, one can drop
the assumption that the subgroup $H$ is reductive,
see \cite[Thms.\ A and B]{xktoshima} for precise formulation.
For instance, the theory of the Whittaker model 
(Kostant--Lynch, H.\ Matumoto) considers the case where
$H$ is a maximal unipotent subgroup $N$.
In this case, $G/N$ is always $G$-real spherical, 
and $G_{\mathbb{C}}/N_{\mathbb{C}}$ is $G_{\mathbb{C}}$-spherical
if and only if $G$ is quasi-split. Thus
Theorems \ref{thm:HP} and \ref{thm:HB} (in a generalized form) 
can be applied.
}}
\end{remark}

\begin{remark}
{\rm{
Theorem \ref{thm:HB}  includes
the new discovery that the property of the ``uniform boundedness of multiplicity''
is determined only by the complexification $(G_\mathbb{C}, X_\mathbb{C})$ 
and is independent of a real form $(G, X)$.
This observation suggests that an analogous result should 
hold more generally
for reductive algebraic groups over non-archimedean local fields as well.
Recently, Sakellaridis--Venkatesh \cite{SaVe17} has obtained 
some positive results in this direction.
See also \cite{Kobayashi21, Tauchi(a21)}
for an alternative approach to the proof 
(ii) $\Rightarrow$ (i) by using holonomic $\mathcal{D}$-modules.
}}
\end{remark}

Theorems \ref{thm:HP} and \ref{thm:HB} 
provide nice settings of global analysis in which 
the \lq\lq{grip strength}\rq\rq\ of representation theory 
is  \lq\lq{firm}\rq\rq\ on the space of functions.
The existing successful theory such as the Whittaker model 
and the analysis on semisimple symmetric spaces
(Section \ref{subsec:L2X}) mentioned above may be thought of as 
the global analysis in this framework (Example \ref{ex:spherical}).
There are also \lq\lq{new}\rq\rq\ settings suggested by 
Theorems \ref{thm:HP} and \ref{thm:HB}, to which the global analysis
has not been paid much attention, and as 
one of such settings, we shall discuss an application to
branching problems in Section \ref{subsec:6.2}.

\subsection{Spectrum of the regular representation $L^2(X)$:
a geometric criterion for temperedness --- new perspective, part 2}
\label{subsec:tempered}
In the previous section \ref{subsec:grip},
we focused on \lq\lq{multiplicity}\rq\rq\ from the perspective of 
the \lq\lq{grip strength}\rq\rq\ of a group on a function space
and proposed (real) sphericity
as
``a nice geometric framework for \textbf{detailed} study of global analysis''
beyond symmetric spaces.
On the other hand, even when
the ``grip strength'' of representation theory is not ``firm'',
we may still expect to analyze $L^2(X)$ in 
a \lq\lq{\textbf{coarse} standpoint}\rq\rq.
In this subsection,
including {\emph{non-spherical cases}},
let us focus on the support of the Plancherel measure
and consider the following problem.

\begin{basicproblem}
\label{prob:6.7}
Suppose that a reductive Lie group $G$ acts on a manifold $X$,
and that there is a $G$-invariant Radon measure on $X$. Determine
a necessary and sufficient condition
on a pair $(G,X)$ for which the regular representation
of $G$ on $L^2(X)$ 
is a tempered representation (Definition \ref{def:coarseLX} (2)).
\end{basicproblem}

Basic Problem \ref{prob:6.7} asks the condition
that any irreducible non-tempered representation 
(e.g.\ a complementary series representation) is not allowed to
contribute to the unitary representation $L^2(X)$.

\begin{observation}
\label{rem:6.15}
\phantom{a}
\begin{enumerate}
\item[{\rm{(1)}}]
In the case where $G/H$ is a semisimple symmetric space,
the Plancherel-type theorem is known
\cite{xdelorme, Oshima2002},
however,
it seems that 
a necessary and sufficient condition
on which $L^2(G/H)$ is tempered
had not been found until the general theory 
\cite{BK2015} is proved.
In fact, 
if one tried to apply the Plancherel-type formula to find an answer to 
Problem \ref{prob:6.7},
one would encounter a problem to find a precise
(non-)vanishing condition
of discrete series representations for $G/H$
with singular parameters, or equivalently, that of
certain cohomologies (Zuckerman derived functor modules)
after crossing a number of \lq\lq{walls}\rq\rq\
and such a condition is combinatorically
complicated in many cases,
see \cite[Chaps.\ 4 and 5]{Kobayashi92}
and \cite{Trapa01} for instance. 
\vskip 0.05in

\item[{\rm{(2)}}]
More generally, when $X_{\mathbb{C}}$ 
is not necessarily a spherical variety of $G_{\mathbb{C}}$,
the ring ${\mathbb{D}}_G(X)$ of $G$-invariant differential operators on $X$
is not commutative as seen in the equivalence (ii) $\iff$ (iii)
of Theorem \ref{thm:HB},
and so we cannot use effectively the existing method 
on non-commutative harmonic analysis based on an expansion
of the functions on $X$ into joint eigenfunctions 
with respect to the commutative ring ${\mathbb{D}}_G(X)$
as was the case of symmetric spaces.
\end{enumerate}
\end{observation}

As observed above,
to tackle Basic Problem \ref{prob:6.7},
we need to develop a completely new method itself.
For this,
we bring an idea of geometric group theory mentioned in Section \ref{sec:1}.
Let us start with an observation of an easy example.
If $H$ is a compact subgroup of $G$, then $L^2(G/H) \subset L^2(G)$ holds.
The following can be readily drawn from this observation.

\begin{example}
\label{ex:propertemp}
If the action of a group $G$ on $X$ is proper (Definition \ref{def:proper}),
then the regular representation in $L^2(X)$ is tempered.
\end{example}

Therefore, in the study of Basic Problem \ref{prob:6.7},
the non-trivial case is when the action of $G$ on $X$ is not proper.
Properness of the action is \textbf{qualitative property},
namely, a non-proper action means that
there exists a compact subset $S$ of $X$ such that 
the set 
$
\{g \in G: g S \cap S \ne \emptyset\}
$
is not compact.
In order to \textbf{quantify} this property,
we highlight the volume $\operatorname{vol}(g S \cap S)$
with respect to a Radon measure on $X$.
Viewed as a function of $g \in G$
\begin{equation}\label{eqn:vol(Kubo)}
G \ni g \mapsto \operatorname{vol}(gS \cap S) \in \mathbb{R}
\end{equation}
is a continuous function on $G$.
Definition \ref{def:proper} tells us that the $G$-action on $X$
is not proper if and only if the support of the function \eqref{eqn:vol(Kubo)}
is non-compact for some compact subset $S$ of $X$.
This suggests that the ``decay'' of the function \eqref{eqn:vol(Kubo)} at infinity
may be considered as  
a \lq\lq{quantification}\rq\rq\ of non-properness of the action.
By pursuing this idea, Basic Problem \ref{prob:6.7} is settled in 
Benoist--Kobayashi \cite{BK2015, BK2017}
when $X$ is a homogeneous space of a real reductive group $G$. 
To describe the solution, let us introduce a piecewise linear function
associated to a finite-dimensional representation of a Lie algebra.

\begin{defin}
[Moment of a representation {\cite[Sect.\ 2.3]{BK2017}}]
For a representation 
$
   \sigma \colon {\mathfrak{h}} \to {\operatorname{End}}_{\mathbb{R}}(V)
$
of  a Lie algebra ${\mathfrak{h}}$ on a finite-dimensional real vector space $V$,
the {\bf{moment}} $\rho_V$ is a function 
on ${\mathfrak{h}}$ defined as
\begin{align*}
  \rho_V \colon {\mathfrak{h}} \to {\mathbb{R}},&
\quad
 \text{$Y \mapsto$
 the sum of the absolute values of the real parts}\\
 &\text{\hspace{0.4in} of the eigenvalues of $\sigma(Y)$ on
 $V \otimes_{\mathbb{R}} {\mathbb{C}}$.}
\end{align*}
\end{defin}

The function $\rho_V$ is uniquely determined by the restriction to 
a maximal abelian split subalgebra  ${\mathfrak{a}}$ of $\mathfrak{h}$.
Further, the restriction $\rho_V|_{\mathfrak{a}}$ is piecewise linear,
in the sense that there exist finitely many convex polyhedral cones which cover
$\mathfrak{a}$ and on which $\rho_V$ is linear.
When $(\sigma,V)$ is the adjoint representation 
of a semisimple Lie algebra ${\mathfrak{h}}$,
the restriction $\rho_{\mathfrak{h}}|_{\mathfrak{a}}$ can be computed by using
a root system and coincides with a scalar multiple of
the usual ``$\rho$'' on the dominant
Weyl chamber.

With the use of the functions $\rho_V$ for the adjoint actions of 
$\mathfrak{h}$ on $V=\mathfrak{h}$ and $\mathfrak{g}/\mathfrak{h}$,
one can give a necessary and sufficient condition for Basic Problem \ref{prob:6.7}.

\begin{theorem}
[Criterion on the temperedness of $L^2(X)$]
\label{thm:BK}
Let 
$G$ be a real reductive Lie group and
$H$ a connected closed subgroup of $G$.
Also let 
${\mathfrak{g}}$ and ${\mathfrak{h}}$ be the Lie algebras of $G$ and $H$, respectively.
Then the following two conditions on a pair $(G,H)$ are equivalent.
\par\indent
{\rm{(i)}}
{\rm{(global analysis)}}\enspace
The regular representation $L^2(G/H)$ is tempered.
\par\indent
{\rm{(ii)}}
{\rm{(combinatorial geometry)}}\enspace
$\rho_{{\mathfrak{h}}} \le \rho_{{\mathfrak{g}}/{\mathfrak{h}}}$.  
\end{theorem}

The implication (i) $\Rightarrow$ (ii) follows from a local estimate of
the asymptotic behavior of the volume 
$\operatorname{vol}(g S \cap S)$, and 
the converse implication (ii) $\Rightarrow$ (i) is much more involved.
We note that Basic Problem \ref{prob:6.7} makes sense even when
there is no $G$-invariant Radon measure on $X$ 
by replacing $L^2(X)$ with the Hilbert space of $L^2$-sections 
for the half-density bundle over $X$. Theorem \ref{thm:BK} holds
in this generality.
See {\cite{BK2015, BK2017}} for details,
\cite{BK2021} for the classification of the pairs
$(G,H)$ of real reductive groups for which $L^2(G/H)$ is 
non-tempered, and
\cite{BK(a20)} for some connection with other disciplines such as
the orbit philosophy and the limit algebras.

\begin{remark}
\label{rem:6.14}
If 
$G$ is an algebraic group
and $X$ is an algebraic $G$-variety,
then, even when $X$ is not a homogeneous space of $G$,
one can give an answer to Basic Problem \ref{prob:6.7} by 
applying Theorem \ref{thm:BK} to generic $G$-orbits
\cite{BK2017}.
\end{remark}

\subsection{The sizes of a group $G$ and a manifold in global analysis}
\label{subsec:6.4}

Let us mention yet another new perspective
in global analysis via representation theory.
To give its flavor, we begin with
``coarse comparison'' of the ``size'' of 
the transformation group $G$ with that of the geometry $X$.
Suppose that a Lie group $G$ acts on a manifold $X$. If there
are at most finitely many $G$-orbits in $X$
(in particular, if $G$ acts transitively), 
one may regard that the size of $G$ is ``comparable'' with $X$
which we write symbolically as
\[
\text{group $G$ $\approx$ manifold $X$.}
\]

Homogeneous spaces $X=G/H$ are typical examples for the 
relation $G\approx X$. We may think of the main results of
Sections \ref{subsec:grip} and \ref{subsec:tempered} as 
a refinement of the 
``relation $G\approx X$'' by introducing
some kind of ``smallness of $X=G/H$ relative to $G$''
or ``largeness of $H$''
from the following points of view:
\smallskip

\begin{itemize}
\item Grip strength $\cdots$
multiplicity (Theorems \ref{thm:HP} and \ref{thm:HB}):
The ``larger'' $H$ is (or the ``smaller'' $X$ is), 
the better $G$ controls the function space on $X$.
\vskip 0.05in

\item Spectrum $\cdots$  temperedness criterion (Theorem \ref{thm:BK}):
The ``larger'' $H$ is (or the ``smaller'' $X$ is), the less likely the regular representation of $G$ on 
$L^2(X)$ becomes tempered.
\end{itemize}
\smallskip

As we have seen, the two notions of ``smallness of $X$'' 
are alike in appearance but quite different in nature 
among the case $G\approx X$.
In the rest of this section, we discuss
new directions of representation theory and global analysis
beyond the setting $G\approx X$:
\par\noindent
(1)\enspace
group $G$ \,\, $\gg$ \,\, manifold $X$ 
\quad
($G$ is \lq\lq{too large}\rq\rq\ to act on $X$ non-trivially.)
\par\noindent
(2)\enspace
group $G$ \,\, $\ll$ \,\, manifold $X$ 
\quad
(The dimensions of all $G$-orbits are smaller than \\
\hspace{2in}\;\;\,that of $X$.)

\vskip 0.5pc
{\bf{Global analysis of minimal representations}}
---
an example for the case that the size of $G$ $\gg$ the size of a manifold $X$ 

Let $e(G)$ be the smallest value of the codimensions of proper subgroups 
of a Lie group $G$. For example, if $G=GL(n,{\mathbb{R}})$, then $e(G)=n-1$.
This means that if the dimension of a manifold $X$ is less than $e(G)$, then 
any (continuous) action of $G$ on $X$ must be trivial.
In this way, if a Lie group $G$ is \lq\lq{too large}\rq\rq\ 
compared to a manifold $X$, for instance, if $e(G) > \dim X$,
then $G$ cannot act on $X$ geometrically.
We write group $G$ $\gg$ manifold $X$ in this case.
Even when $G \gg X$, 
we may perform global analysis on $X$ from a different perspective if
one can define
a natural representation of the group 
$G$ on the space of functions on $X$ although the representation
does not arise from a geometric
action of $G$ on $X$.
In such a case the \lq\lq{grip strength}\rq\rq\ of $G$ on the space of functions 
will be extremely firm, hence one may expect
that the representation theory plays a powerful role in the global analysis,
even more powerful than 
in the analysis on homogeneous spaces \cite{sato50}.
One of the examples is that $X$ is a Lagrangian submanifold of the minimal nilpotent coadjoint orbit (Section \ref{subsec:orbit}) and in this point of view,
``global analysis with minimal representations as a motif'' has emerged \cite{sato50}.
Global analysis in new directions based on 
the Schr{\"o}dinger model of a minimal representation beyond the 
Segel--Shale--Weil representation
\cite{xhkm, xkmanoAMS, KO3}
has been developed rapidly in recent years.
It includes
the construction and the unitarization of a minimal representation by the use of 
conformal geometry \cite{KO1, KO2}, 
the theory of special functions associated with fourth-order differential 
equations \cite{HKMM11, frenkel60}
and the deformation theory of the Fourier transform 
such as the $(k,a)$-generalized
Fourier transform \cite{BKO, CDBL18}.

\vskip 0.05in

{\bf{Visible action}}
---
an example for the case that the size of $G$ $\ll$ the size of a manifold $X$ 

We consider the opposite extremal case where 
group $G$ $\ll$ manifold $X$ in the sense that 
there is no open $G$-orbit $X$, and in particular, $G$ has continuously
many orbits in $X$.
Even in such a case,
there is still a possibility to 
develop global analysis with a reasonable \lq\lq{control}\rq\rq\ by group 
representations if one imposes further 
constraints such as geometrically defined differential equations.
The theory of visible actions on complex manifolds
and the propagation of multiplicity-free property \cite{xrims40, mfbundl}  
is a new attempt in the direction (see Section \ref{subsec:visible}).

\section{Program for Branching Problems in Representation Theory}
\label{sec:Stage}

In the latter part of this article we discuss
branching problems that ask how
the restrictions of irreducible representations 
behave when restricted to subgroups, 
e.g., finding their irreducible decompositions
(the \emph{branching laws}) 
for infinite-dimensional representations of Lie groups.
In spite of its potential importance,
the study of 
the restriction of irreducible representations
of a reductive Lie group 
to non-compact subgroups
was still underdeveloped in the 1980s,
except for some specific cases
such as the theta correspondence \cite{xhowe},
highest weight modules, the compact case,
or the $SL_2$-cases \cite{Mol, Re79}.
The main difficulty seemed to be the lack of 
promising perspectives for the general study of 
the restrictions of representations.
In fact, if one tries to find
a branching law for a group larger than $S L_2$,
then bad phenomena appear such as
infinite multiplicity in the branching laws, which we express as 
\lq\lq{the grip strength of the subgroup is not firm enough}\rq\rq.
For instance, 
when $n \ge 3$,
the tensor product of two principal series representations of 
$S L_n({\mathbb{R}})$ \lq\lq{contains}\rq\rq\ the same
irreducible representation with infinite multiplicity.
In the late 1980s, 
inspired by the new theory
of discontinuous groups beyond the Riemannian setting,
the author discovered a new example of ``good branching laws''
in the sense that an infinite-dimensional irreducible representation $\Pi$
decomposes discretely and multiplicity-freely when restricted
to \emph{non-compact} subgroups in the setting that $\Pi$ is \emph{not}
a highest weight module.\footnote{As its geometric background 
there was an open problem of spectral geometry proposed 
by Toshikazu Sunada at the time. For the details
see \cite{kks} on Sunada's conjecture and also \cite{kob09}.}
This type of branching law had not been previously known 
and an attempt to elucidate this example more generally became a trigger of 
 the general  study of the restriction of representations 
 which had been in a kind of a chaotic state.
Enough tools had been accumulated,
it was poised to take off.
In the following three steps
let us give an overview of some of the progress of the study of branching problems 
over around 30 years from the discovery:
\smallskip
\par
Stage A:\enspace 
Abstract feature of the restriction of representations (Section \ref{sec:A}).
\smallskip
\par
Stage B:\enspace
Branching laws (Section \ref{sec:B}).
\smallskip
\par
Stage C:\enspace
Construction of symmetry breaking operators (Section \ref{sec:C}).
\smallskip
\par\noindent
The name of each stage comes from their initials
\cite{xKVogan2015}.
We shall look into the roles of Stages A, B, and C in 
Sections \ref{subsec:10.1} and  \ref{subsec:SL2} later.
For more details of the program, see \cite{xKVogan2015}.

\section{Program for Branching Problems: Stage A}
\label{sec:A}

The aim of Stage A is to develop an abstract theory on the restriction of 
representations in the setting as general as possible.
Stage A will provide a bird's-eye view to the problem of 
the ``restriction of representations'' in which various
new directions of research may open up.
In particular,
this stage will play a role to single out a nice setting where
one could develop a detailed study of branching laws in Stages B and C.
For instance, in Stage A, we aim to construct a general theory to elucidate the following
properties.

\smallskip
\par\noindent
$\bullet$\enspace
(Existence of the continuous spectrum \cite{xkInvent94})
Does a continuous spectrum appear in the branching law
of the restriction $\Pi\vert_{G'}$ of a unitary representation $\Pi$
of a group $G$
to non-compact subgroups $G'$?
Or, does it decompose discretely?
(Section \ref{subsec:mfbranch}) 

\smallskip\smallskip
\par\noindent
$\bullet$\enspace
(Finiteness/boundedness of multiplicity \cite{Ksuron,xktoshima})
For irreducible representations of a group $G$,
is the multiplicity (i.e.\ the number of times that irreducible 
representations of a subgroup $G'$ appears)
finite or infinite?
In case each multiplicity is finite,
we may also ask, even strongly, if it has
uniform boundedness?
Even further, under what assumptions does a multiplicity-free theorem hold?
We may formulate these problems without assuming that 
the representations are unitary, see Section \ref{subsec:6.2}.

\smallskip\smallskip
\par\noindent
$\bullet$\enspace
(Support of branching laws \cite{BK2015})
Properties on irreducible representations that appear in a branching law.
For instance, 
as an analogue of temperedness criterion in
Section \ref{subsec:tempered}, one may ask if the restriction 
$\Pi\vert_{G'}$ is tempered as a representation of a subgroup $G'$ 
when $\Pi$ is a non-tempered 
irreducible unitary representations of a group $G$.

\subsection{Existence problem of the continuous spectrum in branching laws}
\label{subsec:mfbranch}

If $G'$ is a reductive Lie group, then
with the use of the notion of the direct integral (Section \ref{subsec:irrdeco}),
the restriction of $\Pi|_{G'}$ of a unitary representation $\Pi$
is decomposed uniquely into irreducible representations 
of the subgroup $G'$ as
\begin{equation}
\label{eqn:unibra}
   \Pi|_{G'}
   \simeq
  \int_{\widehat {G'}}^{\oplus} n_{\Pi}(\pi)\, \pi \, d \mu(\pi)
\end{equation}
where $\mu$ is a Borel measure on the unitary dual $\widehat{G'}$
(Theorem \ref{thm:4.1}).
As one of the problems in Stage A, first consider
whether a continuous spectrum appears in the branching law \eqref{eqn:unibra}.
When a continuous spectrum appears in a branching law,
analytic approaches will be natural for a detailed study of the restriction
$\Pi\vert_{G'}$. On the other hand, if one knows 
\emph{a priori} the branching law \eqref{eqn:unibra} is discrete,
one may study the restriction $\Pi\vert_{G'}$ also by
purely algebraic approach, and
expect to develop even more combinatorial techniques
and algorithms that compute branching laws.
Thus we address the following problem.

\begin{basicproblem}
[{\cite{xk:1}}]
\label{prob:deco}
Suppose that $\Pi$ is
an irreducible unitary representation of a group $G$
and that $G'$ is a subgroup of $G$.
Find a criterion on a triple $(G,G',\Pi)$ 
to decompose the restriction $\Pi|_{G'}$ into
the discrete direct sum of irreducible representations of $G'$.
Moreover, find a criterion that the multiplicity of each irreducible
representation is finite in the discrete branching law of the restriction
$\Pi\vert_{G'}$. 
\end{basicproblem}

In the latter case, we say the restriction $\Pi\vert_{G'}$ is
$G'$-\textbf{admissible} (\cite[Sect.~1]{xkInvent94}).
When $G'$ is a maximal compact subgroup $K$ of $G$,
the $K$-admissibility is nothing but Harish-Chandra's
admissibility (Definition \ref{def:adm}).

In the following let 
$G$ be a real reductive Lie group,
$\Pi$ an irreducible unitary representation of $G$,
and $G'$ a reductive subgroup of $G$.

Let us give three elementary examples of discretely decomposable
restrictions that we ask in Basic Problem \ref{prob:deco}.

\begin{example}
If $\Pi$ is finite-dimensional, 
then the restriction $\Pi|_{G'}$ 
is completely reducible, and hence decomposes discretely.
\end{example}

\begin{example}
If $G'$ is compact, 
then the restriction $\Pi|_{G'}$ decomposes discretely.
\end{example}

\begin{example}
[{\cite{xkAnn98, xkdisc}}]
\label{ex:8.2}
If $\Pi$ is a highest weight representation 
and also a pair $G \supset G'$ is of holomorphic type,
then the restriction $\Pi|_{G'}$ decomposes discretely.
\end{example}

Let us recall the terminology in Example \ref{ex:8.2}.
An irreducible representation $\Pi$ of a reductive Lie group $G$
is called a {\bf{highest weight representation}} if
its differential representation contains a nontrivial
subspace invariant under some Borel subalgebra
of the complex Lie algebra 
${\mathfrak{g}}_{\mathbb{C}}
={\mathfrak{g}}\otimes_{\mathbb{R}}{\mathbb{C}}$.
This subspace is automatically one-dimensional because 
$\Pi$ is irreducible.
All finite-dimensional irreducible representations 
are highest weight representations, whereas
highest weight representations are \lq\lq{rare}\rq\rq\ 
among irreducible infinite-dimensional 
representations. Moreover, a simple Lie group $G$
has an infinite-dimensional highest weight representation
if and only if $G$ is of Hermitian type, that is,
$X=G/K$ has the structure of a Hermitian symmetric space.
Here $K$ is the fixed point subgroup of a Cartan involution $\theta$
of $G$.
For a reductive Lie group $G$, we say $G$ is of Hermitian type if 
it is locally isomorphic to the direct product of a compact 
Lie group and simple Lie groups of Hermitian type.
Now, let $G'$ be a reductive subgroup of $G$.
Without loss of generality, we assume the Cartan involution
$\theta$ leaves $G'$ invariant. We set
$K':= G' \cap K$.
We say a pair $G \supset G'$ is 
{\bf{of holomorphic type}}
if both $G$ and $G'$ are of Hermitian type and if 
the natural embedding
$G'/K' \hookrightarrow G/K$ 
is holomorphic when we choose appropriate
$G'$-invariant and $G$-invariant complex structures on
$G'/K'$ and $G/K$, respectively.
See \cite{xkdisc, mf-korea} for 
the list of symmetric pairs $(G,G')$ of holomorphic type.

Highest weight representations 
and finite-dimensional representations are very much alike
and they
share many common properties, whereas 
irreducible infinite-dimensional representations that are not 
highest weight representations to which we refer as
\lq\lq{truly infinite-dimensional representations}\rq\rq,
do not. 
In contrast to Example \ref{ex:8.2} for highest weight representations,
experts seemed to have believed for a long time
that it would not be plausible for
\lq\lq{truly infinite-dimensional}\rq\rq\
representations
to decompose discretely
when they are restricted to a non-compact subgroup.
This \lq\lq{common sense}\rq\rq\ was reversed 
through
a study of discontinuous groups for the indefinite K{\" a}hler
manifold $X=SU(2,2)/U(1,2)$: 
the trigger was a discovery 
that any irreducible representation $\pi$ of $SU(2,2)$ in $L^2(X)$ 
is discretely decomposable 
with respect to the restriction to the subgroup $Sp(1,1) \simeq Spin(4,1)$
though $\pi$ is 
a \lq\lq{truly infinite-dimensional representation}\rq\rq\ 
(1988).
We refer to \cite{kob09} for the details.
The general theory of discretely decomposable 
branching laws has emerged in the 1990s
from the attempts that elucidate this particular example as a general principle.
The series of the three papers \cite{xkInvent94, xkAnn98, xkInvent98}
answer Basic Problem \ref{prob:deco} 
from perspectives in geometry, analysis, and algebra, respectively.
We provide a brief introduction to the main theorem in {\cite{xkAnn98}} here, 
which was obtained by techniques of microlocal analysis.
An alternative proof based on symplectic geometry is given in 
\cite{Kobayashi20}.

\begin{theorem}
[Criterion for the discrete decomposability of unitary representations]
\label{thm:2.2}
Let $\Pi$ be any irreducible unitary representation of a reductive Lie group $G$
and $G'$ a reductive subgroup of $G$.
Then the implication  {\rm{(ii)}} $\Rightarrow$ {\rm{(i)}} on a triple $(\Pi, G, G')$ holds.
\par\indent
{\rm{(i)}}\enspace
The restriction $\Pi|_{G'}$ 
is $G'$-admissible, namely, it decomposes discretely 
and also the multiplicity is finite.
\par\indent
{\rm{(ii)}}\enspace
$
\operatorname{AS}(\Pi)
\cap
\operatorname{Cone}(G')
=
\{0\}.
$
\end{theorem}

Here $\operatorname{AS}(\Pi)$ denotes the asymptotic cone of 
the $K$-types of the representation $\Pi$ \cite{KVcone}
and 
$\operatorname{Cone}(G')$ is the cone determined from the subgroup $G'$
(or more precisely, 
determined only by its maximal compact subgroup $K'$),
which is defined by the structure theory of Lie algebras  \cite{xkAnn98}
or by the moment map 
for the Hamiltonian $K$-action on $T^*(K/K')$
\cite[Thm~6.4.3]{deco-euro}.

\begin{remark}
\label{rem:3.4}
{\rm{
(1) When $G'=K$, one has $\operatorname{Cone}(K) = \{0\}$ 
hence Theorem \ref{thm:2.2} means Harish-Chandra's basic
theorem, see \cite[Thm.\ 3.4.10]{WaI} for instance, that 
$\Pi\vert_{K}$ is $K$-admissible for any irreducible 
unitary representation $\Pi$ of $G$.

(2) When $G'$ is compact, 
the equivalence (i) $\iff$ (ii) in
Theorem \ref{thm:2.2} holds \cite{deco-euro, Kobayashi20}.

(3) For a discrete series representation $\Pi$ of $G$,
the equivalence  (i) $\iff$ (ii) in Theorem \ref{thm:2.2} 
holds \cite{xdv, deco-euro, ZL10}.

(4) Even for nonunitary representations, one may formulate algebraically
the notion of \lq\lq{discrete decomposability}\rq\rq\ of the
restriction of representations \cite{xkInvent98}.
Then (ii) $\Rightarrow$ (i) in Theorem \ref{thm:2.2} still holds.
Moreover the necessary condition
for (algebraic) discrete decomposability 
is also given in \cite[Cor.\ 3.4]{xkInvent98} for the category of 
$({\mathfrak{g}}, K)$-modules
and in \cite[Thm.\ 4.1]{K12} for the category $\mathcal{O}$.
}}
\end{remark}

\vskip 0.5pc
\par\noindent
{\bf{\underline{Classification theory of discretely decomposable branching laws}}}:
For reductive symmetric pairs $(G,G')$,
the triple $(G,G',\Pi)$, for which the underlying $(\mathfrak{g},K)$-module
of an irreducible unitary representation $\Pi$ is discretely decomposable
when restricted to the subgroup $G'$
in the algebraic sense (Remark \ref{rem:3.4} (4)),
has been classified 
in the following settings
by carrying out the combinatorial computations
for the condition (ii) in Theorem \ref{thm:2.2} and for the criterion
by using the associated variety
in \cite{xkInvent98} detecting the opposite direction:
\vskip 0.05in
\begin{enumerate}
\item[$\bullet$]
$\Pi$ is a \lq\lq{geometric quantization}\rq\rq\ of an elliptic orbit
($\iff$ the underlying $({\mathfrak{g}},K)$-module $\Pi_K$
is a Zuckerman derived functor module 
$A_{\mathfrak{q}}(\lambda)$) \cite{decoAq},
\vskip 0.05in
\item[$\bullet$]
$\Pi$  is a \lq\lq{geometric quantization}\rq\rq\ of the minimal nilpotent orbit
($\iff$ $\Pi$ is a minimal representation) \cite{xkyosh13}, and
\vskip 0.05in
\item[$\bullet$]
the tensor product representation of two irreducible representations
 ($\iff$ $(G,G')$ is of the form $(H \times H, \operatorname{diag}(H))$)
\cite{xkyosh13}.
\end{enumerate}

These classification results for discretely decomposable restrictions
for symmetric pairs have been recently extended to 
\emph{non-symmetric pairs} $(G,G')$ by several authors,
see \cite{DGV17, He20, MSV19} for example.

\subsection{Multiplicity of branching laws}
\label{subsec:6.2}
Next let us discuss the multiplicity in branching laws.
In the following, let $G \supset G'$ be a pair of reductive Lie groups.
Here we allow representations to be nonunitary,
and $\Pi \in \operatorname{Irr}(G)$ and $\pi \in \operatorname{Irr}(G')$,
that is, $\Pi$ and $\pi$ be irreducible admissible representations of moderate
growth of $G$ and $G'$, respectively
(Definition \ref{def:irradm}).
We say
\begin{equation}
\label{eqn:multSBO}
   m(\Pi, \pi)
   := \dim_{\mathbb{C}} \invHom {G'}{\Pi|_{G'}}{\pi}
   \in {\mathbb{N}} \cup \{\infty\}
\end{equation}
is the {\bf{multiplicity}} of $\pi$ in the restriction $\Pi|_{G'}$.\footnote{
When $\Pi$ is a unitary representation,
the multiplicity $m(\Pi,\pi)$ for the spaces of $C^{\infty}$ vectors
could be larger than the multiplicity $n_{\Pi}(\pi)$ 
in the direct integral \eqref{eqn:unibra}.}
Unlike the case that $\Pi$ is finite-dimensional,
the second equation of \eqref{eqn:fmult} does not hold generally
when $\Pi$ is infinite dimensional.
Moreover,
the multiplicity $m(\Pi,\pi)$ could be infinite even in a natural situation
that $G'$ is a maximal subgroup of $G$.
For instance, the multiplicity on the tensor product of
infinite-dimensional irreducible representations
is often infinite, see Example \ref{ex:8.6}. 

In view of these observations, let us  illustrate schematically 
the hierarchy of the multiplicity in branching laws.
\begin{enumerate}
\item[]
a.\enspace \text{Some multiplicity is infinite.}
\vskip 0.05in
\item[]
b.\enspace \text{All the multiplicity is finite.} 
\vskip 0.05in
\item[]
c.\enspace \text{The multiplicity is uniformly bounded.}
\vskip 0.05in
\item[]
d.\enspace \text{multiplicity-free.}
\end{enumerate}
\vskip 0.05in
\noindent
The \lq\lq{grip strength}\rq\rq\ of the subgroup is larger in order of 
$a \prec b \prec c \prec d$ and
it is expected that one could do more detailed
analysis of branching laws accordingly.

By applying the criterion for
the \lq\lq{grip strength of representations in global analysis}\rq\rq\
on homogeneous spaces
(Theorem \ref{thm:HP}) discussed in Section \ref{subsec:grip}
to the homogeneous space $(G \times G')/\operatorname{diag}(G')$,
one obtains
a necessary and sufficient condition that
the multiplicity in the branching law is always finite
(\cite{Ksuron, xkProg2014, xktoshima}):

\begin{theorem}
[Finiteness of the multiplicity of branching laws]
\label{thm:pp}
The following two conditions on reductive Lie groups $G\supset G'$
are equivalent.
\par\indent
{\rm{(i)}}
{\rm{(representation theory)}}\enspace
For any $\Pi \in {\operatorname{Irr}}(G)$ and
$\pi \in {\operatorname{Irr}}(G')$,
$
   m(\Pi,\pi) < \infty.
$
\par\indent
{\rm{(ii)}}
{\rm{(geometry)}}\enspace
$(G \times G')/\operatorname{diag}(G')$ is 
a real spherical variety (Definition \ref{def:6.3}).
\end{theorem}

If one requires the ``uniform boundedness'' of
the multiplicity of branching laws, then
the following characterization holds {\cite{Ksuron, xkProg2014, xktoshima}}.

\begin{theorem}
[Uniform boundedness of multiplicity]
\label{thm:BB}
The following three conditions on a pair $(G, G')$ of reductive Lie groups
are equivalent.
\begin{enumerate}
\item[{\rm{(i)}}]
{\rm{(representation theory)}}\enspace
There exists a constant $C>0$  such that
\begin{equation*}
   m (\Pi, \pi) \le C
\quad
({}^{\forall} \Pi \in \operatorname{Irr}(G), 
{}^{\forall} \pi \in \operatorname{Irr}(G')).  
\end{equation*}

\item[{\rm{(ii)}}]
{\rm{(complex geometry)}}\enspace
$(G_{\mathbb{C}} \times G_{\mathbb{C}}')/\operatorname{diag}(G_{\mathbb{C}}')$
is a spherical variety.
\vskip 0.05in

\item[{\rm{(iii)}}]
{\rm{(ring theory)}}\enspace
The subalgebra $U({\mathfrak{g}})^{G'}$ of the universal enveloping algebra $U({\mathfrak{g}})$
of the Lie algebra  ${\mathfrak{g}}$  is commutative.
\end{enumerate}
\end{theorem}

As variants of Theorems \ref{thm:pp} and \ref{thm:BB},
one may fix $\Pi \in \operatorname{Irr}(G)$ and also ask a criterion
for the triple $(G,G',\Pi)$ with uniformly bounded multiplicity
$m(\Pi,\pi)$ where $\pi \in \operatorname{Irr}(G')$ varies
\cite[Prob.\ 6.2]{xKVogan2015}. See \cite[Thm.\ 7.6]{Kobayashi(a19)}
and \cite[Thm.\ 4.2]{Kobayashi21} for 
necessary and sufficient conditions
on the triple $(G,G',\Pi)$ with uniformly bounded multiplicity property
in the setting where $\Pi$ is $H$-distinguished for 
symmetric pairs $(G,H)$ or where $\Pi$ is a (degenerate)
principal series representation of $G$, respectively.

\begin{remark}
{\rm{
Similar to the criterion (Theorem \ref{thm:HB})
for uniform boundedness
in global analysis,
Theorem \ref{thm:BB} includes the discovery that 
the uniform boundedness of the multiplicity of branching laws
is determined only by the complexifications 
$({\mathfrak{g}_{\mathbb{C}}}, {\mathfrak{g}}_{\mathbb{C}}')$
of the Lie algebras and independent of the real forms.
This suggests that similar results may also hold 
for reductive algebraic groups over other local fields.
In fact, the assertion corresponding to (ii) $\Rightarrow$ (i) in 
Theorem \ref{thm:BB} (more strongly $C=1$) is shown
by Aizenbud--Gourevitch--Rallis--Schiffmann 
\cite{AGRS2010}
over a non-Archimedean local field.
}}
\end{remark}

\vskip 1pc
\par\noindent
{\bf{\underline{Classification theory of pairs of reductive Lie groups for which the branching}}}\\
{\bf{\underline{laws have always finite multiplicity}}}:
The easy-to-check geometric condition in Theorems \ref{thm:pp} and \ref{thm:BB}
allow us to extract settings in which branching laws behave nicely 
in terms of multiplicity.
\vskip 0.05in

\par\noindent
(1)\enspace
Since the criterion (ii) in
Theorem \ref{thm:BB}
is determined by the complexifications $(G_{\mathbb{C}}, G'_{\mathbb{C}})$,
one may check when $G$ is a 
compact group,
and thus the criterion (ii) in Theorem \ref{thm:BB} 
coincides with the one which already appeared
in finite-dimensional representation theory.
In fact, the complexified pairs $(G_{\mathbb{C}}, G_{\mathbb{C}}')$
satisfying (ii) were classified in the 1970s,
that is, such pairs $(G_{\mathbb{C}}, G_{\mathbb{C}}')$
are locally isomorphic to
the direct product of the pairs
$(SL(n, {\mathbb{C}}), GL({n-1}, {\mathbb{C}}))$, 
$(SO({n}, {\mathbb{C}}), SO({n-1}, {\mathbb{C}}))$, 
or pairs of abelian Lie algebras
(Kostant (unpublished) and Kr{\"a}mer \cite{xkramer}).
Example \ref{ex:mfSZ} on a pair of reductive Lie groups is their real forms.

\begin{example}
\label{ex:mfSZ}
(1) The constant $C$ in Theorem \ref{thm:BB} (i) can be taken to be $C=1$
for many of real forms
$(G, G')$ of $(SL(n,\mathbb{C}), GL(n-1,\mathbb{C}))$ or 
$(SO(n,\mathbb{C}), SO(n-1, \mathbb{C}))$ such as 
$(G, G')=(SL({n}, {\mathbb{R}}), GL({n-1}, {\mathbb{R}}))$, $(SO(p,q),SO(p-1,q))$,
see Aizenbud--Gourevitch \cite{AG2009} and Sun--Zhu \cite{xsunzhu}.
\par\noindent
(2)\enspace
The symmetric pairs $(G,G')$ that satisfy the criterion (ii) in Theorem \ref{thm:pp}
or in other words, for which $(G \times G')/\operatorname{diag} G'$ 
is real spherical were classified in 2013 (Kobayashi--Matsuki \cite{xKMt}).
This is also a generalization of the following 
earlier example in 1995 by the author \cite{Ksuron}, see also
\cite[Cor.\ 4.2]{xkProg2014}.

\end{example}

\begin{example}
\label{ex:8.6}
{\rm{
For a simple Lie group $G$,
the following four conditions are equivalent.
\par\noindent
(i)\enspace (invariant trilinear form)\enspace
For any $\pi_1$, $\pi_2$, $\pi_3 \in {\operatorname{Irr}}(G)$,
the space of invariant trilinear forms
$
\operatorname{Hom}_G(\pi_1 {\otimes} \pi_2 {\otimes} \pi_3,
                     {\mathbb{C}})
$
is finite-dimensional.
\par\noindent
(ii)\enspace(tensor product representation)\enspace
For any $\pi_1$, $\pi_2$, $\pi_3\in {\operatorname{Irr}}(G)$,
$
m(\pi_1 {\otimes} \pi_2, \pi_3)
<\infty.  
$
\par\noindent
(iii)\enspace(geometric condition)\enspace
$(G \times G \times G)/\operatorname{diag}(G)$ 
is real spherical.
\par\noindent
(iv)\enspace(classification)\enspace
${\mathfrak{g}} \simeq {\mathfrak {o}}(n,1)$
 ($n \ge 2$)
or
$G$ is compact.
}}
\end{example}

See also \cite{xKVogan2015, Kobayashi21} for 
finer classification results of the triples $(G,G',\Pi)$
rather than the pairs $(G,G')$ for the uniformly bounded multiplicity
restriction $\Pi\vert_{G'}$.

Example \ref{ex:8.6} suggests 
that invariant trilinear forms could be investigated explicitly
when $G=O(n,1)$. Indeed, invariant trilinear forms have been studied 
in detail for this group in recent years, not only algebraically but also analytically
(Stage C)
(for $n=2$,
Bernstein--Reznikov \cite{BR};
for general $n$,
Deitmar, Clerc--Kobayashi--{\O}rsted--Pevzner \cite{CKOP}, 
Clerc \cite{Clerc1617}, etc.).

In recent years, rapid progress in  \lq\lq{Stage C}\rq\rq\ 
for the branching problems has been made in 
the construction and classification problems of symmetry breaking operators 
\cite{Frahm, Juhl, KKP16, xtsbon, xksbonvec} in
the \lq\lq{good framework\rq\rq\ suggested by Theorem \ref{thm:pp} 
or more strongly by Theorem \ref{thm:BB}. 
Some of them interact with parabolic geometry such as conformal 
geometry and also with the theory of automorphic forms.
These topics will be discussed in Section \ref{subsec:10.3}.

\vskip 1pc
\par\noindent
{\bf{\underline{Case of infinite multiplicity}}}: 
When the multiplicity $m(\Pi,\pi) = \infty$ for irreducible representations $\Pi$ and $\pi$ of $G$ and its subgroup $G'$, respectively,
we have viewed 
\lq\lq{the grip strength of the subgroup $G'$ in the restriction 
$\Pi|_{G'}$ is weak}\rq\rq.
Apparently, the branching problems are 
\lq\lq{uncontrollable}\rq\rq.
However, even in this case, if we find an external algebraic structure to control
the (infinite-dimensional) space ${\operatorname{Hom}}_{G'}(\Pi|_{G'}, \pi)$ of symmetry breaking operators
(Section \ref{sec:C}), then,
 by using the structure as a good clue,
 it would be still possible
to investigate the restriction of representations.
A plausible candidate of such a structure is the algebra 
$U({\mathfrak{g}})^{G'}$ (see Theorem \ref{thm:BB} (iii)),
which acts on ${\operatorname{Hom}}_{G'}(\Pi|_{G'}, \pi)$ naturally.
Loosely speaking, the property $m(\Pi, \pi)=\infty$
can be described (in the level of the actions on the representation space) 
as follows.
\[
\text{$G'$ is relatively small.}
\,\,
\iff
\,\,
\text{The ring $U({\mathfrak{g}})^{G'}$ is large.}
\]

One of the motivations for Kitagawa's thesis
 \cite{Kitagawa} is to understand
${\operatorname{Hom}}_{G'}(\Pi|_{G'}, \pi)$
algebraically  as a $U({\mathfrak{g}})^{G'}$-module 
beyond the case  
where $\operatorname{Hom}_{G'}(\Pi\vert_{G'},\pi)$ is finite-dimensional.

\subsection{Visible actions and multiplicity-freeness}
\label{subsec:visible}

In Section \ref{subsec:6.2},
we discussed mainly the multiplicity $m(\Pi,\pi)$
for \underline{all} irreducible representations $\Pi$ of a group $G$
and for \underline{all} $\pi$ of its subgroup $G'$, and gave geometric
criteria for the pairs  $G \supset G'$  of groups that 
guarantee the finiteness of the multiplicity
(Theorem \ref{thm:pp}) and the uniform boundedness
(Theorem \ref{thm:BB}).
In this section we discuss an estimate of
the multiplicity of branching laws for the \underline{individual} 
irreducible representations $\Pi$ in more detail,
see also \cite[Problems 6.1 and 6.2]{xKVogan2015}.

\begin{basicproblem}
\label{prob:mf}
Classify triples $(G,G',\Pi)$
for which the restriction of an irreducible representation $\Pi$
of a group $G$ to its subgroup $G'$
is multiplicity-free.
\end{basicproblem}

The following two formulations are possible
for Basic Problem \ref{prob:mf} based on 
the definitions of the \lq\lq{multiplicity}\rq\rq.
\begin{enumerate}
\item[$\bullet$]
the case of unitary representations:
the multiplicity $n_{\Pi}(\pi)$ in the irreducible decomposition 
by using the direct integral, see \eqref{eqn:unibra}.
\smallskip
\item[$\bullet$]
the case of representations that are not necessarily unitary:
the multiplicity $m(\Pi,\pi)$ as the dimension of the space of 
symmetry breaking operators, see \eqref{eqn:multSBO}.
\end{enumerate}
In this subsection we consider the former case and 
introduce a new geometric principle that gives the multiplicity-freeness
of branching laws.

\begin{defin}
[{\cite[Def.\ 3.3.1]{xrims40}}]
\label{def:8.10}
{\rm{
Let $X$ be a connected complex manifold. A biholomorphic
action of a Lie group $G$ on $X$ is said to be 
\textbf{strongly visible} if
there exist a non-empty  $G$-invariant open set $X'$,
an anti-holomorphic diffeomorphism $\sigma$ of $X'$, and 
real submanifold $S$ (\emph{slice}) such that
\begin{equation*}
 \sigma |_S = \operatorname{id}
\quad
\text{and}
\quad
 G \cdot S= X'.
\end{equation*}
}}
\end{defin}

A strongly visible action is visible 
\cite[Thm.\ 4]{xrims40},
where we recall that the $G$-action on a complex manifold
is \textbf{visible} if there exist a non-empty
$G'$-invariant open set $X'$ and a totally real submanifold
$S$ in $X$ such that $G\cdot S = X'$ and 
$J_x(T_x S) \subset T_x(G\cdot x)$ for all $x \in X$ 
\cite[Def.\ 2.3]{xk2004}.

On the multiplicity-freeness of representations,
the strong visibility gives a new \lq\lq{mechanism}\rq\rq\
that produces systematically
from simple examples 
(for instance, one-dimensional representations, which are clearly multiplicity-free)
to more complicated examples
 (for instance, multiplicity-free infinite-dimensional representations).
This mechanism is formulated as a 
\lq\lq{propagation theorem of multiplicity-freeness}\rq\rq,
for which we describe in a slightly simplified way by omitting
some small technical compatibility conditions about the isotropy action on fibers
(that are automatically satisfied in many situations)
so the next theorem mentions only the main assumptions and conclusions.
For precise statements, see  \cite{mfbundl}.

\begin{theorem}
[Propagation theorem of multiplicity-freeness]
\label{thm:mfbundl}
Suppose that a group $G$ acts on a holomorphic vector bundle
${\mathcal{V}}$ over a complex manifold $X$
and that the action of $G$ on the base space $X$ is strongly visible.
If the isotropy representation of the isotropy subgroup at a
generic point on the fiber is 
multiplicity-free, then any unitary representation $\Pi$ of $G$
realized on subspaces of the space ${\mathcal{O}}(X,{\mathcal{V}})$ of 
holomorphic sections is multiplicity-free.
\end{theorem}

Here are two examples of the applications.

\begin{example}
[Highest weight representations]
\label{ex:MFholo}
{\rm{
Suppose that $(G,G')$ is a symmetric pair
and $\Pi$ is a unitary highest weight representation of $G$.
If the minimal $K$-type of $\Pi$ is one-dimensional,
then the irreducible decomposition of the restriction $\Pi|_{G'}$
is multiplicity-free (Kobayashi \cite{xrims40}).
This theorem can be derived from
the fact that the subgroup $G'$ acts on the Hermitian symmetric space 
$G/K$ strongly visibly \cite{K2007a} and from Theorem \ref{thm:mfbundl}.
The restriction $\Pi\vert_{G'}$ may or may not contain 
continuous spectrum, cf \cite[Sect.\ 5]{xkInvent98}.
See Theorem \ref{thm:HKSK} for an explicit branching law
in the discretely decomposable setting.
}}
\end{example}

\begin{example}
[Tensor product representations]
\label{example:6.16}
{\rm{
Consider the pairs $(\pi_1, \pi_2)$ of irreducible finite-dimensional representations
of the unitary group $U(n)$, for which the tensor product representation
$\pi_1 \otimes \pi_2$ decomposes multiplicity-freely.
Such pairs $(\pi_1, \pi_2)$ include the classical cases
when $\pi_2 \simeq \Exterior^j(\mathbb{C}^n)$ or 
$S^j(\mathbb{C}^n)$ (Pieri's rule) and were classified by
a combinatorial method (Stembridge \cite{Stem2001}, 2001).
All such pairs $(\pi_1, \pi_2)$
can be reconstructed geometrically
from a $(G\times G)$-equivariant holomorphic vector bundle over 
a double flag variety $X=G/P_1 \times G/P_2$ with strongly visible
action of $G$ via the diagonal action (Kobayashi {\cite{xk2004}}).
This geometric interpretation based on the theory of visible actions 
is extended to a reconstruction of all multiplicity-free tensor product representations of $SO(n)$ (Tanaka \cite{Ty2015}).
}}
\end{example}

\vskip 0.6pc
\par\noindent
{\bf{\underline{Classification theory of visible actions}}}:
Whereas (strongly) visible actions are defined on \underline{complex manifolds},
there are analogous notions in other geometries: 
polar actions of isometry groups on
\underline{Riemannian manifolds} and 
coisotropic actions on
\underline{symplectic manifolds}
(Guillemin--Sternberg and Huckleberry--Wurzbacher \cite{HW90}).
For K{\"a}hler manifolds on which complex, symplectic and Riemannian
structures are defined in a compatible fashion, 
these three notions are close to each other to some extent,
see Podest{\`a}--Thorbergsson \cite{PT02} and the author
\cite[Thms.\ 7 and 8]{xrims40} for precise formulation.
For polar actions of compact groups on Riemannian manifolds,
the classification theory has been developed over decades.
On the other hand, the classification theory for
visible actions has started only recently
(see Kobayashi \cite{K2007a, K2007b},
Sasaki \cite{Sa2008},
Tanaka \cite{Ty2015}, 
and the references therein),
and it
would be interesting to pursue further developments
and discoveries in this \lq\lq{young}\rq\rq\ area.
Whereas the classification theory of polar actions has focused mainly
on the compact setting from topological viewpoints historically,
the classification theory of visible actions will be useful
also for \underline{non-compact} transformation groups $G$
because it will yield new families of 
\underline{infinite-dimensional} multiplicity-free representations
of non-compact Lie groups.
We note that for a compact Lie group $G$, strongly visible action of the group $G$ 
is essentially equivalent to sphericity of the complexified group $G_{\mathbb{C}}$,
see Tanaka \cite{Ty2}, and thus
the classification theory
of strongly visible actions of a \textbf{compact} $G$ is essentially the same with
that of spherical varieties for the complexified reductive group $G_{\mathbb{C}}$.
In the expository paper \cite{xrims40}, we have presented various
classification results of visible actions including those of 
\textbf{non-compact} Lie groups, which have led us, via 
Theorem \ref{thm:mfbundl},
to new multiplicity-free theorems
as well as
a unified and geometric 
proof of the existing multiplicity-free theorems
for some families of representations 
that were found in the past by individual arguments.

\section{Branching Laws: Stage B}
\label{sec:B}
In Stage B of the program, we aim to find explicitly 
the irreducible decomposition of the restriction 
of representations (\textbf{branching law}). 
Here Stage A (Section \ref{sec:A}) serves
as a guideline to single out a \lq\lq{nice setting}\rq\rq\ in which
we could expect a simple and detailed study of branching laws 
through
a priori estimate.
In this section we describe typical examples for 
branching laws with focus on the following two 
\lq\lq{nice settings}\rq\rq:
multiplicity-free cases (Section \ref{subsec:Bmf})
and discretely decomposable cases 
(Section \ref{subsec:Bdiscdeco}).

\subsection{Multiplicity-free representations}
\label{subsec:Bmf}

Multiplicity-free representations are often hidden
in classical analysis, even though we usually do not notice
that there are representations behind.
For example,
the Fourier expansion, Taylor expansion, and
spherical harmonics expansion were already useful tools
in analysis, historically much before the notion of groups
and representations emerged. From representation-theoretic 
viewpoints, these expansions may
be regarded as irreducible decompositions
of multiplicity-free representations.
We observe here the
\lq\lq{multiplicity-freeness}\rq\rq\ 
is the underlying algebraic structure of 
these expansions in the sense that multiplicity-freeness
assures that the irreducible decomposition is canonical.
The expansion by the Gelfand--Tsetlin basis is 
also defined by this principle. More generally, we may
utilize multiplicity-free representations as a driving force
not only for studying 
the explicit formulas of branching laws (Stage B) 
but also for pursuing global analysis 
through canonical expansions of functions via
representation theory (Stage C).

In this subsection 
we consider the setting of Example \ref{ex:MFholo} as a case
in which we know in advance that 
the branching law is multiplicity-free.
Moreover, we assume that $(G,G')$ is of holomorphic type,
by which we know a priori
that the branching law is further discretely decomposable
(Example \ref{ex:8.2}).
In this case the explicit formulas of branching laws (Stage B) 
should take a simple form, which we describe now.

First we set up some notation.
Let 
$G$ be a real simple Lie group of Hermitian type
(Example \ref{ex:8.2}) 
and $(G,G')$ a symmetric pair of holomorphic type
defined by an involution $\sigma$ of $G$.
For simplicity, we take $G'$ to be
the connected component $G^\sigma_0$
containing the identity element
of $G^{\sigma}:=\{g\in G:\sigma g=g\}$.
Take a Cartan involution $\theta$ of $G$ commuting with $\sigma$
and we write $K$ for the fixed point group of the Cartan involution $\theta$.
Then the complexification ${\mathfrak {g}}_{\mathbb{C}}$ 
of the Lie algebra of $G$
can be decomposed into the direct sum
${\mathfrak {g}}_{\mathbb{C}}=
{\mathfrak {k}}_{\mathbb{C}}+{\mathfrak {p}}_+
+ {\mathfrak {p}}_-$
as a $K$-module, and $G/K$ carries a structure of a Hermitian symmetric
space with holomorphic tangent space $T_{eK}(G/K) \simeq \mathfrak{p}_+$.
We take a Cartan subalgebra $\mathfrak{j}$ of $\mathfrak{k}$ such that 
$\mathfrak{j}^\sigma$ is also a Cartan subalgebra of 
$\mathfrak{k}^\sigma=\{X\in \mathfrak{k}: \sigma X= X\}$.
We fix compatible positive systems 
$\Delta^+(\mathfrak{k}_\mathbb{C}, \mathfrak{j}_\mathbb{C})$
and 
$\Delta^+(\mathfrak{k}^\sigma_\mathbb{C}, \mathfrak{j}^\sigma_\mathbb{C})$.
Since $\sigma\theta = \theta \sigma$,
we have $(\sigma\theta)^2=\text{id}$, 
hence $\mathfrak{g}^{\sigma \theta}_\mathbb{C}$ is a
reductive Lie algebra.
We decompose  ${\mathfrak {g}}_{\mathbb{C}}^{\sigma \theta}$
into the direct sum of the simple Lie algebras ${\mathfrak {g}}_{\mathbb{C}}^{(i)}$
$(1 \leq i \leq N)$ and the abelian ideal ${\mathfrak {g}}_{\mathbb{C}}^{(0)}$.
For each $i$ $(\ne 0)$, we write $\{\nu_1^{(i)}, \ldots, \nu_{k_i}^{(i)}\}$
for a maximal set of strongly orthogonal roots of
$\Delta({\mathfrak {p}}_+^{-\sigma} \cap {\mathfrak {g}}_{\mathbb{C}}^{(i)})
(\subset ({\mathfrak {j}}_{\mathbb{C}}^{\sigma})^{\ast})$
with $\nu^{(i)}_k$ the lowest among the elements in 
$\Delta({\mathfrak {p}}_+^{-\sigma} \cap {\mathfrak {g}}_{\mathbb{C}}^{(i)})$
that are strongly orthogonal to $\nu^{(i)}_1, \ldots, \nu^{(i)}_{k-1}$.

We parametrize holomorphic discrete series representations as follows. 
Any holomorphic discrete series representation $\Pi$ of $G$ 
is determined by its unique minimal $K$-type.
We denote $\Pi$ by  $\Pi^G(\lambda)$ if
$\lambda \in \mathfrak{j}_\mathbb{C}^*$ is 
the highest weight of the minimal $K$-type of $\Pi$.
Similarly, 
holomorphic discrete series representations of the subgroup 
$G'=G^\sigma_0$ are expressed as $\pi^{G'}(\mu) \in \widehat {G'}$ 
if $\mu \in ({\mathfrak {j}}_{\mathbb{C}}^{\sigma})^{\ast}$
is the highest weight of the minimal $K'$-type with respect to 
$\Delta^+(\mathfrak{k}^\sigma_\mathbb{C}, \mathfrak{j}^\sigma_\mathbb{C})$,
where $K'=K\cap G' = K^\sigma_0$.

\begin{theorem}
[Hua--Kostant--Schmid--Kobayashi]
\label{thm:HKSK}
Let $(G,G')$ be a symmetric pair of holomorphic type.
For any holomorphic discrete series representation $\Pi^G(\lambda)$ of scalar type,
the following multiplicity-free decomposition holds.
\begin{equation}
\label{eqn:HKSK}
  \Pi^{G}(\lambda)|_{G'}
\simeq 
\bigoplus_{i=1}^N 
{\sum_{a_1^{(i)} \ge \cdots \ge a_{k_i}^{(i)} \ge 0}}^{\hskip-1.8pc\oplus}
\pi^{G'} (\lambda|_{{\mathfrak {j}}^{\sigma}}
-
\sum_{j=1}^{k_i} a_j^{(i)} \nu_j^{(i)})
\qquad
\emph{(Hilbert direct sum)}.
\end{equation}
\end{theorem}

\begin{remark}
{\rm{
If $G'$ is a maximal compact subgroup of $G$, 
then $\mathfrak{g}^{(0)}_\mathbb{C} = \{0\}$, $N=1$,
$\mathfrak{g}^{\sigma \theta} = \mathfrak{g}$,
and any irreducible summand
$\pi^{G'}(\mu)$ is finite-dimensional. In this case, the
formula \eqref{eqn:HKSK} is due to 
Hua (classical groups), Kostant (unpublished), and
Schmid \cite{xschthe}.
The proof for the branching law \eqref{eqn:HKSK} 
in the general setting where $G'$ is non-compact 
can be found in {\cite[Thm.\ 8.3]{mf-korea}}.
}}
\end{remark}

\begin{remark}
The \lq\lq{geometric quantization}\rq\rq\
in Section \ref{subsec:orbit} commutes with the reduction
in this case, and
the \lq\lq{classical limit}\rq\rq\ 
of the branching law \eqref{eqn:HKSK}
is given by the Corwin--Greenleaf function
for coadjoint orbits \cite{KN03, xKN, Paradan15}.
\end{remark}

Yet another important multiplicity-free results in the branching laws
are for more special pairs $(G,G')$, 
but for more general representations. 
Typical cases are those which we have seen
in Section \ref{subsec:6.2}, that is, 
the real forms $(G, G')$ of 
$(G L (n, {\mathbb{C}}), G L({n-1}, {\mathbb{C}}))$
or $(S O(n, {\mathbb{C}}), S O({n-1}, {\mathbb{C}}))$,
such as $(GL(n,\mathbb{R}), GL(n-1,\mathbb{R})$ or
$(SO(p,q), SO(p-1,q))$, have the property that the multiplicity
$m(\Pi,\pi)=\dim_\mathbb{C}\operatorname{Hom}_{G'}(\Pi\vert_{G'}, \pi)$ 
is either 0 or 1 for any $\Pi \in \operatorname{Irr}(G)$
and $\pi \in \operatorname{Irr}(G')$ even when $\Pi$ is neither
a unitarizable representation nor a highest weight representation
(Example \ref{ex:mfSZ}).
For these pairs,
when $\Pi$ and $\pi$ are both tempered representations,
the description of $m(\Pi,\pi)$ is predicted by the (local)
Gan--Gross--Prasad conjecture \cite{GP, xgrwaII}. 
In Section \ref{subsec:10.3},
we shall consider the case in which
some aspect of this conjecture 
(and more generally, non-tempered cases)
is connected with a particular problem in 
conformal geometry.

\subsection{Discretely decomposable branching laws}
\label{subsec:Bdiscdeco}
When 
the restriction is discretely decomposable,
it is expected that 
algebraic approaches would be useful in finding
explicit branching laws.
The first general theory of discretely decomposable restriction
was established in 1990s
in \cite{xkInvent94, xkAnn98, xkInvent98},
see e.g.\ Theorem \ref{thm:2.2}
for the criterion
of the triples $(G,G',\Pi)$ with discretely 
decomposable restriction $\Pi\vert_{G'}$.
In parallel,
there have been various attempts to find concrete
branching laws
in the framework of discretely decomposable restriction.
Among them,
especially important is 
the case that $\Pi$ is the geometric quantization of an 
elliptic orbit (see Section \ref{subsec:orbit}),
namely,
the case that the $({\mathfrak{g}}, K)$-module $\Pi_K$
is a Zuckerman derived functor module $A_{\mathfrak {q}}(\lambda)$,
where $\mathfrak{q}$ is a $\theta$-stable parabolic subalgebra,
see Knapp--Vogan \cite{xknappv} for instance, 
for the definition of $A_{\mathfrak {q}}(\lambda)$.
The first explicit branching laws 
for $\Pi_K = A_{\mathfrak{q}}(\lambda)$ 
in the general setting where $\Pi_K$ do not have highest weights
were given for any of the
adjacent pairs $(G,G')$ in the following diagram
(Kobayashi \cite{xk:1, xkInvent94}).
\allowdisplaybreaks[0]
\begin{alignat*}{4}
& O(4p, 4q) && \supset \,\, &&O(4k) && \times O(4p-4k, 4q)
\\
&  \qquad\cup     &&         &&      && \quad\cup
\\
& U(2p, 2q) && \supset &&U(2k) && \times U(2p-2k, 2q)
\\
&  \qquad\cup     &&         &&      && \quad\cup
\\
& Sp(p, q) && \supset && Sp(k) && \times Sp(p-k, q)
\end{alignat*}

The method in obtaining the branching laws
\cite{xk:1, xkInvent94} was to use a structural 
theorem of the ring of
invariant differential operators on homogeneous spaces $X$
with \lq\lq{overgroups}\rq\rq\ and to realize
$A_{\mathfrak{q}}(\lambda)$ in the space of functions on $X$,
see \cite{kob09} for a related geometric problem and \cite{KK19, K16b}
for the method in full generality.
Ever since \cite{xk:1,xkInvent94},
various methods have been developed for finding
explicit formulas of discretely decomposable branching laws
(Stage B) in the other settings by
Gross--Wallach \cite{xgrwaII}, 
Loke, 
J.-S. Li, 
Huang--Pand{\v z}i{\'c}--Savin, 
{\O}rsted--Speh \cite{xspeh}, 
Duflo--Vargas \cite{xdv},
Sekiguchi \cite{xhsekiijm},
Y.\ Oshima \cite{yophd},
Kobayashi \cite{mf-korea, K16b}, and so on.
Among them,
 Duflo--Vargas \cite{xdv} develops
 the idea of the orbit method and symplectic geometry
 when $\Pi$ is a discrete series representation,
 and Y.\ Oshima \cite{yophd} makes 
 use of the theory of ${\mathcal{D}}$-modules
 for Zuckerman derived functor modules 
 $\Pi_K= A_{\mathfrak{q}}(\lambda)$.
 
\section{Program for the Theory of Branching Laws: Stage C}
\label{sec:C}
In Stage C, we consider not only abstract branching laws 
(the decomposition of representations) but also
how they decompose (the decomposition of vectors).
For the latter purpose, a crucial step is to construct
$G'$-intertwining operators
from $\Pi$ to $\pi$ ({\bf{symmetry breaking operators}})
in a geometric model of 
irreducible representations $\Pi$ and $\pi$ 
of a group $G$ and its subgroup $G'$, respectively.
One can also consider $G'$-intertwining operators
in the opposite direction, namely, from $\pi$ to $\Pi$
({\bf{holographic operators}} \cite{Kobayashi(a19), KP3}) as a dual notion,
but we do not discuss holographic operators in this article.
Let us start with an elementary example.

\subsection{The regular representation of ${\mathbb{R}}$ and Fourier transform}
\label{subsec:10.1}
By using the regular representation $L^2({\mathbb{R}})$ of the additive group ${\mathbb{R}}$,
we illustrate Stages A, B, and C in Section  \ref{sec:Stage}.

\par\noindent
{\bf{Stage A.}}\enspace
The irreducible decomposition of
the regular representation 
$L^2({\mathbb{R}})$ is multiplicity-free and has only
continuous spectrum.
\smallskip
\par\noindent
{\bf{Stage B.}}\enspace
The regular representation $L^2({\mathbb{R}})$
is decomposed into the direct integral 
(Theorem \ref{thm:4.1})
of one-dimensional Hilbert spaces 
${\mathbb{C}} e^{ix\xi}$.
\smallskip
\par\noindent
{\bf{Stage C.}}\enspace
The irreducible decomposition of $L^2({\mathbb{R}})$ is 
realized by the Fourier transform \eqref{eqn:F} concretely.
\smallskip

Apparently, this example on the classical harmonic analysis 
looks nothing to do with branching laws,
however, one may interpret it as an example of restriction problems.
For example, let $\Pi$ be a unitary principal series representation of 
$G=SL(2, {\mathbb{R}})$, and $N(\simeq \mathbb{R})$
a maximal unipotent subgroup. Then the restriction $\Pi\vert_{N}$
is unitarily equivalent to the regular representation 
$L^2({\mathbb{R}})$ of the additive group ${\mathbb{R}}$
(see, for instance, \cite[Prop.~3.3.2]{deco-euro}).
Thus the above case may be interpreted as Stages A--C
for the restriction problem of $SL(2, {\mathbb{R}}) \downarrow {\mathbb{R}}$.
Needless to say, since 
we already know the theory of the Fourier transform on $L^2(\mathbb{R})$
which corresponds to Stage C for this particular example $\Pi\vert_{N}$,
Stages A and B stay in the background.

\subsection{Tensor product representation of $SL(2, {\mathbb{R}})$}
\label{subsec:SL2}
Next, 
we illustrate Stages A--C of branching laws by 
another example where the groups are 
\lq\lq{highly non-commutative}\rq\rq\
this time rather than abelian subgroups as in the previous case.

Let ${\mathcal{O}}({\mathcal{H}})$ denote
the space of holomorphic functions on the upper half plane
$
 {\mathcal{H}}
  =
 \{z \in {\mathbb{C}}: \operatorname{Im} z >0\}.
$
Then, there are a family of linear actions of 
the group $G=SL(2, {\mathbb{R}})$ 
on ${\mathcal{O}}({\mathcal{H}})$ 
with integer parameters $\lambda \in {\mathbb{Z}}$
as follows.
\[
   (\pi_{\lambda}\begin{pmatrix} a & b \\ c & d \end{pmatrix} f)(z)
  :=
  (-cz+a)^{-\lambda} f\left(\frac{dz-b}{-cz+a}\right).
\]
Moreover, if $\lambda >1$, then 
$V_{\lambda}
:={\mathcal{O}}(\mathcal{H}) \cap L^2({\mathcal{H}}, y^{\lambda-2}d x d y)$
is an infinite-dimensional Hilbert space and the action
$\pi_{\lambda}$ on $V_{\lambda}$
defines an irreducible and unitary representation of $G=SL(2, {\mathbb{R}})$,
referred to as the holomorphic discrete series representation.
We consider the tensor product $\pi_{\lambda'} \otimes \pi_{\lambda''}$
of two such representations, which may be interpreted as an example of 
the restriction from the direct product group $G \times G$ to the 
diagonal subgroup $\operatorname{diag}(G) \simeq G$.

\begin{enumerate}
\item[$\bullet$]
In Stage A,
we have the following a priori estimate on \lq\lq{abstract features}\rq\rq\ 
of the branching law:

For $\lambda', \lambda''>1$, the tensor product representation 
$\pi_{\lambda'} \otimes \pi_{\lambda''}$
decomposes discretely and multiplicity-freely into
irreducible unitary representations of $G$.
Here the symbol for the tensor product 
$\otimes$ of two Hilbert spaces means taking the Hilbert 
completion of the algebraic tensor product.
(Each property \lq\lq{discrete decomposability}\rq\rq\ 
and \lq\lq{multiplicity-freeness}\rq\rq\
of the tensor product representation
$\pi_{\lambda'}\otimes \pi_{\lambda''}$ is a special 
case of the general results of the restriction in Stage A as we have seen in 
Theorems \ref{thm:2.2} and \ref{thm:mfbundl}, respectively.)
\vskip 0.05in

\item[$\bullet$]
In Stage B, we determine a concrete branching law:

If $\lambda', \lambda'' >1$, then
one has the following irreducible decomposition
\begin{equation}
\label{eqn:Repka}
\pi_{\lambda'} \otimes \pi_{\lambda''}
\simeq
{\sum_{a \in {\mathbb{N}}}}^{\oplus}
\pi_{\lambda'+ \lambda''+2a}
\qquad
\text{(Hilbert direct sum).}
\end{equation}
The formula \eqref{eqn:Repka} in 
the $SL(2,\mathbb{R})$ case is due to
Molchanov \cite{Mol} and Repka \cite{Re79}.
Theorem \ref{thm:HKSK} is a generalization of the formula
\eqref{eqn:Repka} to arbitrary 
semisimple symmetric pairs of holomorphic type.
\vskip 0.05in

\item[$\bullet$]
In Stage C, we construct {\bf{symmetry breaking operators}}.

For given $\lambda'$, $\lambda''$, $\lambda''' \in {\mathbb{Z}}$,
a linear map
$
   R \colon {\mathcal{O}}({\mathcal{H}}) \otimes {\mathcal{O}}({\mathcal{H}})
            \to {\mathcal{O}}({\mathcal{H}})
$
satisfying
\begin{equation}
\label{eqn:Rinv}
  R(\pi_{\lambda'}(g) f_1 \otimes \pi_{\lambda''}(g)f_2)
 =\pi_{\lambda'''}(g) R(f_1 \otimes f_2)
\quad
 ({}^{\forall} g \in G)
\end{equation}
is a symmetry breaking operator
from the tensor product representation 
$\pi_{\lambda'} \otimes \pi_{\lambda''}$
to $\pi_{\lambda'''}$ with respect to the restriction $G\times G \downarrow G$ 
where $G=SL(2,\mathbb{R})$.
The next theorem interprets the classical Rankin--Cohen
bidifferential operator \cite{xcmz97, xrankin},
which was originally used to construct
modular forms of higher weight
from those of lower weight, 
as a symmetry breaking operator in
branching problems of
representation theory.
\end{enumerate}

\begin{theorem}
[Rankin--Cohen bidifferential operator
\cite{xcmz97, KP2, xrankin}]
\label{thm:RC}
Suppose that $\lambda'''-\lambda'-\lambda''$ 
is a non-negative even integer.
If we set $2a:=\lambda'''-\lambda'-\lambda''$ 
($a$ is a natural number),
then the linear map 
$
RC_{\lambda',\lambda''}^{\lambda'''}
\colon
{\mathcal{O}}({\mathcal{H}})
\otimes 
{\mathcal{O}}({\mathcal{H}})
\to
{\mathcal{O}}({\mathcal{H}})
$
defined by
\begin{equation}
\label{eqn:RC}
RC_{\lambda',\lambda''}^{\lambda'''}(f_1 \otimes f_2)(z)
:=
\sum_{\ell=0}^a
\frac{(-1)^\ell \Gamma(\lambda'+a) \Gamma(\lambda''+a)}
     {\ell! (a-\ell)! \Gamma(\lambda'+a-\ell) \Gamma(\lambda''+\ell)}
\frac{\partial^{a-\ell} f_1}{\partial z^{a-\ell}}
\frac{\partial^{\ell} f_2}{\partial z^\ell}
\end{equation}
satisfies \eqref{eqn:Rinv}, hence is a symmetry breaking operator
from $\pi_{\lambda'} \otimes \pi_{\lambda''}$ to $\pi_{\lambda'''}$.
\end{theorem}

\begin{remark}
{\rm{(1)}}\enspace
{\rm{
It turns out that
the coefficients appeared in the finite sum \eqref{eqn:RC} 
coincide with those of a Jacobi polynomial.
This fact can be checked if we know the formula \eqref{eqn:RC},
which is found, for example, by recurrence relations
that reflect the intertwining property \eqref{eqn:Rinv}.
But more intrinsically, 
the \lq\lq{F-method}\rq\rq, 
a method that the author and his collaborators introduced
\cite{xkhelgason, Eastwood60, KP1},
reveals directly why the Jacobi polynomial shows up,
see \cite[Sect.\ 9]{KP2} for example.}}
\par\noindent
{\rm{(2)}}\enspace
{\rm{Theorem \ref{thm:RC} does not assert the uniqueness of the 
Rankin--Cohen bidifferential operators. 
In fact, it turns out quite recently that
there exist symmetry breaking operators other than  
the Rankin--Cohen operator $R C_{\lambda', \lambda''}^{\lambda'''}$
for exceptional (negative) parameters $(\lambda', \lambda'', \lambda''')$,
and we gave the complete classification in \cite[Cor.\ 9.3]{KP2} (2015).
The main machinery of the proof is the \lq\lq{F-method}\rq\rq, 
which connects the following different topics:
\vskip 0.05in
\begin{itemize}
\item
the dimension of the polynomial solutions to the hypergeometric differential equation,
\vskip 0.05in
\item
the determination of the composition series of
the tensor product of two reducible Verma modules.
\end{itemize}
\smallskip

\noindent
The name \lq\lq{F-method}\rq\rq\ originated from the fact that it utilizes 
the \lq\lq{algebraic Fourier transform 
of Verma modules}\rq\rq\ \cite{xkhelgason, KP1}.
This method is also applied to the construction of differential symmetry breaking operators
for groups of higher rank \cite{FJS20, KKP16, KOSS, KP2}.
}}
\end{remark}

\subsection{
Classification theory of symmetry breaking operators in conformal geometry}
\label{subsec:10.3}
We end this article with yet another example for interactions of the theory of branching
laws with different fields of mathematics, this time with conformal geometry.

Consider the following problem:
given a Riemannian manifold $X$ and its submanifold $Y$,
find \lq\lq{conformally covariant}\rq\rq\ operators from the space of 
functions on $X$ to that on the submanifold $Y$.
We may also consider a generalization of this problem,
e.g.\ from \lq\lq{functions}\rq\rq\ to \lq\lq{differential forms}\rq\rq\
or from \lq\lq{Riemannian manifolds}\rq\rq\ to 
\lq\lq{pseudo-Riemannian manifolds}\rq\rq.
To formulate the problem rigorously,
we introduce the following notation.
\begin{align*}
G:&=\operatorname{Conf}(X): 
\hspace{12pt}
\text{the group of conformal transformations of $X$,}
\\
G':&= \operatorname{Conf}(X;Y):
\text{the subgroup of $G$ consisting of elements that preserve $Y$.}
\end{align*}
For any Riemannian manifold $X$, 
one can form a family of 
$\operatorname{Conf}(X)$-equivariant line bundles
${\mathcal{L}}_{\lambda}$ 
($\lambda \in {\mathbb{C}}$) 
over $X$, hence obtain a natural family of representations $\Pi_\lambda$
of the conformal group $G$ on the space $\Gamma(X,\mathcal{L}_{\lambda})$
of smooth sections for the line bundle $\mathcal{L}_{\lambda}$.
Since these line bundles are trivial topologically,
one may realize the family $\Pi_{\lambda}$ of representations of $G$
as multiplier representations
on the vector space 
$C^{\infty}(X)$
\cite[Sect.\ 2]{KO1}.
Since the subgroup $G'$ acts conformally 
on the submanifold $Y$ equipped with the induced metric,
a similar family $\pi_{\nu}$ ($\nu \in {\mathbb{C}}$) of representations of $G'$
can be defined on $C^{\infty}(Y)$ via the canonical group
homomorphism $\operatorname{Conf}(X;Y) \to \operatorname{Conf}(Y)$.
These representations are extended to 
the representations
$\Pi_{\lambda}^{(i)}$ and $\pi_{\nu}^{(j)}$ on
the spaces ${\mathcal{E}}^i(X)$ and ${\mathcal{E}}^j(Y)$
of differential forms, respectively.

\begin{basicproblem}
[Symmetry breaking operators in conformal geometry 
\cite{FJS20, KKP16, xksbonvec}]
\label{prob:confSBO}
Let 
$X$ be a Riemannian manifold
and $Y$ a submanifold of $X$.
For what parameters $(i,j,\lambda,\nu)$,
does a non-zero continuous operator $T \colon {\mathcal{E}}^i(X) \to {\mathcal{E}}^j(Y)$ 
satisfying
\[
  \pi_{\nu}^{(j)}(h) \circ T = T \circ \Pi_{\lambda}^{(i)}(h)
\quad
 {}^{\forall}h \in {\operatorname{Conf}}(X;Y) 
\]
exist? Further, find an explicit formula for such $T$.
\end{basicproblem}

For Basic Problem \ref{prob:confSBO},
if there exists a 
\lq\lq{universal construction}\rq\rq\
of such an operator $T$ for any pair $(X,Y)$ of Riemannian manifolds,
then the \lq\lq{solution}\rq\rq\
must survive in the model space 
$(X,Y)=(S^{n}, S^{n-1})$,
which has \lq\lq{large symmetry}\rq\rq\ in the sense that
the dimension of the group $\operatorname{Conf}(X;Y)$
attains its maximum among all pairs $(X,Y)$ of Riemannian
manifolds with $\dim Y= n-1(\geq 3)$.
Recently,
the construction and classification theory of 
conformally covariant
symmetry breaking operators 
for the case $(X,Y)=(S^n, S^{n-1})$ have been
rapidly developed as follows and completed in \cite{xksbonvec}.
\smallskip

\begin{enumerate}
\item[$\bullet$]
($i=j=0$; 
 $T$: differential operator)\enspace
Juhl 
constructed all conformally covariant, differential 
symmetry breaking operators $T$ in the flat model
by determining
the coefficients of $T$ 
using recurrence relations (Book \cite{Juhl}, 2009).
Afterwards, a short proof on the construction and classification of $T$ was given by a different approach (F-method)
(Kobayashi--{\O}rsted--Somberg--Sou{\v c}ek \cite{KOSS}).
\vskip 0.05in

\item[$\bullet$]
($i=j=0$; $T$: general)\enspace
In general, 
there is much more possibility of having symmetry breaking 
operators if we allow
integral operators or singular integral operators
other than just differential operators
\cite{Eastwood60}.
For the scalar-valued case $(i=j=0)$,
all the symmetry breaking operators (including integral operators
and singular integral operators) were classified with explicit construction
of the distribution kernels of operators $T$ by 
Kobayashi--Speh (Book \cite{xtsbon}, 2015).
\vskip 0.05in

\item[$\bullet$]
($i$, $j$: general ; 
$T$: differential operator)\enspace
By enhancing the F-method to the matrix-valued case,
the construction and classification of $T$ in the matrix-valued 
case ($i$, $j$: general) 
were shown by Kobayashi--Kubo--Pevzner
(Book \cite{KKP16}, 2016).
See also Fischmann--Juhl--Somberg (Book \cite{FJS20}, 2020) which 
uses the F-method, too.
\vskip 0.05in

\item[$\bullet$]
($i$, $j$: general; $T$: general)\enspace
The classification was completed by Kobayashi--Speh
(Book \cite{xksbonvec}, 2018).
For $j=i+1$ or $i-2$, all the symmetry breaking operators 
are differential operators \cite[Thm.\ 3.6]{xksbonvec}, 
whereas all the differential symmetry breaking operators
for $j=i$, or $i-1$ are obtained as the 
\lq\lq{residues}\rq\rq\ of integral symmetry breaking operators
with meromorphic parameter \cite{Kobayashi18}.

\end{enumerate}
\vskip 0.05in

In this way,
Basic Problem \ref{prob:confSBO} 
(construction and classification of symmetry breaking operators)
for the model space $(X,Y)=(S^n,S^{n-1})$
was completely solved in \cite{xksbonvec} 
based on the results of \cite{KKP16} and \cite{xtsbon}.
The three books \cite{KKP16, xtsbon, xksbonvec} over 600 pages long in total
develop a general idea for the construction and the classification
of symmetry breaking operators 
for a pair $G \supset G'$ of reductive Lie groups, and then 
apply the idea to
the particular case $(G,G') = (O(n+1,1),O(n,1))$, which 
is also important in conformal geometry, giving a proof 
of the construction and classification of those operators.
Functional equations of such operators yield some new results
to the (local) Gross--Prasad conjecture in number theory 
and its generalization to the non-tempered 
case (Kobayashi--Speh \cite{xksGP}).
It is not easy to explain the methods in \cite{KKP16, xtsbon, xksbonvec}
in a few lines, but we try to give its flavor here
at the last part of this article.

The pair $(X,Y)=(S^n,S^{n-1})$ is regarded as a pair $(G/P,G'/P')$ of 
real flag varieties for the pair
$(G,G')=(O(n+1,1),O(n,1))$ of real reductive Lie groups.
We have chosen this specific pair $(G,G')$
in the article \cite{xtsbon}, as the first test case of 
an explicit construction and complete classification
of (non-local) symmetry breaking operators between principal series representations.
The pair $(G,G')$ satisfies the sphericity condition (ii) in Theorem \ref{thm:BB},
hence the uniform boundedness of multiplicity in the branching laws 
is a priori guaranteed (Stage A).
In this sense, 
various results of \cite{KKP16, xtsbon, xksbonvec} 
can be interpreted 
as a step forward to branching problems in Stages B and C in the general
setting with uniformly bounded (in particular, finite)
multiplicity property. With this in mind, we explain
the ideas and methods for the proofs 
given in the three books in this general setting.
First, it follows from Theorem \ref{thm:pp} that 
$(G \times G')/\operatorname{diag}(G')$ is real spherical;
in particular, the number of $G'$-orbits 
on the real flag variety $(G\times G')/(P\times P')$ under the 
diagonal action is finite, where $P$ and $P'$
are minimal parabolic subgroups of $G$ and $G'$, respectively.
Second, the support of the distribution kernel for any 
symmetry breaking operator from a principal series 
representation of $G$ to that of the subgroup $G'$ is  a closed
$\operatorname{diag}(G')$-invariant subset of the
real flag variety
$(G \times G')/(P\times P')$.
Thus there are only finitely many possibilities of the support.
Third, we construct distribution kernels of symmetry 
breaking operators for each orbit 
and show meromorphic continuation and functional equations.
The residue of the meromorphic family of symmetry breaking operators
are symmetry breaking operators with 
distribution kernels of \emph{smaller support}.
Then we proceed by induction on the stratification
 for the closure relations 
in $\operatorname{diag}(G')\backslash (G\times G')/(P\times P') \simeq 
P'\backslash G/P$.
The first step of the induction is the symmetry breaking operators
that can be described as \lq\lq{differential operators}\rq\rq,
which correspond to the smallest orbit, namely, 
the unique closed orbit \cite[Lem.\ 2.3]{KP1}.
Symmetry breaking operators that can be described as differential operators
could appear not only as a \lq\lq{series}\rq\rq\ but also \lq\lq{sporadically}\rq\rq\
\cite{KP2, xksbonvec}.
The former one may be obtained as the \lq\lq{residues}\rq\rq\ of the 
meromorphic continuation of other operators such as integral symmetry
operators (e.g. \cite{Kobayashi18}), but the latter is more involved.
By the F-method, 
the construction of all such operators can be reduced 
to a problem of determining some polynomials
(\lq\lq{special polynomials}\rq\rq)
 that satisfy a system of differential equations.
The classification for differential symmetry breaking operators
is completed by solving such a system of 
equations \cite{KKP16, KOSS}.
The support of other symmetry breaking operators is strictly
larger than the closed orbit.
Thus an inductive argument 
with respect to the closure relation of 
$\operatorname{diag}(G')$-orbits on the real flag variety
$(G\times G')/(P\times P')$ exhausts all 
symmetry breaking operators, with the last one being
\lq\lq{regular symmetry breaking operators}\rq\rq\
that are obtained by the analytic continuation of 
integral symmetry breaking operators.
See \cite[Chap.\ 3, Sect.\ 3]{xksbonvec} for further details.

\vskip 1pc
{\bf{Acknowledgements:}} 
The author would like to express his sincere gratitude to
his collaborators of the various projects mentioned in this article.
He is also grateful to the referees for their careful reading of this article.

This work was partially supported by Grant-in-Aid
for Scientific Research (A) (18H03669), Japan Society for the Promotion of
Science.
\vskip 0.1in

{\bf{Added in translation, November, 2021}}.

\vskip 0.1in
\noindent
The author is grateful to Professor Toshihisa Kubo for translating
the article which originally appeared in Japanese in Sugaku \textbf{71}
(2019) from the Mathematical Society of Japan.

\bibliographystyle{amsplain}

\begin{thebibliography}{00}

\bibitem{ALTV20}
J. D. Adams, M. A. A. van Leeuwen, P. E. Trapa, D. A. Vogan, Jr,
Unitary Representations of Real Reductive Groups, 
Ast{\' e}risque
{\bf{417}}  (2020), {viii+188 pp.}

\bibitem{AG2009}
A.~Aizenbud, D.~Gourevitch,
\textit{Multiplicity one theorem
 for $(GL_{n+1}({\mathbb{R}}), GL_{n}({\mathbb{R}}))$},
Selecta Math. (N.S.) \textbf{15} (2009), pp.~271--294. 

\bibitem{AGRS2010}
A.~Aizenbud, D.~Gourevitch, S.~Rallis, G.~Schiffmann,
\textit{Multiplicity one theorems},
Ann.\ of Math.\ (2),
{\bf{172}} (2010), 
pp.~1407--1434.

\bibitem{Bn1996}
Y.~Benoist, 
\textit{Actions propres sur les espaces homog{\`e}nes r{\'e}ductifs}, 
Ann.\ of Math. (2),
{\bf{144}} 
(1996),
pp.~315--347.

\bibitem{BK2015}
Y.~Benoist, T.~Kobayashi, 
\textit{Temperedness of reductive homogeneous spaces}, 
J.\ Eur.\ Math.\ Soc.\ {\bf{17}} (2015), 
{pp.~3015--3036}.

\bibitem{BK2017}
Y.~Benoist, T.~Kobayashi, 
Tempered homogeneous spaces II,
In: Dynamics, Geometry, Number Theory: 
The Impact of Margulis on Modern Mathematics 
(eds.\ D.\ Fisher, D.\ Kleinbock, and G.\ Soifer),
The University of Chicago Press, 2022, to appear.\
Available also at arXiv:1706.10131.  

\bibitem{BK2021}
Y. Benoist, T. Kobayashi, 
\textit{Tempered homogeneous spaces III}, 
J.\ Lie Theory {\bf{31}} (2021), 
{pp.~833--869}.

\bibitem{BK(a20)}
Y. Benoist, T. Kobayashi, 
\textit{Tempered homogeneous spaces IV}, 
J.\ Inst.\ Math.\ Jussieu {\bf{22}} (2023), {pp.~2879--2906}.

\bibitem{BKO}
S.~Ben Sa{\"i}d, T.~Kobayashi, B.~\O rsted, 
\textit{Laguerre semigroup and Dunkl operators},
Compos.\ Math.\, {\bf{148}} (2012), 
{pp.~1265--1336}.

\bibitem{ber}
M. Berger, 
{\it Les espaces sym\'{e}triques non compacts}, 
Ann.\ Sci.\ \'{E}cole\ Norm.\ Sup.\ (3) {\bf 74} (1957), 
{pp.~85--177}.

\bibitem{BR}
J.~Bernstein, A.~Reznikov,
\textit{Analytic continuation of representations and estimates of automorphic 
forms},  
Ann.\ of Math.\ (2)
{\bf{150}} (1999), 
pp.~329--352.

\bibitem{xbrion}
M. Brion, 
\textit{Classification des espaces homog{\`e}nes sph{\'e}riques}, 
Compos. Math. {\bf{63}} (1987), pp.~189--208.

\bibitem{CM}
E. Calabi, L. Markus,
\textit{Relativistic space forms},
Ann.\ of Math. (2),
{\bf{75}} 
(1962), 
pp.~63--76,  

\bibitem{Clerc1617}
J.-L. Clerc,
\textit{Singular conformally invariant trilinear forms, I: 
The multiplicity one theorem; II: The higher multiplicity case}, 
Transform. Groups {\bf{21}}  (2016), {pp.~619--652};
ibid, {\bf{22}}  (2017), {pp.~651--706}.

\bibitem{CKOP}
J.-L.~Clerc, 
T.~Kobayashi,
B.~{\O}rsted, and M.~Pevzner, 
\textit{Generalized
Bernstein--Reznikov integrals},  
Math. Ann.,
{\bf{349}} 
(2011), 
{pp.~395--431}.

\bibitem{xcmz97}
P.~B. Cohen, Y.~Manin, D.~Zagier, 
\textit{Automorphic pseudodifferential operators}, 
In Memory of Irene Dorfman, (eds.\ A.\ S.\ Fokas and I.\ M.\ Gelfand),
Progr. Nonlinear Differential Equations Appl., 
\textbf{26}, 
Birkh\"auser, 1997, 
pp.~17--47.  
%

\bibitem{CDBL18}
D. Constales, H. De Bie, P. Lian,
\textit{Explicit formulas for the Dunkl dihedral kernel
and the $(k, a)$-generalized Fourier kernel}, 
J.\ Math.\ Anal.\ Appl.\ 
{\bf{460}}  (2018), {pp.~900--926}.


\bibitem{xdelorme}
P. Delorme,
\textit{Formule de Plancherel pour les espaces sym\'etriques r\'eductifs},
Ann.\ of Math. (2), {\textbf{147}} (1998),  pp.~417--452.
%
\bibitem{Duflo82}
M.~Duflo, 
Th{\'e}orie de Mackey pour les groupes
 de Lie alg{\'e}briques, 
Acta Math.~{\bf{149}} (1982), 
153--213.
%

\bibitem{DGV17}
M. Duflo, E. Galina, J. A. Vargas,
\textit{Square integrable representations of reductive 
Lie groups with admissible restriction to $SL_2(\mathbb{R})$}, 
J.\ Lie Theory
{\bf{27}}  (2017), {pp.~1033--1056}.

\bibitem{xdv}
M. Duflo, J. A. Vargas, 
\textit{Branching laws for square integrable representations},
Proc. Japan Acad. Ser. A,
Math. Sci., 
\textbf{86} 
(2010),
pp.~49--54.  

\bibitem{FJS20}
M. Fischmann, A. Juhl, P. Somberg,
Conformal Symmetry Breaking Differential Operators on Differential Forms,
Mem. Amer. Math. Soc. (2020), \textbf{268} 
{no. 1304}, v+112 pages.  

\bibitem{Frahm}
J.~Frahm, 
\textit{Symmetry breaking operators for strongly spherical reductive pairs},
Publ.\ Res.\ Inst.\ Math.\ Sci.\ {\bf{59}} (2023), pp.~259--337.

\bibitem{F2010}
H.~Fujiwara,
Unitary representations of exponential solvable Lie groups---the orbit method,
Sugaku-no-mori {\bf{1}}, Sugakushobo, 2010 (Japanese);
see also \cite{FL15}.

\bibitem{FL15}
H. Fujiwara, J. Ludwig,
Harmonic Analysis on Exponential Solvable Groups,
Springer Monogr.\ Math., Springer, 2015, xii+465 pages.  

\bibitem{GGPW12}
W.\ T.\ Gan, B.\ H.\ Gross, D.\ Prasad, J.-L.\ Waldspuruger,
Sur les conjectures de Gross et Prasad.\ I,
Ast{\'e}risque \textbf{346} (2012), 
Soci{\'e}t{\'e} Math{\'e}matique de France,
Paris, 2012. xi+318 pages.

\bibitem{GP} 
B.~Gross, D.~Prasad,  
\textit{On the decomposition of a representation of $SO_n$
 when restricted to $SO_{n-1}$}, 
 Canad. J. Math. {\bf{44}} 
 (1992), 
 pp.~974--1002.  

\bibitem{xgrwaII}
B. Gross, N. Wallach,
\textit{Restriction of small discrete series representations
        to symmetric subgroups},
        In: The Mathematical Legacy of Harish-Chandra,
        (eds, R.\ S.\ Doran and V.\ S.\ Varadarajan),
\textit{Proc. Sympos. Pure Math.},
       \textbf{68} (2000),  Amer. Math. Soc., pp.~255--272.

\bibitem{He20}
H. He,
\textit{A criterion for discrete branching laws for 
Klein four symmetric pairs and its applications to $E_{6(-14)}$}, 
Internat.\ J.\ Math. 
{\bf{31}}  (2020), 2050049, 15 pp.

\bibitem{Helgason}
S. Helgason, 
Groups and Geometric Analysis. 
Integral Geometry, Invariant Differential Operators, 
and Spherical functions. 
Math.\ Surveys Monogr., {\bf{83}} 
Amer.\ Math.\ Soc.\ 
2000. xxii+667 pp.

\bibitem{Helgason2}
S. Helgason, 
Differential Geometry, Lie Groups, and Symmetric Spaces.
Corrected reprint of the 1978 original,
Grad.\ Stud.\ Math., {\bf{34}} 
Amer.\ Math.\ Soc.\ 
2001. xxvi+641 pages.

\bibitem{HKMM11}
J. Hilgert, T. Kobayashi, G. Mano, J. M{\"o}llers,
\textit{Special functions associated with a certain
fourth-order differential equation}, 
Ramanujan J.\ 
{\bf{26}}  (2011), {pp.~1--34}.

\bibitem{xhkm}
J. Hilgert, T. Kobayashi, and J. M\"ollers. Minimal representations via 
Bessel operators. 
J. Math. Soc. Japan {\bf{66}} (2014), 
pp.\ 
{349--414}.

\bibitem{xhowe}
R. Howe,
\textit{$\theta$-series and invariant theory},
In: 
Automorphic Forms,
Representations, and $L$-functions, (eds.\ A.\ Borel and W.\ Casselman),
Proc. Symp. Pure Math.
\textbf{33}. Part I, 
(1979),
Amer. Math. Soc.,
pp.~275--285.

\bibitem{HW90}
A. T. Huckleberry, T. Wurzbacher,
\textit{Multiplicity-free complex manifolds}, 
Math.\ Ann.\ 
{\bf{286}}  (1990), {pp.~261--280}.


\bibitem{Juhl}
A. Juhl, 
Families of Conformally Covariant Differential Operators, $Q$-curvature and Holography. 
Progr. Math., 
{\bf{275}} 
Birkh\"auser, 
2009. 

\bibitem{KVcone}
M. Kashiwara, and M. Vergne, 
$K$-types and singular spectrum.
In: Noncommutative Harmonic Analysis,
(eds.\ J.\ Carmona and M.\ Vergne),
Lecture Notes in Math., {\bf{728}} Springer, Berlin, 1979,
pp. 177--200.

\bibitem{Adv16}
 F. Kassel, T. Kobayashi, 
\textit{Poincar\'e series for non-Riemannian locally 
symmetric spaces}, 
Adv. Math. {\bf{287}} (2016), 
{pp.~123--236}.

\bibitem{KK19}
F. Kassel, T. Kobayashi,
\textit{Invariant differential operators on 
spherical homogeneous spaces with overgroups}, 
J.\ Lie  Theory
{\bf{29}}  (2019), {pp.~663--754}.

\bibitem{KaKo20}
 F. Kassel, T. Kobayashi, 
\textit{Spectral analysis on pseudo-Riemannian locally symmetric spaces}, 
Proc.\ Japan Acad.\ Ser.\ A Math.\ Sci.\ 
{\bf{96}} (2020), 
{pp.~69--74}.
See also arXiv:1912.12601 for the full paper, to appear in 
Lecture Notes in Mathematics, Springer--Nature.

\bibitem{kirillov62}
A.~A.~Kirillov, 
{\it{Unitary representations of nilpotent Lie groups}}, 
Uspehi Mat.~Nauk {\bf{17}} 1962, 
57--110.

\bibitem{kirillov}
A.~A.~Kirillov, 
{\it Lectures on the Orbit Method}, 
Grad.\ Stud.\ Math., {\bf 64} 
Amer.~Math.~Soc.~2004. 

\bibitem{Kitagawa}
M.~Kitagawa
\textit{Algebraic structure
 on the space of intertwining operators}, 
Ph.D.\ dissertation, The University of Tokyo, 2016.  

\bibitem{xknappv}
A. W. Knapp, D. Vogan, Jr., 
Cohomological Induction and Unitary Representations,
Princeton Math.\ Ser.\ {\bf{45}},
Princeton Univ.\ Press, 1995.  
%

\bibitem{KZ}
A. W. Knapp and G. J. Zuckerman,
Classification of irreducible tempered representations of semisimple 
groups, 
Ann.\ of Math.\ (2) {\bf{116}} (1982),
389--455; 
II. ibid., 457--501.

\bibitem{kob89}
T. Kobayashi, 
\textit{Proper action on a homogeneous space of 
reductive type}, 
Math. Ann.~{\bf{285}} 
(1989), 
{pp.~249--263}.

\bibitem{Kobayashi92}
T. Kobayashi,
Singular Unitary Representations and Discrete Series for Indefinite Stiefel Manifolds $U(p,q;\mathbb{F})/U(p-m,q;\mathbb{F})$,
Mem. Amer. Math. Soc. (1992), \textbf{95} 
{no. 462}, v+106 pages.  

\bibitem{xk:1}
T. Kobayashi, 
\textit{The restriction of $A_{\mathfrak{q}}(\lambda)$
            to reductive subgroups},
           Proc. Japan Acad. Ser.\ A Math.\ Sci.,
           {\textbf{69}}
           (1993),
{pp.~262--267}.  

\bibitem{xkInvent94}
T.~Kobayashi, 
{\textit{Discrete decomposability of the restriction of
             $A_{\frak q}(\lambda)$
            with respect to reductive subgroups and its applications}}, 
Invent. Math.,
{\bf{117}} 
(1994), 
{pp.~181--205}.  
%
\bibitem{K94b}
T.~Kobayashi, 
\textit{Harmonic analysis on homogeneous manifolds of reductive type
         and unitary representation theory},
Sugaku \textbf{46} (1994), pp. 124--143 (Japanese);
In: Selected Papers on Harmonic Analysis, Groups, and Invariants,
(ed.\ K.\ Nomizu), Amer.\ Math.\ Soc.\ Transl.\ Ser.\ 2, 
\textbf{183}, Amer.\ Math.\ Soc., 
Providence, RI, 1998, {pp.\ 1--31} (English translation).
%

\bibitem{Ksuron}
T.~Kobayashi,
\textit{Introduction to harmonic analysis
 on real spherical homogeneous spaces},
In: Proceedings of the 3rd Summer School on Number Theory
\lq\lq{Homogeneous Spaces and Automorphic Forms}\rq\rq\
in Nagano, (ed.\ F.\ Sato), 1995, pp.~22--41.  
%

\bibitem{kob96}
T. Kobayashi, 
\textit{Criterion for proper actions on 
homogeneous spaces of reductive groups}, 
J. Lie Theory~{\bf{6}} (1996), 
{pp.~147--163}.
%

\bibitem{xkAnn98}
T.~Kobayashi, 
{\textit{Discrete decomposability of the restriction of
             $A_{\frak q}(\lambda)$
            with respect to reductive subgroups {\rm{II}}---micro-local analysis and asymptotic $K$-support}}, 
Ann.\ of Math. (2), 
{\bf {147}} 
(1998), 
{pp.~709--729}.  
%
\bibitem{xkInvent98}
T.~Kobayashi, 
{\textit{Discrete decomposability of the restriction of
             $A_{\frak q}(\lambda)$
            with respect to reductive subgroups {\rm{III}}---restriction of Harish-Chandra modules
 and associated varieties}}, 
Invent. Math., {\bf{131}} (1998), 
{pp.~229--256}.  
%
\bibitem{xkdisc}
   T. Kobayashi, 
\textit{Discrete series representations for the orbit spaces
     arising from two involutions of real reductive Lie groups},
     J. Funct. Anal.,
     \textbf{152} (1998), 
{pp.~100--135}.
%
\bibitem{K1998d}
T. Kobayashi, 
\textit{Deformation of compact Clifford--Klein forms of indefinite-Riemannian 
homogeneous manifolds}, 
 Math. Ann., {\bf{310}} (1998), 
{pp.~395--409}.
%
\bibitem{xk2001}
T. Kobayashi,
Discontinuous groups for non-Riemannian homogeneous spaces,
Mathematics Unlimited - 2001 and Beyond 
(B.~Engquist and W.~Schmid, eds.), 
Springer-Verlag, 2001, pp.~723--747. 
Expanded Japanese translation:
Discontinuous groups for non-Riemannian homogeneous spaces, In: Advanced
Mathematics, a challenge to the 21st century, Maruzen, 2002, pp.~18--73.
%
\bibitem{xk2004}
T.~Kobayashi, 
\textit{Geometry of multiplicity-free representations of $GL(n)$,
 visible actions on flag varieties, and triunity},
Acta Appl. Math. {\bf{81}} (2004), 
{pp.~129--146.}

\bibitem{deco-euro}
T. Kobayashi, 
\textit{Restrictions of unitary representations of real reductive groups},
Progr. Math. \textbf{229} 
{pp.~139--207},
Birkh\"auser, 2005.
%
\bibitem{xrims40}
T. Kobayashi, 
\textit{Multiplicity-free representations and visible actions
on complex manifolds}, 
Publ. Res. Inst. Math. Sci. 
{\textbf{41}} (2005),
{pp.~497--549},
special issue commemorating the fortieth anniversary of the founding of RIMS.

\bibitem{mf-korea}
T. Kobayashi, 
\textit{Multiplicity-free theorems of the restrictions of unitary 
highest weight modules with respect to reductive symmetric pairs}, 
Progr. Math., 
\textbf{255} 
{pp.~45--109}, 
Birkh\"auser, 2008. 
%
\bibitem{K2007a}
T. Kobayashi, 
\textit{Visible actions on symmetric spaces}, 
Transform. Groups, \textbf{12} (2007), 
{pp.~671--694}.  

\bibitem{K2007b}
T. Kobayashi,
\textit{A generalized Cartan decomposition
 for the double coset space 
$(U(n_1)\times U(n_2)\times U(n_3))\backslash U(n)/(U(p)\times U(q))$}, 
J. Math. Soc. Japan \textbf{59} (2007), 
{pp.~669--691}.  

\bibitem{kob09}
T. Kobayashi, 
\textit{Hidden symmetries and spectrum of the 
Laplacian on an indefinite Riemannian manifold}, 
In: Spectral Analysis in Geometry and Number Theory
---In Honor of T.\ Sunada (eds.\ M.\ Kotani, H.\ Naito, and T.\ Tate),
Contemp. Math.~{\bf{484}} 
{pp.~73--87}, 
Amer. Math. Soc., 
2009.

\bibitem{sato50}
T. Kobayashi, 
\textit{Algebraic analysis of minimal representations}, Publ. Res.\ Inst.\ Math.\ Sci.\ 
{\bf{47}} (2011), 
pp.
{585--611}, 
Special issue in commemoration of the golden jubilee of algebraic 
analysis.

\bibitem{Zuckerman60}
T. Kobayashi, 
\textit{Branching problems of Zuckerman derived functor modules}, 
In: Representation Theory and Mathematical Physics
---In Honor of G.~Zuckerman,
(eds.\ J.\ Adams, B.\ Lian, and S.\ Sahi), 
Contemp. Math., 
{\bf{557}} 
{pp.~23--40}, 
Amer. Math. Soc., Providence, RI, 2011.

\bibitem{K12} 
T.~Kobayashi, 
\textit{Restrictions of generalized Verma modules to symmetric 
pairs}, 
Transform. Groups, 
{\bf{17}} 
 (2012), 
{pp.~523--546}.
%
\bibitem{xkhelgason}
T.~Kobayashi, 
\textit{$F$-method for constructing equivariant 
 differential operators},
In: Geometric Analysis and 
Integral Geometry---AMS Special Session on Radon Transforms 
and Geometric Analysis in Honor of Sigurdur Helgason's 85th 
Birthday, January 4--7, 2012 Boston, MA,
(eds.\ E.\ T.\ Quinto, F.\ Gonzalez and J.\ G.\ Christensen),
 Contemp. Math., {\bf {598}},
Amer. Math. Soc., Providence, RI,
2013,
{pp.~139--146}.

\bibitem{mfbundl}
T.~Kobayashi, 
\textit{Propagation of multiplicity-freeness property for
holomorphic vector bundles}, 
In: Lie Groups: Structure, Actions, 
and Representations---In Honor of Joseph A.\ Wolf on the 
Occasion of his 75th Birthday, Bochum, 2012,
(eds.\ A.\ Huckleberry, I.\ Penkov and G.\ Zuckerman),
Progr. Math., {\bf{306}},
Birkh\"{a}user, 
2013, 
{pp.~113--140}. 

\bibitem{Eastwood60}
T. Kobayashi, 
\textit{F-method for symmetry breaking operators},
 Differential Geom. Appl. 
{\bf{33}}, suppl.\ 
(2014), 
{pp.~272--289},
Special issue in honor of M.~Eastwood.  

\bibitem{xkProg2014}
T. Kobayashi,
\textit{Shintani functions, real spherical manifolds, and symmetry breaking operators}, 
Dev.\ Math., 
{\bf{37}} 
(2014), 
{pp.~127--159}.  
Springer.

\bibitem{frenkel60}
T.~Kobayashi,
\textit{Special functions in minimal representations},
In: Perspectives in Representation Theory in honor of Igor Frenkel
 on his 60th birthday, 
Comtemp. Math., {\bf{610}} 
{pp.~253--266}.
Amer. Math. Soc.,  2014.

\bibitem{xKVogan2015}
T.~Kobayashi, 
\textit{A program for branching problems in the representation 
theory of real reductive groups}, 
In: Representations of Reductive 
Groups---In Honor of the 60th Birthday of David A.\ Vogan,
Jr., (eds.\ M.\ Nevins and P.\ E.\ Trapa),
Progr. Math., 
{\bf{312}}, 
Birkh{\"a}user, 
2015,
{pp.~277--322}.

\bibitem{xkIntrinsic}
T. Kobayashi, 
\textit{Intrinsic sound of anti-de Sitter manifolds}, 
Springer Proc.~Math.~Stat., 
{\bf{191}}, Springer,
2016, 
pp.~83--99.  

\bibitem{xkmsj70}
T. Kobayashi, 
Birth of New Branching Problems. 70th anniversary lecture (featured invited talk), 
MSJ Autumn Meeting 2016. Kansai University, Japan, 15-18 September 2016
(Japanese).

\bibitem{K16b}
T. Kobayashi, 
\textit{Global analysis by hidden symmetry}, 
In: Representation Theory, Number Theory, 
and Invariant Theory---In Honor of Roger Howe
on the Occasion of His 70th Birthday, 
New Haven, CT, 2015,
(eds.\ J.\ Cogdell,
J.-L.\ Kim and C.-B.\ Zhu),
Progr.\ Math.,
{\bf{323}},
Birkh{\"a}user,
2017, 
pp.\ 359--397.

\bibitem{Kobayashi18}
T. Kobayashi
\textit{Residue formula for regular symmetry branching operators}, 
Contemp.\ Math.\ 
{\bf{714}}  (2018), {pp.~175--197}, Amer.\ Math.\ Soc.

\bibitem{Kobayashi(a19)}
 T. Kobayashi, 
\textit{Branching laws of unitary representations associated to minimal
elliptic orbits for indefinite orthogonal group $O(p,q)$}, 
Adv.\ Math.\ 
{\bf{388}}  (2021), Paper No.\ 107862, 38 pages.

\bibitem{Kobayashi21}
T. Kobayashi,
\textit{Bounded multiplicity theorems for induction and 
restriction}, 
J.\ Lie Theory
{\bf{32}} 
 (2022), 
pp.\ 197--236.
Available also at arXiv:2109.14424.

\bibitem{Kobayashi20}
 T. Kobayashi, 
\textit{Admissible restrictions of irreducible representations of 
reductive Lie groups: symplectic geometry and discrete decomposability}, 
To appear in the special issue in memory of Bertram Kostant,
Pure Appl.\ Math.\ Q.\
Available also at arXiv:1907.12964.

\bibitem{KKP16}
T. Kobayashi, T. Kubo, M. Pevzner, 
Conformal Symmetry Breaking Operators
 for Differential Forms on Spheres, 
Lecture Notes in Math., 
{{\bf{2170}}} 
Springer, 
2016, 
viii $+$ 192 pages.  

\bibitem{xkmanoAMS}
T.~Kobayashi, G.~Mano,
The Schr{\"o}dinger Model for the Minimal
Representation of the Indefinite Orthogonal Group $O(p,q)$,
Mem. Amer. Math. Soc. (2011), \textbf{212} 
{no. 1000}, vi+132 pages.  

\bibitem{xKMt}
T. Kobayashi, T. Matsuki, 
\textit{Classification of finite-multiplicity symmetric pairs},
Transform. Groups, 
 {\bf{19}} 
(2014), 
{pp.~457--493}, 
Special issue in honor of Dynkin
 for his 90th birthday. 

\bibitem{KN03}
T. Kobayashi, S. Nasrin,
\textit{Multiplicity one theorem in the orbit method}, 
Lie Groups and Symmetric Spaces: In memory of 
Professor F.\ I.\ Karpelevi{\v c} 
(ed. S.\ Gindikin), Amer.\ Math.\ Soc.\ Transl.\ Ser. 2
{\bf{210}}  (2003), {pp.~161--169}.

\bibitem{xKN}
T. Kobayashi, S. Nasrin, 
\textit{Geometry of coadjoint orbits and multiplicity-one branching laws
 for symmetric pairs}, 
Algebr.\ Represent.\ Theory 
{\bf{21}} 
 (2018), 
pp.\ 1023--1036, 
Special issue in honor of Alexandre Kirillov.  

\bibitem{kks}
T.\ Kobayashi, K.\ Ono, T.\ Sunada, 
\textit{Periodic Schr{\"o}dinger operators on a manifold}, 
Forum Math. 
{{\bf{1}}} 
(1989), 
pp.~69--79.  

\bibitem{KO1}
T.~Kobayashi and B.~{\O}rsted,
\textit{
{Analysis on the
  minimal representation
  of\/ {${\rm O}(p,q)$}. {\rm{I}}. Realization via conformal
  geometry}},
Adv. Math. \textbf{180} (2003), 486--512.


\bibitem{KO2}
T.~Kobayashi and B.~{\O}rsted,
\textit{
{Analysis
  on the minimal representation
 of\/ {${\rm O}(p,q)$}. {\rm{II}}. Branching laws}},
  Adv. Math. \textbf{180} (2003), 513--550.

\bibitem{KO3}
T.~Kobayashi and B.~{\O}rsted,
\textit{
{Analysis
  on the minimal representation
 of\/ {${\rm O}(p,q)$}. {\rm{III}}. Ultrahyperbolic
  equations on {$\mathbb{R}\sp {p-1,q-1}$}}},
Adv. Math. \textbf{180} (2003),
  551--595.

\bibitem{KOSS}
T.~Kobayashi, B.~{\O}rsted,
P.~Somberg, V.~Sou\v{c}ek, 
\textit{Branching laws for Verma modules and applications in
parabolic geometry}, I, 
Adv. Math., {\bf{285}} 
(2015), 
{pp.~1796--1852}. 

\bibitem{xktoshima}
T.~Kobayashi, T.~Oshima, 
\textit{Finite multiplicity theorems for induction and restriction}, 
Adv. Math., \textbf{248} (2013), 
{pp.~921--944}. 

\bibitem{decoAq}
T. Kobayashi, Y. Oshima, 
\textit{Classification of discretely decomposable 
$A_{\mathfrak {q}}(\lambda)$
 with respect to reductive symmetric pairs}, 
Adv. Math., 
{\bf{231}}
 (2012), 
{pp.~2013--2047}.
%
\bibitem{xkyosh13}
T. Kobayashi, Y. Oshima, 
\textit{Classification of symmetric pairs with 
discretely decomposable restrictions of $({\mathfrak{g}},K)$-modules}, 
J.\ Reine Angew.\ Math.,  
{\bf{2015}} 
(2015), 
no.703, 
{pp.~201--223}.  

\bibitem{KP1}
T.~Kobayashi, M.~Pevzner,
\textit{Differential symmetry breaking operators}. 
I. \textit{General theory and F-method}, 
Selecta Math.
 (N.S.) {\bf{22}} (2016), 
{pp.~801--845}.

\bibitem{KP2}
T.~Kobayashi, M.~Pevzner,
\textit{Differential symmetry breaking operators}. 
II. \textit{Rankin--Cohen operators for symmetric pairs}, 
Selecta Math.
(N.S.) {\bf{22}} (2016), pp.\ 847--911.

\bibitem{KP3}
T.~Kobayashi, M.~Pevzner,
\textit{Inversion of Rankin--Cohen operators
 via holographic transform}, 
Ann.\ Inst.\ Fourier (Grenoble) 
{\bf{70}} (2020), pp.\ 2131--2190.
  

\bibitem{xtsbon}
T. Kobayashi, B. Speh,
Symmetry Breaking
 for Representations of Rank One
 Orthogonal Groups, 
(2015), 
 Mem.\ Amer.\ Math.\ Soc. 
 {\bf{238}} 
{no.1126}, 
 118 pages.  

\bibitem{xksGP}
T.~Kobayashi, B.~Speh, 
\textit{Symmetry breaking for orthogonal groups and a conjecture
 by B.~Gross and D.~Prasad}. 
In: Geometric Aspects of the Trace Formula,
 (eds.\ W.\ M{\"u}ller, S.\ W.\ Shin and N.\ Templier),
 Simons Symp., 
 Springer, 
2018,
pp.\ 245--266.  

\bibitem{xksbonvec}
T.~Kobayashi, B.~Speh, 
Symmetry Breaking for Representations of Rank One Orthogonal Groups, 
II, 
Lecture Notes in Math., {\bf{2234}} 
Springer, 2018.  
xv$+$342 pages.  
%
\bibitem{xkramer}
M.~Kr{\"a}mer, 
\textit{Multiplicity free subgroups
 of compact connected Lie groups},
 Arch. Math. (Basel)
 {\bf{27}} (1976), pp.~28--36.  
%

\bibitem{MSV19}
O. M{\'a}rquez, S. Simondi, J. A. Vargas,
\textit{Branching laws: some results and new examples}, 
Rev.\ Un.\ Mat.\ Argentina
{\bf{60}}  (2019), {pp.~45--59}.

\bibitem{Mautner50}
F. I. Mautner,
\textit{Unitary representations of locally compact groups I, II}, 
Ann.\ of Math. (2)
{\bf{51}} (1950), {pp.~1--25}; 
{\bf{52}} (1950), {pp.~528--556}.

\bibitem{xmik}
I. V. Mikityuk,
\textit{Integrability of invariant Hamiltonian systems with homogeneous 
configuration spaces}, 
Math. USSR-Sbornik {\bf{57}} (1987),
pp.~527--546.  

\bibitem{Mol} 
V.F. 
Mol{\v c}hanov, Tensor products of unitary representations of the three-dimensional Lorentz group.
\emph{Math. USSR. Izv.} {\bf{15}} (1980), pp. 113--143.

\bibitem{xspeh}
B. \O rsted, B. Speh,
\textit{Branching laws for some unitary representations
 of $SL(4,{\mathbb{R}})$},
 SIGMA Symmetry Integrability Geom.\ Methods Appl.\ 
 {\bf{4}} (2008), 017.
%
\bibitem{Oshima2002}
T.~Oshima,
\textit{Harmonic analysis on semisimple symmetric spaces}, 
Sugaku {\bf{37}} (1985) pp.\ 97--112 (Japanese);
Sugaku Expositions, 
{\bf{15}} (2002), 
pp.~151--170, 
Amer. Math. Soc.
(English translation).

\bibitem{yophd}
Y.~Oshima,
\textit{Discrete branching laws of Zuckerman's derived functor modules},
Ph.D.\ dissertation, The University of Tokyo, 2013. 

\bibitem{Paradan15}
P.-E. Paradan,
\textit{Quantization commutes with reduction in the non-compact setting:
the case of holomorphic discrete series}, 
J.\ Eur.\ Math.\ Soc.\ (JEMS)
{\bf{17}}  (2015), {pp.~955--990}.


\bibitem{PT02}
F. Podest{\`a}, G. Thorbergsson,
\textit{Polar and coisotropic actions on K{\"a}hler manifolds}, 
Trans.\ Amer.\ Math.\ Soc.\ 
{\bf{354}}  (2002), {pp.~1759--1781}.

\bibitem{xrankin}
R.~A.~Rankin,
\textit{The construction of automorphic forms from the derivatives
 of a given form},
J. Indian Math. Soc. (N.S.),  
{\bf{20}} 
(1956), 
pp.~103--116.  
%
\bibitem{Re79}
J. Repka, 
\textit{Tensor products of holomorphic discrete series representations}, 
Canad. J. Math. {\textbf{31}} (1979),
{pp.~836--844}.

\bibitem{SaVe17}
Y.~Sakellaridis, A.~Venkatesh,
{Periods and {H}armonic {A}nalysis on {S}pherical {V}arieties},
{Ast{\' e}risque}, {\bf{396}}, Soc.\ Math.\ France, Paris, 2017.

\bibitem{Sa2008}
A.~Sasaki,
\textit{Multiplicity-free representations and visible actions},
Sugaku, {\bf{74}} (2022), no.\ 3, {pp.~225--252}; English transl.\
to appear in Sugaku Exposition, Amer.\ Math.\ Soc.

\bibitem{xschthe}
W. Schmid,
\textit{Die Randwerte holomorphe Funktionen
 auf hermitesch symmetrischen R{\" a}umen},   
Invent. Math.
{\textbf{9}} 
(1969--70),
pp.~61--80.  

\bibitem{xhsekiijm}
H.~Sekiguchi, 
\textit{Branching rules of singular unitary representations
 with respect to symmetric pairs $(A_{2n-1}, D_n)$},
Internat.\ J.\ Math.\ {\bf{24}} (2013), no. 4, 1350011, 25 pp. 

\bibitem{Stem2001}
J.~R.~Stembridge,
\textit{Multiplicity-free products of Schur functions}, 
Ann. Comb., {\bf{5}} 
(2001), 
pp.~113--121.  
%
\bibitem{Sugiura}
M.~Sugiura,
Introduction to Unitary Representations,
with Commentary by T.\ Kobayashi,
TokyoTosho Co.,Ltd. 2018,
viii+262 pages.

\bibitem{xsunzhu}
B.~Sun, C.-B.~Zhu,
\textit{Multiplicity one theorems:
the Archimedean case},
Ann.\ of Math. (2), 
{\bf{175}} 
(2012), 
pp.~23--44.  
%
%
\bibitem{Ty2015}
Y.~Tanaka,
\textit{Visible actions of reductive algebraic groups 
 on complex algebraic varieties}, 
Ph.D.\ dissertation, The University of Tokyo, 2015.
%

\bibitem{Ty2}
Y. Tanaka,
\textit{Visible actions of compact Lie groups
on complex spherical varieties}, 
J.\ Diff.\ Geom.\
to appear.

\bibitem{xtatsuuma}
N.~Tatsuuma
A Duality Theorem for Topological Groups,
Kinokuniya Sugaku Sosho {\bf{32}}
Kinokuniya Co. Ltd., 1994. (Japanese)
%

\bibitem{Tauchi(a21)}
T. Tauchi,
\textit{A generalization of the Kobayashi--Oshima uniformly 
bounded multiplicity theorem}, 
to appear in Internat.\ J.\ Math.\ 
(Ph.D.\ dissertation, The University of Tokyo, 2019).
Available also at arXiv:2108.02139.


\bibitem{Trapa01}
P. E. Trapa,
\textit{Annihilators and associated varieties of 
$A_{\mathfrak{q}}(\lambda)$ modules for $U(p,q)$}, 
Compositio Math.\ 
{\bf{129}}  (2001), {pp.~1--45}.

\bibitem{Vinberg01}
E. B. Vinberg,
\textit{Commutative homogeneous spaces and co-isotropic symplectic actions}, 
Russian Math.\ Surveys
{\bf{56}}  (2001), {pp.~1--60}.

\bibitem{Vogan81}
D. A. Vogan, Jr., 
Representations of Real Reductive Lie Groups, 
Progr. Math.\
{\bf{15}} 
Birkh{\"a}user, 
1981. 

\bibitem{Vogan84}
D. A. Vogan, Jr., 
\textit{Unitarizability of certain series of representations}, 
Ann.\ of Math. (2), {\bf{120}} (1984), pp.~141--187.  

\bibitem{WaI}
N. R. Wallach,
Real reductive groups. I, II, 
Pure Appl.\ Math.\ 
{\bf{132}} 
 Academic Press, Inc., Boston, MA, 1988;
{\bf{132}}-II, ibid, 1992. 
%

\bibitem{Weil64}
A. Weil,
\textit{Remarks on the cohomology of groups}, 
Ann.\ of Math. (2), {\bf{80}} (1964), pp.~149--157.  
%

\bibitem{Wolf07}
J. A. Wolf,
Harmonic Analysis on Commutative Spaces,
Math.\ Surveys Monogr.\
{\bf{142}} 
Amer.\ Math.\ Soc.\ 2007, xvi+387.


\bibitem{Wolpert}
S. A. Wolpert,
\textit{Disappearance of cusp forms
 in special families}, 
 Ann.\ of Math. (2),
 {\bf{139}} 
 (1994), 
pp.~239--291.  

\bibitem{Wong}
H.-W. Wong,
\textit{Dolbeault cohomological realization of Zuckerman modules associated with 
finite rank representations.}
J.\ Funct.\ Anal.\ {\bf{129}} (1995), no. 2, 428--454. 

\bibitem{ZL10}
F. Zhu, K. Lian,
\textit{On a branching law of unitary representations and 
a conjecture of Kobayashi}, 
C.\ R.\ Acad.\ Sci.\ Paris
{\bf{348}}  (2010), {pp.~959--962}.

\end{thebibliography}

\end{document}